\documentclass{article}
\usepackage{amssymb,amsfonts,amsmath,amsthm,amsopn,amstext,amscd,latexsym,xy,color,mathrsfs,etoolbox}
\usepackage{hyperref}

\usepackage{stmaryrd}

\input xy
\xyoption{all}

\newtheorem{theorem}{Theorem}[section]
\newtheorem*{theorem*}{Theorem}
\newtheorem{lemma}[theorem]{Lemma}

\newtheorem{proposition}[theorem]{Proposition}
\newtheorem{propdef}[theorem]{Proposition-Definition}
\newtheorem{corollary}[theorem]{Corollary}
\newtheorem{definition}[theorem]{Definition}

\theoremstyle{remark}

\newtheorem*{remark}{Remark}
\numberwithin{equation}{section}

\renewcommand{\(}{\textup{(}}
\renewcommand{\)}{\textup{)}}

\DeclareMathOperator{\Mod}{Mod}
\newcommand{\Modsm}{\Mod_\text{sm}}

\DeclareMathOperator{\GL}{GL}

\DeclareMathOperator{\Supp}{Supp}

\DeclareMathOperator{\Frob}{Frob}

\DeclareMathOperator{\Hom}{Hom}

\DeclareMathOperator{\Sp}{Sp}
\DeclareMathOperator{\Gal}{Gal}

\DeclareMathOperator{\Stab}{Stab}

\DeclareMathOperator{\diag}{diag}

\DeclareMathOperator{\ind}{ind}
\DeclareMathOperator{\End}{End}
\DeclareMathOperator{\tr}{tr}

\DeclareMathOperator{\Lie}{Lie}

\DeclareMathOperator{\Res}{Res}

\DeclareMathOperator{\Ind}{Ind}
\DeclareMathOperator{\Spa}{Spa}

\DeclareMathOperator{\Sh}{Sh}
\DeclareMathOperator{\Art}{Art}
\DeclareMathOperator{\Inf}{Inf}

\newcommand{\cA}{{\mathcal A}}
\newcommand{\cB}{{\mathcal B}}
\newcommand{\cC}{{\mathcal C}}

\newcommand{\cF}{{\mathcal F}}
\newcommand{\cG}{{\mathcal G}}
\newcommand{\cH}{{\mathcal H}}
\newcommand{\cI}{{\mathcal I}}
\newcommand{\cJ}{{\mathcal J}}

\newcommand{\cN}{{\mathcal N}}
\newcommand{\cO}{{\mathcal O}}

\newcommand{\cS}{{\mathcal S}}

\newcommand{\cV}{{\mathcal V}}

\newcommand{\ffrm}{{\mathfrak m}}

\newcommand{\bbA}{{\mathbb A}}

\newcommand{\bbC}{{\mathbb C}}

\newcommand{\bbF}{{\mathbb F}}

\newcommand{\bbL}{{\mathbb L}}

\newcommand{\bbP}{{\mathbb P}}
\newcommand{\bbQ}{{\mathbb Q}}
\newcommand{\bbR}{{\mathbb R}}

\newcommand{\bbT}{{\mathbb T}}

\newcommand{\bbV}{{\mathbb V}}

\newcommand{\bbZ}{{\mathbb Z}}

\newcommand{\plim}{\mathop{\varprojlim}\limits}

\newcommand{\wt}{{\widetilde{\tau}}}

\newcommand{\Fp}{{F^+}}

\title{Torsion Galois representations over CM fields and Hecke algebras in the derived category}
\author{James Newton and Jack A. Thorne}

\begin{document}
\setcounter{tocdepth}{2}
\maketitle

\tableofcontents

\section{Introduction}

In this paper we study the relation between Galois representations and the cohomology of arithmetic locally symmetric spaces. Let $F$ be a number field and let $n \geq 2$ be an integer. Associated to any open compact subgroup $U \subset \GL_n(\bbA_F^\infty)$ is the topological space defined as a double quotient
\[ X^U_{\GL_n} = \GL_n(F) \backslash \GL_n(\bbA_F) / U \times \bbR^\times U_\infty, \]
where $U_\infty$ is a fixed choice of maximal compact subgroup of $\GL_n(F \otimes_\bbQ \bbR)$. If $U$ is neat (a condition that can always be achieved by replacing $U$ by a finite index subgroup), then $X_{\GL_n}^U$ is naturally an orientable smooth manifold, and we now assume this. If $F = \bbQ$ and $n = 2$, then $X_{\GL_n}^U$ can be identified with the set of complex points of a classical modular curve. In general, however, the space $X_{\GL_n}^U$ has no direct link to algebraic geometry.

Nevertheless, several mathematicians (see e.g.\ \cite{Ash02}) have conjectured an explicit relation between the cohomology of the spaces $X_{\GL_n}^U$ and the representations of the absolute Galois group $G_F = \Gal(\overline{F}/F)$. A remarkable feature of this conjectured correspondence is that it should take into account torsion in the cohomology groups $H^\ast(X_{\GL_n}^U, \bbZ)$, which falls outside the scope of the theory of automorphic forms and, for example, earlier conjectures of Langlands and Clozel (see e.g.\ \cite{Clo90}).

Let us now assume that $F$ is an imaginary CM field (for example, an imaginary quadratic field). In a recent breakthrough work \cite{Sch14}, Scholze has established this torsion correspondence, in a form that we now describe. We first introduce some helpful notation.  It is enough to work `one prime at a time', so we fix a prime $p$. We suppose that our choice of level subgroup $U$ splits as a product $U = \prod_v U_v$ over the finite places $v$ of $F$, where each $U_v$ is an open compact subgroup of $\GL_n(\cO_{F_v})$. We let $S$ be a finite set of finite places of $F$, containing all the places dividing $p$, such that for all $v \not\in S$, we have $U_v = \GL_n(\cO_{F_v})$. 

We also introduce coefficients. Let $E$ be a finite extension of $\bbQ_p$ large enough to contain all embeddings of $F$ in $\overline{\bbQ}_p$, and let $\cO$ be its ring of integers, $k$ its residue field. We can associate to any tuple $\boldsymbol{\lambda} = (\lambda_\tau) \in (\bbZ^n)^{\Hom(F, E)}$ satisfying the condition 
\[ \lambda_{\tau, 1} \geq \lambda_{\tau, 2} \geq \dots \geq \lambda_{\tau, n} \]
for each $\tau \in \Hom(F, E)$ a local system $M_{\boldsymbol{\lambda}}$ of finite free $\cO$-modules on $X_{\GL_n}^U$. (The precise definition is given in \S \ref{sec_hecke_algebras} below, in terms of the algebraic representations of $\GL_n$ associated to the dominant weights $\lambda_\tau$. In the body of the paper, $M_{\boldsymbol{\lambda}}$ is denoted by the symbol $\underline{A(\GL_n; \boldsymbol{\lambda})}_{\GL_n}^U$ in order to keep track of its relation to other objects.) Then the cohomology groups
\[ H^\ast(X_{\GL_n}^U, M_{\boldsymbol{\lambda}}) \]
are finite $\cO$-modules, and for each finite place $v\not\in S$ of $F$ we can define a family of Hecke operators $T_v^1, \dots, T_v^n$ in terms of double cosets. We write $\bbT^S(H^\ast(X_{\GL_n}^U, M_{\boldsymbol{\lambda}}))$ for the (commutative) $\cO$-subalgebra of $\End_\cO(H^\ast(X_{\GL_n}^U, M_{\boldsymbol{\lambda}}))$ generated by these operators. We can now state one consequence of Scholze's results as follows (\cite[Theorem 5.4.1]{Sch14}):
\begin{theorem}\label{thm_intro_scholze_theorem}
There exists an integer $N = N(d, n)$ depending only on $n$ and $d = [F : \bbQ]$, an ideal $I \subset \bbT^S(H^\ast(X_{\GL_n}^U, M_{\boldsymbol{\lambda}}))$ satisfying $I^N = 0$, and a continuous group determinant
\[ D : G_{F, S} \to \bbT^S(H^\ast(X_{\GL_n}^U, M_{\boldsymbol{\lambda}}))/I \]
such that for each finite place $v \not\in S$ of $F$, the characteristic polynomial of $D(\Frob_v)$ is 
\begin{equation}\label{eqn_intro_characteristic_polynomial_of_frobenius}
 X^n - T_v^1 X^{n-1} + \dots + (-1)^j q_v^{j(j-1)/2} T_v^j X^{n-j} + \dots + (-1)^n q_v^{n(n-1)/2} T_v^n
\end{equation}
mod $I$.
\end{theorem}
Since group determinants are in bijective correspondence with isomorphism classes of semi-simple representations over algebraically closed fields, we deduce:
\begin{corollary}\label{cor_intro_scholze_result} \begin{enumerate}
\item Let $\phi \in H^\ast(X_{\GL_n}^U, M_{\boldsymbol{\lambda}}) \otimes_\cO \overline{\bbQ}_p$ be an eigenvector for $\bbT^S(H^\ast(X_{\GL_n}^U, M_{\boldsymbol{\lambda}}))$, in the sense that for all $T_v^i$, we have $T_v^i \phi = a_v^i \phi$ for some numbers $a_v^i \in \overline{\bbQ}_p$. Then there exists a continuous representation $\rho_\phi : G_{F, S} \to \GL_n(\overline{\bbQ}_p)$ such that for each finite place $v\not\in S$ of $F$, the characteristic polynomial of $\rho_\phi(\Frob_v)$ is $\sum_j (-1)^j q_v^{j(j-1)/2} a_v^j X^{n-j}$.
\item Let $\phi \in H^\ast(X_{\GL_n}^U, M_{\boldsymbol{\lambda}}) \otimes_\cO \overline{\bbF}_p$ be an eigenvector for $\bbT^S(H^\ast(X_{\GL_n}^U, M_{\boldsymbol{\lambda}}))$, in the sense that for all $T_v^i$, we have $T_v^i \phi = a_v^i \phi$ for some numbers $a_v^i \in \overline{\bbF}_p$. Then there exists a continuous representation $\overline{\rho}_\phi : G_{F, S} \to \GL_n(\overline{\bbF}_p)$ such that for each finite place $v\not\in S$ of $F$, the characteristic polynomial of $\overline{\rho}_\phi(\Frob_v)$ is $\sum_j (-1)^j q_v^{j(j-1)/2} a_v^j X^{n-j}$.
\end{enumerate}
\end{corollary}
The aim of this paper is to improve Theorem \ref{thm_intro_scholze_theorem} in ways that will be useful for applications to modularity of Galois representations, following the schema outlined by Calegari--Geraghty \cite{Cal14}. The first goal is to try to remove the nilpotent ideal $I$; indeed, it seems natural to expect that one should always have $I = 0$. The second goal is to replace the Hecke algebra $\bbT^S(H^\ast(X_{\GL_n}^U, M_{\boldsymbol{\lambda}}))$ by a derived variant that has $\bbT^S(H^\ast(X_{\GL_n}^U, M_{\boldsymbol{\lambda}}))$ as a quotient, but a priori could be larger.

Let us now discuss these goals in more detail. We first choose a maximal ideal 
\[ \ffrm \subset \bbT^S(H^\ast(X_{\GL_n}^U, M_{\boldsymbol{\lambda}})) \]
such that the associated Galois representation $\overline{\rho}_\ffrm$ (which exists by Corollary \ref{cor_intro_scholze_result}) is absolutely irreducible. (In the body of the paper, we refer to such an ideal as a \emph{non-Eisenstein} maximal ideal.) We will work after localization at $\ffrm$. Since one of our main motivations is the possibility of applying our results in the context of $R = \bbT$ theorems, this seems like a natural simplifying step. 

Now we define our derived Hecke algebra.\footnote{We find it convenient in this paper to use the terminology `derived Hecke algebra', which refers to an enhancement of the usual notion of Hecke algebra living in the derived category. However, we wish to emphasize that this is not the same as the derived Hecke algebra considered in recent works of Venkatesh, in which additional `derived' Hecke operators are considered which act on cohomology by shifting degrees. It is clear that there is a common generalization of these two notions, but we do not discuss this here.} We replace the groups $H^\ast(X_{\GL_n}^U, M_{\boldsymbol{\lambda}})$ by the complex $R \Gamma(X_{\GL_n}^U, M_{\boldsymbol{\lambda}})$, which lives in the derived category $\mathbf{D}(\cO)$ of $\cO$-modules, and recovers $H^\ast(X_{\GL_n}^U, M_{\boldsymbol{\lambda}})$ after taking cohomology. There is a natural way to lift the operators $T_v^i$ to endomorphisms of the complex $R \Gamma(X_{\GL_n}^U, M_{\boldsymbol{\lambda}})$ in $\mathbf{D}(\cO)$, and we define the algebra $\bbT^S(R \Gamma(X_{\GL_n}^U, M_{\boldsymbol{\lambda}}))$ to be the (commutative) $\cO$-subalgebra of 
\[ \End_{\mathbf{D}(\cO)}(R \Gamma(X_{\GL_n}^U, M_{\boldsymbol{\lambda}})) \]
generated by these operators. Then $\bbT^S(R \Gamma(X_{\GL_n}^U, M_{\boldsymbol{\lambda}}))$ is a finite $\cO$-algebra, and taking cohomology gives rise to a surjective homomorphism
\[ \bbT^S(R \Gamma(X_{\GL_n}^U, M_{\boldsymbol{\lambda}})) \to \bbT^S(H^\ast(X_{\GL_n}^U, M_{\boldsymbol{\lambda}})), \]
which has nilpotent kernel. We consider $\bbT^S(R \Gamma(X_{\GL_n}^U, M_{\boldsymbol{\lambda}}))$ to be the more natural object of study for a number of reasons. First, as our results show, it also receives Galois representations. Second, for any $m \geq 1$ there is a surjective map
\[ \bbT^S(R \Gamma(X_{\GL_n}^U, M_{\boldsymbol{\lambda}})) \to \bbT^S(H^\ast(X_{\GL_n}^U, M_{\boldsymbol{\lambda}} \otimes_\cO \cO/\lambda^m)). \]
Since patching together finite quotients of Hecke algebras plays an essential role in the Taylor--Wiles method, this is a desirable property. For this in action, together with conjectures about existence of Galois representations in this context, see the joint work of Khare and the second named author \cite[Conjecture 6.18]{Tho14}.

We now state our first main theorem:
\begin{theorem}[Theorem \ref{thm_existence_of_galois_in_general_case}]\label{thm_intro_theorem_no_assumption_on_p}
Let $F$ be an imaginary CM field, let $U \subset \prod_v \GL_n(\cO_{F_v})$ be a small\footnote{This condition can be ensured by making $U$ smaller at the places $v|q$  for a rational prime $q \ne p$, see Definition \ref{defsuffsmall}} open compact subgroup, and let $\boldsymbol{\lambda} = (\lambda_\tau)_{\tau \in \Hom(F, E)} \in (\bbZ^n)^{\Hom(F, E)}$. Let $\ffrm \subset \bbT^S(R \Gamma(X_U, M_{\boldsymbol{\lambda}}))$ be a non-Eisenstein maximal ideal. 

Suppose that the $p$-adic places of the maximal totally real subfield $F^+$ of $F$ are all unramified in $F$. Then there exists an ideal $I \subset \bbT^S(R \Gamma(X_{\GL_n}^U, M_{\boldsymbol{\lambda}}))_\ffrm$ satisfying $I^4 = 0$ and a continuous representation 
\[ \rho_\ffrm : G_{F, S} \to \GL_n(\bbT^S(R \Gamma(X_{\GL_n}^U, M_{\boldsymbol{\lambda}}))_\ffrm/I) \]
satisfying the following condition: for each finite place $v\not\in S$ of $F$, the characteristic polynomial of $\rho_\ffrm(\Frob_v)$ is equal to $X^n - T_v^1 X^{n-1} + \dots + (-1)^j q_v^{j(j-1)/2} T_v^j X^{n-j} + \dots + (-1)^n q_v^{n(n-1)/2} T_v^n$ mod $I$.
\end{theorem}
With a stronger assumption on $\boldsymbol{\lambda}$ relative to $p$, we can eliminate the nilpotent ideal $I$ completely, as in our second main theorem:
\begin{theorem}[Theorem \ref{thm_existence_of_galois_in_regular_case}]\label{thm_intro_theorem_p_regular}
Let $F$ be an imaginary CM field in which the prime $p$ is unramified, and let $U = \prod_v U_v \subset \prod_v \GL_n(\cO_{F_v})$ be a small open compact subgroup such that $U_v = \GL_n(\cO_{F_v})$ for each place $v | p$. Let $c \in \Gal(F/F^+)$ denote complex conjugation, and let $\widetilde{I}_p$ denote a set of embeddings $\widetilde{\tau} : F \hookrightarrow E$ such that $\widetilde{I}_p \coprod \widetilde{I}_p c = \Hom(F, E)$. Let $\boldsymbol{\lambda} = (\lambda_\tau)_{\tau \in \Hom(F, E)} \in (\bbZ^n)^{\Hom(F, E)}$, and suppose that for each $\tau \in \Hom(F, E)$, we have
\[ \lambda_{\tau, 1} > \lambda_{\tau, 2} > \dots > \lambda_{\tau, n} \]
and that the condition 
\begin{equation}\label{eqn_intro_p_regularity} [F^+ : \bbQ] n(n + 6 + \sup_{\widetilde{\tau} \in \widetilde{I}_p} (\lambda_{\widetilde{\tau}, 1} + \lambda_{\widetilde{\tau} c, 1}) ) + \sum_{\widetilde{\tau} \in \widetilde{I}_p} \sum_{i=1}^n \left( \lambda_{\widetilde{\tau} , i} - \lambda_{\widetilde{\tau} c, i} - 2 \lambda_{\widetilde{\tau}, n} \right) < p 
\end{equation}
holds. Let $\ffrm \subset \bbT^S(R \Gamma(X_{\GL_n}^U, M_{\boldsymbol{\lambda}}))$ be a non-Eisenstein maximal ideal. Then there exists a continuous representation 
\[ \rho_\ffrm : G_{F, S} \to \GL_n(\bbT^S(R \Gamma(X_{\GL_n}^U, M_{\boldsymbol{\lambda}}))_\ffrm)\]
satisfying the following condition: for each finite place $v\not\in S$ of $F$, the characteristic polynomial of $\rho_\ffrm(\Frob_v)$ is equal to $X^n - T_v^1 X^{n-1} + \dots + (-1)^j q_v^{j(j-1)/2} T_v^j X^{n-j} + \dots + (-1)^n q_v^{n(n-1)/2} T_v^n$.
\end{theorem}
The condition (\ref{eqn_intro_p_regularity}) comes from a theorem of Lan--Suh that will be applied during the proof, as will be explained below.

We now describe the strategy of the proof. We follow Scholze (and the earlier work \cite{Har13}) in first looking at the arithmetic locally symmetric space of the group $G$, the quasi-split unitary group in $2n$ variables over $F^+$ associated to the quadratic extension $F/F^+$. The group $G$ admits a parabolic subgroup $P$ with Levi quotient $M = \Res^F_{F^+} \GL_n$. Writing $U \subset G(\bbA_{F^+}^\infty)$ for a sufficiently small open compact subgroup, $U_P = P(\bbA_{F^+}^\infty) \cap U$, and $U_M$ for the image of $U_P$ in $M(\bbA_{F^+}^\infty)$, we have a diagram of spaces
\begin{equation}\label{eqn_intro_diagram_of_spaces}\begin{aligned} \xymatrix{ X_P^{U_P} \ar[r] \ar[d] & \partial \overline{X}_G^U \ar[r] & \overline{X}_G^U \\ X_M^{U_M} & & X_G^U \ar[u].} 
\end{aligned}
\end{equation}
Here we write $\overline{X}_G^U$ for the Borel--Serre compactification of $X_G^U$, and $\partial \overline{X}_G^U$ for its boundary. Let us write $\bbT^S_G = \cO[ U^S \backslash G(\bbA_{F^+}^{\infty, S}) / U^S]$ for the `abstract' unramified Hecke algebra of $G$, and $\bbT^S_M = \cO[U_M^S \backslash \GL_n(\bbA_F^{\infty, S}) / U_M^S]$ for the abstract unramified Hecke algebra of $\Res^F_{F^+} \GL_n$. If $\mathbf{a} = (a_\tau) \in (\bbZ^{2n})^{\Hom(F^+, E)}$ is a tuple satisfying the condition
\[ a_{\tau, 1} \geq a_{\tau, 2} \geq \dots \geq a_{\tau, 2n} \]
for each $\tau \in \Hom(F^+, E)$, then there is an associated local system $M_{\mathbf{a}}$ of $\cO$-modules on $X_G^U$ (denoted $\underline{A(G; \mathbf{a})}_G^U$ in the body of this article), and the first step is to use the diagram of spaces (\ref{eqn_intro_diagram_of_spaces}) to construct a diagram
\begin{equation}\begin{aligned}\label{eqn_intro_diagram_of_hecke_algebras} \xymatrix{ \bbT^S_G \ar[r] \ar[d]_\cS & \End_{\mathbf{D}(\cO)}(R \Gamma_{\partial \overline{X}_G^U} M_{\mathbf{a}}) \ar[d] \\
\bbT^S_M \ar[r] & \Hom_{\mathbf{D}(\cO)}( R \Gamma_{X_M^{U_M}, c} M_{\boldsymbol{\lambda}}, R \Gamma_{X_M^{U_M}} M_{\boldsymbol{\lambda}})} 
\end{aligned}
\end{equation}
for appropriate choices of $\mathbf{a}$ and $\boldsymbol{\lambda}$. (In order to save space, we have now switched notation from $R \Gamma(X_G^U, -)$ to $R \Gamma_{X_G^U}$.) The map $\cS : \bbT^S_G \to \bbT^S_M$ is the unnormalized Satake transform, given at the level of groups by the slogan `restriction to $P$ and integration along the fibers of $P \to M$'. 

We then show that the natural map $R \Gamma_{X_M^{U_M}, c} M_{\boldsymbol{\lambda}} \to R \Gamma_{X_M^{U_M}} M_{\boldsymbol{\lambda}}$ in $\mathbf{D}(\cO)$ becomes an isomorphism after localizing at $\ffrm$; equivalently, the cohomology of the boundary of the Borel--Serre compactification of $X_M^{U_M}$ vanishes after localization at $\ffrm$. This implies the existence of a homomorphism
\begin{equation}\label{eqn_intro_unnormalized_satake_on_cohomology} \bbT^S_G(R \Gamma_{\partial \overline{X}_G^U} M_{\mathbf{a}})_{\cS^\ast(\ffrm)} \to \bbT^S_M(R \Gamma_{{X}_M^{U_M}} M_{\boldsymbol{\lambda}})_\ffrm. 
\end{equation}
The next step is to construct a Galois group determinant valued in the Hecke algebra $\bbT^S_G(R \Gamma_{\partial \overline{X}_G^U} M_{\mathbf{a}})_{\cS^\ast(\ffrm)}$, or some quotient by a nilpotent ideal. We accomplish this using the exact triangle  in $\mathbf{D}(\cO)$:
\begin{equation}\label{eqn_intro_exact_triangle} \xymatrix@1{ R \Gamma_{X_G^{U}, c} M_{\mathbf{a}} \ar[r] & R \Gamma_{X_G^{U}} M_{\mathbf{a}} \ar[r] & R \Gamma_{\partial \overline{X}_G^{U}} M_{\mathbf{a}}  \ar[r] & R \Gamma_{X_G^{U}, c} M_{\mathbf{a}} [-1].} 
\end{equation}
By reworking Scholze's arguments slightly, we find Galois group determinants valued in $\bbT_G^S(R \Gamma_{X_G^{U}, c} M_{\mathbf{a}})$ and $\bbT_G^S(R \Gamma_{X_G^{U}} M_{\mathbf{a}})$. This leads to a Galois group determinant valued in  $\bbT^S_G(R \Gamma_{\partial \overline{X}_G^U} M_{\mathbf{a}})$, at least at the cost of a nilpotent ideal of square 0, and pushing it along the map (\ref{eqn_intro_unnormalized_satake_on_cohomology}) essentially completes the proof of Theorem \ref{thm_intro_theorem_no_assumption_on_p}. 

To prove Theorem \ref{thm_intro_theorem_p_regular}, we make appeal to the results of Lan--Suh \cite{Lan13}. The main theorems of \emph{op. cit.} imply that under the conditions of Theorem \ref{thm_intro_theorem_p_regular}, the groups $H^i(X_G^U, M_\mathbf{a})$ vanish for $i < D = \frac{1}{2}\dim_\bbR X_G^U = [F^+ : \bbQ] n^2$, and consquently there is an isomorphism of truncations $\tau_{\leq D - 2} R \Gamma_{\partial \overline{X}_G^{U}} M_{\mathbf{a}} \cong \tau_{\leq D - 1}(R \Gamma_{X_G^{U}, c} M_{\mathbf{a}})[-1]$, using the exact triangle (\ref{eqn_intro_exact_triangle}). The diagram (\ref{eqn_intro_diagram_of_hecke_algebras}) is compatible with this truncation, and the map $\tau_{\leq D-2} R \Gamma_{X_M^{U_M}} M_{\boldsymbol{\lambda}} \to R \Gamma_{X_M^{U_M}} M_{\boldsymbol{\lambda}}$ is a quasi-isomorphism since $\dim X_M^{U_M} = D - 1$ and $X_M^{U_M}$ is non-compact. This is enough to give Theorem \ref{thm_intro_theorem_p_regular}.

We note that in all of the theorems proved here, we work with Hecke algebras only after localization at a non-Eisenstein maximal ideal. As we show below, the natural map from compactly cohomology of the $\GL_n$-symmetric space to usual cohomology becomes a quasi-isomorphism after such a localization. On the other hand, Scholze works primarily with interior cohomology (i.e.\ the image of compactly supported cohomology in usual cohomology), which does not seem to have a good derived analogue. Since it is imperative for us to be able to work at the level of complexes rather than at the level of cohomology groups, it seems difficult to avoid this non-Eisenstein condition.

We now describe the structure of this paper. In \S \ref{sec_preliminaries}, we carry out the groundwork necessary to be able to work in a derived setting. In \S \ref{sec_groups_and_symmetric_spaces} we introduce the locally symmetric spaces associated to reductive groups over number fields and discuss their sheaves and cohomology groups. In \S \ref{sec_i_dont_see_any_boundary}, we carry out the important step of showing that the cohomology of the boundary of the $\GL_n$ locally symmetric space vanishes after localizing at a non-Eisenstein maximal ideal. This has been sketched elsewhere, but we give the full details of the argument. Finally, in \S \ref{sec_main_arguments}, we combine all of these ingredients to prove Theorems \ref{thm_intro_theorem_p_regular} and Theorem \ref{thm_intro_theorem_no_assumption_on_p} by carrying out Scholze's perfectoid $p$-adic interpolation argument at the derived level (\S \ref{sec_james}), giving us group determinants at the level of derived Hecke algebras $\bbT^S_G(R \Gamma_{X_G^U} M_{\mathbf{a}})$, and then using the other arguments sketched above to obtain the desired Galois representations for $\GL_n$ (\S\S \ref{sec_application_to_galois_I} -- \ref{sec_application_to_galois_II}). 

\subsection{Notation}

We fix some notation relating to number fields and their Galois groups. A base number field $F$ having been fixed, we will fix an algebraic closure $\overline{F}$ and algebraic closures $\overline{F}_v$ of the completion $F_v$ for every place $v$ of $F$. We also fix embeddings $\overline{F} \hookrightarrow \overline{F}_v$. Writing $G_F = \Gal(\overline{F}/F)$ and $G_{F_v} = \Gal(\overline{F}_v / F_v)$, these embeddings determine continuous embeddings $G_{F_v} \hookrightarrow G_F$ for every place $v$. If $S$ is a finite set of finite places of $F$, then we write $F_S$ for the maximal subfield of $\overline{F}$ unramified outside $S$, and set $G_{F, S} = \Gal(F_S / F)$. It is a quotient of $G_F$. If $v$ is a finite place of $F$, then we will write $\cO_{F_v}$ for the ring of integers of $F_v$, $\varpi_v \in \cO_{F_v}$ for a choice of uniformizer, $k(v) = \cO_{F_v} / (\varpi_v)$ for the residue field, and $q_v = \# k(v)$.

A prime $p$ having been fixed, we will fix an algebraic closure $\overline{\bbQ}_p$ of $\bbQ_p$ and view finite extensions $E/\bbQ_p$ as being subfields of $\overline{\bbQ}_p$. If $E/\bbQ_p$ is such an extension, then we will generally write $\cO$ for its ring of integers, $\pi \in \cO$ for a choice of uniformizer, and $k = \cO / (\pi)$ for the residue field. If $F$ is a field of characteristic 0, then we will write $\epsilon : G_F \to \bbZ_p^\times$ for the usual cyclotomic character.

\subsection{Acknowledgments} 

We thank Peter Scholze for useful conversations about his work \cite{Sch14}. This research was partially conducted during the period that Jack Thorne served as a Clay Research Fellow. James Newton is supported by the ERC Starting Grant 306326.

\section{Preliminaries}\label{sec_preliminaries}

In this section, we will discuss Hecke algebras of locally profinite groups, their module categories, and categories of equivariant sheaves on spaces. We also set up some machinery which constructs natural objects in derived categories of smooth representations for a profinite group, whose cohomology groups are the `completed cohomology' groups (see \cite{emcalsumm}) of a tower of arithmetic locally symmetric spaces, or compactifications of such.

\subsection{Homological algebra}

We first fix notation for derived categories. If $\cA$ is an abelian category with enough injectives, then we write $\mathbf{K}(\cA)$ for the homotopy category of complexes in $\cA$, and $\mathbf{D}(\cA)$ for the corresponding derived category, if it exists. Our normalizations are always cohomological, i.e.\ differentials increase degrees. For other notations (shifts, truncations, etc.\) we follow the conventions of Weibel \cite{Wei94}. We write $\mathbf{K}^+(\cA) \subset \mathbf{K}(\cA)$ for the full subcategory with objects the bounded below complexes, and $\mathbf{D}^+(\cA)$ for its corresponding derived category; it can be identified with the full subcategory of $\mathbf{D}(\cA)$ with objects the bounded below complexes (\cite[Example 10.3.15]{Wei94}). 

If $\cB$ is another abelian category with enough injectives and $F : \cA \to \cB$ is a left exact functor, then the derived functor $RF : \mathbf{D}^+(\cA) \to \mathbf{D}^+(\cB)$ exists (\cite[Theorem 10.5.6]{Wei94}), and is characterized by the following universal property. Let $q_\cA : \mathbf{K}^+(\cA) \to \mathbf{D}^+(\cA)$ and $q_\cB : \mathbf{K}^+(\cB) \to \mathbf{D}^+(\cB)$ be the usual projections, and let $\mathbf{K}F : \mathbf{K}^+(\cA) \to \mathbf{K}^+(\cB)$ be the induced functor on homotopy categories of complexes. Then $RF$ comes equipped with a natural transformation $\xi : q_\cB \mathbf{K} F \to RF q_\cA$ such that for any other functor $G : \mathbf{D}^+(\cA) \to \mathbf{D}^+(\cB)$ equipped with a natural transformation $\zeta : q_\cB \mathbf{K} F \to G q_\cA$, there is a unique natural transformation $\eta : RF \to G$ such that $\zeta_X =  \eta_{q_{\cA}(X)} \circ \xi_X$ for all $X \in \mathbf{K}^+(\cA)$. 

We will often use this universal property in order to compare different functors between derived categories, as in the following lemma.
\begin{lemma}\label{lem_universal_property_of_derived_functor}
Let $\cA, \cB, \cC$ be abelian categories with enough injectives, and let $F : \cA \to \cC, G : \cB \to \cC$ be left exact functors, $i : \cA \to \cB$ an exact functor. Suppose given a natural transformation $\alpha : F \to G\circ i$. Then there is a canonical natural transformation $\eta : RF \to RG \circ i$ (since $i$ is exact, we write $i = Ri$). 
\end{lemma}
\begin{proof}
Let $\xi_F : q_\cC \mathbf{K} F \to RF q_\cA$, $\xi_G : q_\cC \mathbf{K}G \to RG q_{\cB}$, and $\xi_i : q_\cB \mathbf{K}i \to Ri q_\cA$ be the natural transformations that exist by universality. We write $\zeta : q_\cC \mathbf{K} F \to RG Ri q_\cA$ for the natural transformation whose value on $X \in \mathbf{K}^+(\cA)$ is given by the composite
\[ \xymatrix@1{ q_\cC \mathbf{K} F(X) \ar[r]^-{q_\cC(\mathbf{K}\alpha_X)} & q_\cC \mathbf{K} G \mathbf{K} i(X) \ar[r]^-{\xi_{G, \mathbf{K}i(X)}} & RG q_\cB \mathbf{K} i(X) \ar[r]^-{RG(\xi_{i, X})} & RG Ri q_\cA(X). } \]
By the universal property of $RF$, there is a unique natural transformation $\eta : RF \to RG Ri$ with the property that for all $X \in \mathbf{K}^+(\cA)$, $\zeta_X = \eta_{q_\cA(X)} \circ \xi_{F, X}$. This is the $\eta$ of the lemma.
\end{proof}
We now specialize our discussion. Let $R$ be a ring. We allow $R$ to be non-commutative; we will see group algebras and abstract Hecke algebras as examples of such rings. We will write $\Mod(R)$ for the abelian category of (left) $R$-modules, and we will simplify our notation by writing $\mathbf{K}(R)$ etc.\ instead of $\mathbf{K}(\Mod(R))$. If $G$ is a group, then we will write $\Mod(G)$ for the abelian category of $\bbZ[G]$-modules, $\Mod(G, R)$ for the abelian category of $R[G]$-modules, which each have enough injectives, and $\mathbf{K}(G)$, $\mathbf{K}(G, R)$ etc.\ in a similar way. If $G$ is a profinite group, then we will write $\Modsm(G, R)$ for the abelian category of smooth $R[G]$-modules, and $\Modsm(G) = \Modsm(G, \bbZ)$, $\mathbf{K}_\text{sm}(G, R) = \mathbf{K}(\Modsm(G, R))$, etc.\ 

If $G$ is a group and $H \subset G$ is a subgroup, then there are functors $\Ind_H^G : \Mod(H, R) \to \Mod(G,R)$ and $\Res_H^G : \Mod(G,R) \to \Mod(H,R)$, where $\Ind_H^G M = \{ f : G \to M \mid f(hg) = hf(g) \forall h \in H \}$ and $\Res_H^G$ is the usual restriction. We recall that $\Ind_H^G$ is the right adjoint of $\Res_H^G$, that $\Res_H^G$ is exact, and that $\Ind_H^G$ is exact and preserves injectives. The functor $\Res_H^G$ also has a left adjoint $\ind_H^G :  \Mod(H,R) \to \Mod(G,R)$, where $\ind_H^G M = \{ f \in \Ind_H^G M \mid f \text{ finitely supported mod } H \}$. This functor is also exact, showing that $\Res_H^G$ also preserves injectives.

If $N \subset G$ is a normal subgroup, then there is an inflation functor $\Inf_{G/N}^{G} : \Mod(G/N,R) \to \Mod(G,R)$, left adjoint to the functor $\Gamma_N : \Mod(G,R) \to \Mod(G/N,R)$ of $N$-invariants. Inflation is exact, showing that $\Gamma_N$ preserves injectives.

We will introduce more abelian categories (in particular, categories of modules over Hecke algebras and categories of $G$-equivariant sheaves on a space $X$) in the following sections.
\begin{lemma}\label{lem_derived_homtensor}
Let $B \rightarrow R$ and $B \rightarrow C$ be ring maps, with $R$ Noetherian, $B, C$ commutative and $C$ a flat $B$-algebra. Suppose $X, Y \in \mathbf{D}(R)$ are bounded complexes of $R$-modules, with $X$ a bounded complex of finitely generated $R$-modules. Then the natural map \[C\otimes_B\mathrm{Hom}_{\mathbf{D}(R)}(X,Y) \rightarrow \mathrm{Hom}_{\mathbf{D}(C\otimes_B R)}(C\otimes_B X,C\otimes_B Y)\] is an isomorphism.
\end{lemma}
\begin{proof}
This is essentially \cite[Lemma 3]{Zimmermann} (and is probably well known). We denote $C\otimes_B R$ by $R_C$ and similarly denote the functor $\otimes_B C$ by $(-)_C$ (this is an exact functor from $B$-modules to $C$-modules). First we claim that for $M$ a finitely generated $R$-module and $N$ an $R$-module the natural map 
\begin{equation}\label{homtensor}\Hom_R(M,N)_C \rightarrow \Hom_{R_C}(M_C, N_C)\end{equation} is an isomorphism. In fact, this holds without the Noetherian hypothesis on $R$ so long as $M$ is a finitely presented $R$-module. The claim is shown in \cite[Theorem 2.38]{curtis-reiner}. If we consider the functor from finitely generated $R$-modules to $C$-modules given by $M \mapsto \Hom_R(M,N)_C = \Hom_{R_C}(M_C, N_C)$, then the higher derived functors are given by $\mathrm{Ext}^i_R(M,N)_C$ and $\mathrm{Ext}^i_{R_C}(M_C,N_C)$ (since $(-)_C$ preserves projectives; these derived functors are defined because the category of finitely generated $R$-modules has enough projectives). We conclude that the natural maps 
\begin{equation}\label{exttensor} 
\mathrm{Ext}^i_R(M,N)_C \rightarrow \mathrm{Ext}^i_{R_C}(M_C,N_C)
\end{equation} are also isomorphisms (see \cite[8.16]{curtis-reiner}). Note that since the forgetful functor from $\Mod(B)$ to $\Mod(\bbZ)$ is exact, the Ext groups $\mathrm{Ext}^i_R(M,N)$ naturally acquire $B$-module structures, by identifying them with the image on $N$ of the derived functors of $\Hom_R(M,-): \Mod(R)\rightarrow\Mod(B)$.

Next we claim that for a bounded complex $X$ of finitely generated $R$-modules and an $R$-module $N$, the natural map
\[\Hom_{\mathbf{D}(R)}(X,N[0])_C \rightarrow \Hom_{\mathbf{D}(R_C)}(M_C, N_C[0])\] is an isomorphism. We do this by induction on the length $d$ of the complex $X$. For $d = 1$ the claim holds because of the isomorphism (\ref{exttensor}). For the inductive step we do a d\'{e}vissage using truncation functors. Suppose the highest degree in which $H^*(X)$ has a non-zero term is $i$. We have an exact triangle 
\[\tau_{\le i-1}X \rightarrow X \rightarrow H^i(X)[i]\rightarrow \tau_{\le i-1}X[-1]\] and hence a commutative diagram with exact columns
\[\begin{CD} \Hom_{\mathbf{D}(R)}(\tau_{\le i-1}X[-1],N[0])_C @>>> \Hom_{\mathbf{D}(R_C)}(\tau_{\le i-1}X_C[-1],N_C[0])\\
@VVV @VVV \\
\Hom_{\mathbf{D}(R)}(H^i(X)[i],N[0])_C @>>> \Hom_{\mathbf{D}(R_C)}(H^i(X)_C[i],N_C[0])\\
@VVV @VVV \\
\Hom_{\mathbf{D}(R)}(X,N[0])_C @>>> \Hom_{\mathbf{D}(R_C)}(X_C,N_C[0])_C \\
@VVV @VVV \\
\Hom_{\mathbf{D}(R)}(\tau_{\le i-1}X,N[0])_C @>>> \Hom_{\mathbf{D}(R_C)}(\tau_{\le i-1}X_C,N_C[0])\\
@VVV @VVV \\
\Hom_{\mathbf{D}(R)}(H^i(X)[i-1],N[0])_C @>>> \Hom_{\mathbf{D}(R)}(H^i(X)_C[i-1],N_C[0]).\\
\end{CD}\]
By the inductive hypothesis and the five lemma, we are done. Finally, we take our bounded complexes $X, Y$ as in the statement of the lemma. An induction on the length of the complex $Y$ (using the five lemma as above) completes the proof of the lemma.
\end{proof}

\subsection{Hecke algebras}\label{sec_hecke_algebras}

We now introduce the Hecke algebra of a locally profinite group, and discuss various important maps between Hecke algebras in the context of reductive groups over local fields.

\subsubsection{Abstract Hecke algebras}

Let $G$ be a locally profinite group, and let $U \subset G$ be an open compact subgroup. We write $\cH(G, U)$ for the set of compactly supported, $U$-biinvariant functions $f : G \to \bbZ$. 
\begin{lemma}\label{lem_hecke_algebra_acts_on_fixed_vectors}
\begin{enumerate}
\item The $\bbZ$-module $\cH(G, U)$ is in fact an associative $\bbZ$-algebra under convolution, with unit element $[U]$, the characteristic function of $U$.
\item For any $\bbZ[G]$-module $M$, the space $M^U$ of $U$-invariants admits a canonical structure of $\cH(G, U)$-module. This defines a functor $\Gamma_U : \Mod(G) \to \Mod(\cH(G, U))$.
\end{enumerate}
\end{lemma}
We will write $M \mapsto M^\sim$ for the exact functor $\Mod(\cH(G, U)) \to \Mod(\bbZ)$ given by forgetting the $\cH(G, U)$-action.
\begin{proof}
Note that $\cH(G, U)$ is a free $\bbZ$-module, with basis being given by the characteristic functions $[U \alpha U]$ of double cosets $U \alpha U \subset G$. Let us endow $G$ with the unique left-invariant Haar measure giving $U$ volume 1. We observe that $\cH(G, U) \otimes_\bbZ \bbR$ is the space of compactly supported and locally constant $U$-biinvariant functions $f : G \to \bbR$. For functions $f_1, f_2 \in \cH(G, U) \otimes_\bbZ \bbR$, we define their convolution in $\cH(G, U) \otimes_\bbZ \bbR$ by the formula
\begin{equation}\label{eqn_definition_of_convolution} (f_1 \star f_2)(g) = \int_{x \in G} f_1(x) f_2(x^{-1} g) dx. 
\end{equation}
The usual calculation shows that this gives $\cH(G, U) \otimes_\bbZ \bbR$ the structure of assocative algebra with unit $[U]$ (even in the case where $G$ is not unimodular). We now show that the submodule $\cH(G, U)$ is closed under multiplication. It suffices to check this on elements of the form $[U \alpha U]$, $\alpha \in G$; we compute
\[ [U \alpha U] \star [U \beta U] (\gamma) = \int_{x \in U \alpha U} [U \beta U](x^{-1} \gamma) dx = \text{vol}( U \alpha U \cap \gamma U \beta^{-1} U ) = \# ( U \alpha U \cap \gamma U \beta^{-1} U / U), \]
an integer. This shows the first part of the lemma.

For the second part, we note that if $V$ is an $\bbR[G]$-module, then the algebra $\cH(G, U) \otimes_\bbZ \bbR$ acts on $V^U$ by the formula $(v \in V^U, f \in \cH(G, U) \otimes_\bbZ \bbR$):
\[ f \cdot v = \int_{g \in G} f(g) (g \cdot v) dg. \]
If $f = [U \alpha U]$ and $U \alpha U = \coprod_i \alpha_i U$, then this is easily seen to be equal to $\sum_i \alpha_i \cdot v$. We use the same formula to define the action of $[U \alpha U]$ on $M^U$ for any $\bbZ[G]$-module $M$.

This action is clearly functorial in $M$, so to complete the proof of the lemma we just need to show that it is compatible with multiplication of basis elements in $\cH(G, U)$, i.e.\ that for all $m \in M^U$, we have
\begin{equation}\label{eqn_hecke_action_equation} [U \alpha U] \cdot ( [U \beta U] \cdot m) = ([U \alpha U] \cdot  [U \beta U] ) \cdot m. 
\end{equation}
Choose decompositions $[U \alpha U] = \coprod_i \alpha_i U$, $[U \beta U] = \coprod_j \beta_j U$. We see finally that it is enough to show that $[U \alpha U] \cdot [U \beta U] = \sum_{i, j} [\alpha_i \beta_j U]$ as functions $G \to \bbZ$. Evaluating at an element $\gamma \in G$, this is equivalent to the identity
\[ \#( U \alpha U \cap \gamma U \beta^{-1} U / U )= \# \{ (i, j) \mid \gamma \in \alpha_i \beta_j U \}, \]
and this is an elementary exercise in group theory.
\end{proof}
It will be useful to note that the action of $[U \alpha U] \in \cH(G, U)$ on $M^U$, $M$ a $\bbZ[G]$-module, can also be described as the composite
\begin{equation} M^U \to M^{U \cap \alpha U \alpha^{-1}} \to M^U,
\end{equation}
where the first map is given by $v \mapsto \alpha \cdot v$ and the second by $\tr_{U / U \cap \alpha U \alpha^{-1}}$.

\subsubsection{The case of a reductive group}

Now suppose that $F/\bbQ_p$ is a finite extension, and that $G$ is reductive group over $F$; then $G(F)$ is a locally profinite group. We are going to do homological algebra in $\Mod(G(F))$, $\Mod(\cH(G(F), U))$ and related categories. The reader may object that it would be more natural to work, for example, in the abelian category of smooth $\bbZ[G(F)]$-modules. However, in order to understand Hecke actions on cohomology it will suffice for our purposes to work simply with abstract $G(F)$-modules (cf. Corollary \ref{cor_correct_cohomology_groups}).

\subsubsection{Restriction to parabolic subgroup}

Let $P \subset G$ be a rational parabolic subgroup. Suppose moreover that $U \subset G(F)$ satisfies $G(F) = P(F) U$, and set $U_P = P(F) \cap U$. Then we have (for the left-invariant Haar measures $dg$ on $G(F)$ and $dp$ on $P(F)$ giving $U$ and $U_P$ volume 1, respectively) the formula 
\begin{equation}\label{eqn_iwasawa_integration_formula} \int_{g \in G} f(g) dg = \int_{u \in U} \int_{p \in P(F)} f(pu) dp du.
\end{equation}
(For the proof, see \cite[\S 4.1]{Car79}; the proof uses that $G$ is reductive, so $dg$ is also right invariant.) Restriction of functions defines a map $r_P : \cH(G(F), U) \to \cH(P(F), U_P)$.
\begin{lemma}\label{lem_hecke_restriction_to_parabolic} Let $G, P, U$ be as above.
\begin{enumerate}
\item The map $r_P : \cH(G(F), U) \to \cH(P(F), U_P)$ is an algebra homomorphism. 
\item Let $V$ be a $\bbZ[G(F)]$-module, $W$ a $\bbZ[P(F)]$-module, and $f : \Res^{G(F)}_{P(F)} V \to W$ a homomorphism of $\bbZ[P(F)]$-modules. Then the induced map $V^U \to W^{U_P}$  is $r_P$-equivariant in the following sense: for any $t \in \cH(G(F), U)$, $v \in V^U$, we have $f ( t \cdot v) = r_P(t) \cdot f(v)$.
\item Let $W$ be a $\bbZ[P(F)]$-module, and let $V = \Ind_{P(F)}^{G(F)} W$. Then there is a natural isomorphism $V^U \cong r_P^\ast (W^{U_P})$ of $\cH(G(F), U)$-modules.
\end{enumerate}
\end{lemma}
\begin{proof}
For the first part, we can extend scalars to $\bbR$ and calculate for any $\gamma \in P(F)$, $f_1, f_2 \in \cH(G(F), U)$:
\[\begin{split} (f_1 \star_G f_2)(\gamma) = \int_{x \in G(F)} f_1(x) f_2(x^{-1} \gamma) dx = \int_{p \in P(F)} \int_{u \in U} f_1(p u) f_2(u^{-1} p^{-1} \gamma) du dp \\ = \int_{ p \in P(F) } f_1(p) f_2(p^{-1} \gamma) dp = (f_1 \star_P f_2)(\gamma). \end{split}\]
For the second part, we reduce immediately to the universal case $W = \Res^{G(F)}_{P(F)} V$, and must show the formula $t \cdot v = r_P(t) \cdot v$ for any $v \in V^U$. It suffices to check this on basis elements $[U \alpha U]$. Fix a decomposition $U \alpha U = \coprod \alpha_i U$ with $\alpha_i \in P(F)$. It is enough to show that we have in fact $(U \alpha U) \cap P(F) = \coprod \alpha_i U_P$, but this is clear.

For the third part, we observe that 
\[ V^U = \{ f : G(F) \to W \mid \text{ for all }p \in P(F), u \in U, g \in G(F), f(pgu) = p f(g) \}. \]
There is a map $V^U \to W^{U_P}$ given by $f \mapsto f(1)$. This map is injective (since $G(F) = P(F) U$) and surjective (since $U_P = P(F) \cap U$). We must show that for all $t \in \cH(G(F), U)$, we have $(t \cdot f)(1) = r_P(t) f(1)$. This can be checked on basis elements $[U \alpha U]$. Again writing $U \alpha U = \coprod_i \alpha_i U$ with $\alpha_i \in P(F)$, we see that this follows from the formula $r_P([U \alpha U]) = \sum_i [\alpha_i U_P]$.
\end{proof}
Continuing with the notation of the lemma, we observe that there is a diagram of functors
\[ \xymatrix{ \Mod(G(F))\ar[r]^-{\Gamma_U} \ar[d]_{\Res^{G(F)}_{P(F)}} & \Mod(\cH(G(F), U)) \\
\Mod(P(F)) \ar[r]_-{\Gamma_{U_P}} & \Mod(\cH(P(F), U_P)), \ar[u]_{r_P^\ast} } \]
together with a natural transformation $\Gamma_U \to r_P^\ast \circ \Gamma_{U_P} \circ \Res^{G(F)}_{P(F)}$. The vertical functors are exact and the horizontal functors are left exact. Applying the universal property of the derived functor of $\Gamma_U$, we obtain:
\begin{corollary}\label{cor_derived_hecke_morphism_for_parabolic_restriction}
There is a canonical natural transformation $R \Gamma_U \to r_P^\ast R \Gamma_{U_P} \Res^{G(F)}_{P(F)}$. In particular, for any $V \in \Mod(G(F))$, $W \in \Mod(P(F))$ equipped with a morphism $f : \Res^{G(F)}_{P(F)} V \to W$, there is a canonical induced morphism $R \Gamma_U V \to r_P^\ast R \Gamma_{U_P} W$.
\end{corollary}
\begin{proof}
The morphism $R \Gamma_U V \to r_P^\ast R \Gamma_{U_P} W$ is defined as the composite
\[ R \Gamma_U V \to r_P^\ast R \Gamma_{U_P} \Res^{G(F)}_{P(F)} V \to r_P^\ast R \Gamma_{U_P} W, \]
the first arrow by universality and the second by the existence of $f$.
\end{proof}
There is a variant of this involving induction instead of restriction. Indeed, we observe that there is another diagram of functors, commutative up to natural isomorphism:
\[ \xymatrix{ \Mod(G(F)) \ar[r]^-{\Gamma_U} & \Mod(\cH(G(F), U)) \\
\Mod(P(F)) \ar[r]_-{\Gamma_{U_P}}  \ar[u]_{\Ind^{G(F)}_{P(F)}} & \Mod(\cH(P(F), U_P)) \ar[u]_{r_P^\ast} } \]
The vertical functors are exact and the horizontal functors are left exact. Applying the formula for the composition of derived functors (\cite[Corollary 10.8.3]{Wei94}), we obtain:
\begin{corollary}\label{cor_derived_hecke_morphism_for_parabolic_induction}
There is a natural isomorphism $R \Gamma_U \Ind^{G(F)}_{P(F)} \cong r_P^\ast R \Gamma_{U_P}$. In particular, for any $V \in \Mod(P(F))$, there is a canonical isomorphism $R \Gamma_U \Ind^{G(F)}_{P(F)} V \cong r_P^\ast R \Gamma_{U_P} V$.
\end{corollary}
\subsubsection{Projection to Levi quotient}

We now suppose that $P$ is a not necessarily reductive connected linear algebraic group over $F$, with unipotent radical $N$ and reductive quotient $M = P/N$. Choose a Levi decomposition $P = M \ltimes N$, and suppose given an open compact subgroup $U \subset P$ such that $U = (U \cap M(F)) \ltimes (U \cap N(F)) = U_M \ltimes U_N$, say. (We say that $U$ is decomposed with respect to the fixed Levi decomposition of $P$. In this case, $U \cap M(F)$ is identified with the image of $U$ under the projection $P(F) \to M(F)$.) We can then choose left invariant measures $dp, dm$ and $dn$ on the groups $P(F)$, $M(F)$ and $N(F)$, respectively, giving the groups $U$, $U_M$ and $U_N$ measure 1 and satisfying the identity
\begin{equation}\label{eqn_parabolic_integration_formula} \int_{p \in P(F)} f(p) dp = \int_{m \in M(F)} \int_{n \in N(F)} f(mn) dn dm
\end{equation}
\begin{lemma}\label{lem_hecke_integration_along_fibres} Let $P$, $M$, $N$ and $U$ be as above.
\begin{enumerate}
\item Integration along fibers defines an algebra homomorphism $r_M : \cH(P(F), U) \to \cH(M(F), U_M)$.
\item Let $V$ be a $\bbZ[M(F)]$-module, $W$ a $\bbZ[P(F)]$-module, and let $f : \Inf_{M(F)}^{P(F)} V \to W$ be a homomorphism of $\bbZ[P(F)]$-modules. Then the induced map $V^{U_M} \to W^{U}$ is $r_M$-equivariant, in the sense that for all $v \in V^{U_M}$, $t \in \cH(P(F), U)$, we have $f(r_M(t) \cdot v) = t \cdot f(v)$.
\end{enumerate}
\end{lemma}
\begin{proof}
We define a map $r_M :\cH(P(F), U)\otimes_\bbZ \bbR \to \cH(M(F), U_M)\otimes_\bbZ \bbR$ by the formula $r_M(f)(m) = \int_{n \in N(F)} f(mn) dn$. It follows easily from formula (\ref{eqn_parabolic_integration_formula}) that $r_M$ is an algebra homomorphism. To prove the first part of the lemma, it is enough to show that for any $\alpha \in P(F)$, $m \in M(F)$, we have $r_M([ U \alpha U])(m) \in \bbZ$. We calculate
\[ r_M([ U \alpha U])(m)  = \int_{n \in N(F)} [U \alpha U](mn) dn = \int_{n \in m^{-1} U \alpha U \cap N(F)} dn = \# ( m^{-1} U \alpha U \cap N(F)) / U_N, \]
an integer.

For the second part of the lemma, it is enough to consider the case where $W=\Inf_{M(F)}^{P(F)} V$. Let $\alpha \in P(F)$, and choose a decomposition $U \alpha U = \coprod_i \alpha_i U$. We claim that $r_M([U \alpha U]) = \sum_i [ \overline{\alpha}_i U_M ]$, where $\overline{\alpha}_i$ denotes the image of $\alpha_i$ in $M(F)$. This follows from the easily verified formula
\[ \int_{n \in N(F)} [\alpha_i U](mn) dn = \# ( m^{-1} \alpha_i U_M U_N \cap N(F) ) / U_N = [ \overline{\alpha}_i U_M ](m). \]
For any $v \in V^{U_M}$, we thus have $[U \alpha U] \cdot v = \sum_i \alpha_i \cdot v = \sum_i \overline{\alpha}_i \cdot v = r_M([U \alpha U]) \cdot v$. This completes the proof of the lemma.
\end{proof}
Let us continue with the notation of the above lemma. We have constructed a diagram of functors, commutative up to natural isomorphism:
\[ \xymatrix{ \Mod(P(F))\ar[r]^-{\Gamma_U} & \Mod(\cH(P(F), U)) \\
\Mod(M(F)) \ar[r]_-{\Gamma_{U_M}} \ar[u]^{\Inf_{M(F)}^{P(F)}} & \Mod(\cH(M(F), U_M)) \ar[u]_{r_M^\ast}. } \]
The vertical functors are exact, and the horizontal functors are left exact. We deduce:
\begin{corollary}\label{cor_universality_of_derived_functor}
There is a canonical natural transformation $r_M^\ast \circ R\Gamma_{U_M} \to R\Gamma_U \circ \Inf_{M(F)}^{P(F)}$.
\end{corollary}
\begin{proof}
By Lemma \ref{lem_universal_property_of_derived_functor}, there is a canonical natural transformation
\[ r_M^\ast \circ R \Gamma_{U_M} \cong R(r_M^\ast \circ \Gamma_{U_M}) \cong R(\Gamma_{U} \circ \Inf_{M(F)}^{P(F)}) \to R\Gamma_U \circ \Inf_{M(F)}^{P(F)}. \]
\end{proof}

\subsubsection{Adeles}
All of the results in this section have obvious analogues for Hecke algebras $\cH(G(\bbA_F^S), U^S)$, where now $F$ is a number field, $G$ is a connected linear algebraic group over $F$, $S$ is a finite set of places of $F$ containing the infinite places, and $U^S \subset G(\bbA_F^S)$ is an open compact subgroup. We omit the formulation of these generalizations. If $U^S = \prod_{v \not\in S} U_v$ decomposes as a product, then we have the usual decomposition of this global Hecke algebra as a restricted tensor product of local Hecke algebras:
\[ \cH(G(\bbA_F^S), U^S) = \otimes'_{v \not\in S} \cH(G(F_v), U_v). \]
If $X$ is a complex in $\mathbf{D}(\cH(G(\bbA_F^S), U^S))$, then there is a canonical homomorphism $T_G : \cH(G(\bbA_F^S), U^S) \to \End_{\mathbf{D}(\bbZ)}(X^\sim)$, and similarly with $\bbZ$ replaced by any commutative ring $R$ of coefficients. Indeed, for any $t \in \cH(G(\bbA_F^S), U^S)$, the module structure on $X$ defines a map $X^\sim \to X^\sim$ of complexes, hence an element $T_G(t) \in \End_{\mathbf{D}(\bbZ)}(X^\sim)$. It is easy to check that this is independent of choices in the sense that if $X \to Y$ is a quasi-isomorphism in $\mathbf{D}(\cH(G(\bbA_F^S), U^S))$, the elements $T_G(t)$ of $\End_{\mathbf{D}(\bbZ)}(X^\sim) \cong \End_{\mathbf{D}(\bbZ)}(Y^\sim)$ are identified. We will use this observation in our construction of Hecke algebras.
\subsubsection{Application when $G$ is unramified}\label{sec_application_when_G_unramified}
To obtain situations where the results of this section apply, let us now assume again that $F$ is a finite extension of $\bbQ_p$, and consider an unramified reductive group $G$ over $F$ with a reductive model $\underline{G}$ over $\cO_{F}$. Thus $\underline{G}$ is affine and smooth over $\cO_{F}$ with connected reductive fibres. We fix a choice $\underline{S} \subset \underline{G}$ of maximal $\cO_{F}$-split torus, as well as a choice $\underline{B} \subset \underline{G}$ of Borel subgroup containing $\underline{T} = Z_{\underline{G}}(\underline{S})$. Then the group $X^\ast(\underline{S})$ of $\cO_{F}$-rational characters is a finite free $\bbZ$-module, and contains the subset $\Phi(G, S) = \Phi(\underline{G}, \underline{S})$ of $F$-rational roots. The choice of Borel subgroup determines a root basis $R^\text{rat} \subset \Phi(\underline{G}, \underline{S})$, and the $G(F)$-conjugacy classes of parabolic subgroups are in bijection with the subsets $I \subset R^\text{rat}$. A representative of the conjugacy class corresponding to a given $I$ is given by the generic fiber of the closed subgroup
\begin{equation}\label{eqn_integral_levi_decomposition} \underline{P}_I = \underline{M}_I \ltimes \underline{N}_I 
\end{equation}
of $\underline{G}$, where $\underline{P}_I$ contains $\underline{B}$ and $\underline{M}_I$ is the unique Levi subgroup of $\underline{P}_I$ containing $\underline{T}$. (The existence of the tori $\underline{S}$ and $\underline{T}$ can be justified using the usual structure theory for the special fibre of $\underline{G}$ and then \cite[Corollary B.3.5]{Con14} (lifting of tori). The existence and properties of the parabolic subgroups then follow from \cite[Theorem 4.1.7]{Con14}.)
\begin{lemma}
Let $U = \underline{G}(\cO_{F})$. Then $U$ is a hyperspecial maximal compact subgroup of $G(F)$ satisfying the following conditions:
\begin{enumerate}
\item For each $I \subset R^\text{rat}$, the subgroup $U_{P_I} = U \cap P_I(F)$ is decomposed with respect to the given Levi decomposition, i.e.\ $U_{P_I} = U_{M_I} \ltimes U_{N_I}$.
\item We have $G(F) = P(F) U$. 
\end{enumerate}
\end{lemma}
\begin{proof}
The first part is immediate from the decomposition (\ref{eqn_integral_levi_decomposition}), since $U_{P_I} = \underline{P}_I(\cO_F)$. The second part is the Iwasawa decomposition \cite[3.3.2]{Tit79}.
\end{proof}
In the situation of the lemma, we therefore obtain (using Lemma \ref{lem_hecke_restriction_to_parabolic} and Lemma \ref{lem_hecke_integration_along_fibres}) a canonical homomorphism 
\[ \cS = r_{M_I} \circ r_{P_I} : \cH(G(F), U) \to \cH(M_I(F), U_M), \]
which we call the unnormalized Satake transform. 

\subsubsection{Representations of $U$}\label{sec_representations_of_u}

Keeping the notation of the previous section, we now move in a slightly different direction and describe some interesting $\bbZ[U]$-modules that will later be used to define Hecke-equivariant coefficient systems on arithmetic locally symmetric spaces. Thus $F$ is a finite extension of $\bbQ_p$ and $\underline{G}$ is a reductive group over $\cO_F$ with generic fibre $G$, equipped with maximal torus $\underline{T}$ and Borel subgroup $\underline{B}$. 

Let $E/F$ be a finite extension that splits $T$ (and therefore $G$). Then the choice of Borel subgroup $T_E \subset B_E$ determines a root basis $R^\text{abs} \subset \Phi(G_E, T_E)$, and we write $X^\ast(T_E)^+ \subset X^\ast(T_E)$ for the set of $B$-dominant weights, i.e.\ the set of $\lambda \in X^\ast(T_E)$ satisfying the condition $\langle \lambda, \alpha^\vee \rangle \geq 0$ for all $\alpha \in R^\text{abs}$. We write $X^\ast(T_E)^{++} \subset X^\ast(T_E)^{+}$ for the set of regular dominant weights, i.e.\ satisfying the condition $\langle \lambda, \alpha^\vee \rangle > 0$ for all $\alpha \in R^\text{abs}$. We also define
\[ X^\ast(T_E)^{< p} = \{ \lambda \in X^\ast(T_E) \mid \langle \lambda + \rho, \alpha^\vee \rangle \leq p \text{ }\forall \alpha \in \Phi^\text{abs}\}, \]
and $X^\ast(T_E)^{+, <p} = X^\ast(T_E)^+ \cap X^\ast(T_E)^{< p}$, and similarly for $X^\ast(T_E)^{++, <p}$.

Let $\cO$ denote the ring of integers of $E$, $(\pi) \subset \cO$ its maximal ideal, and $k = \cO/(\pi)$ its residue field. If $\lambda \in X^\ast(T_E) = X^\ast(\underline{T}_\cO)$, then we write $B(G; \lambda)$ for the functor defined on $\cO$-algebras $R$ by the formula
\[ B(G; \lambda)(R) = (\text{ind}_{\underline{B}^-}^{\underline{G}} \lambda)(R) = \{ f \in R[\underline{G}] \otimes_\cO \cO(\lambda) \mid \text{ for all }R \to A, f \otimes_R A \in (A[G] \otimes_\cO \cO(\lambda))^{\underline{B}^-(A)} \}. \]
We write $\underline{B}^-$ for the opposite Borel subgroup to $\underline{B}$. This functor is defined and studied in \cite[I.3.3]{Jan03}. In particular, $A(G; \lambda) = B(G; \lambda)(\cO)$ is an $\cO[U]$-module, finite free as $\cO$-module (it is finitely generated by \cite[I.5.12(c)]{Jan03}, and is then clearly free); and if $\lambda \in X^\ast(T_E)^+$, then it follows from \cite[II.4.5]{Jan03} that for any $\cO$-algebra $R$ the natural map 
\[ A(G; \lambda) \otimes_\cO R \to B(G; \lambda)(R) \]
is an isomorphism.
\begin{proposition}\label{prop_direct_sum_decomposition_on_restriction_to_levi} Let $\lambda \in X(T_E)^+$.
\begin{enumerate}
\item The module $A(G; \lambda) \otimes_\cO E$ is an absolutely irreducible $E[U]$-module.
\item Let $I \subset R^\text{rat}$, so the parabolic subgroup $\underline{P}_I = \underline{M}_I \ltimes \underline{N}_I \subset \underline{G}$ is defined. Then there is a direct sum decomposition
\[ \Res_{U_{M_I}}^U A(G; \lambda) = A(M_I; \lambda) \oplus K \]
of $\cO[U_{M_I}]$-modules, and $U_{N_I}$ acts trivially on $A(M_I, \lambda)$.
\end{enumerate}
\end{proposition}
\begin{proof}
The first part is a consequence of highest weight theory in characteristic 0 and the Zariski density of $U \subset \underline{G}(E)$. For the second part, we observe that since $\underline{T}_\cO$ is a torus, there is a decomposition
\[ \Res_{\underline{T}_\cO}^{\underline{G}_\cO} B(G; \lambda) = \oplus_{\mu} M_\mu \]
as a sum of finitely many non-zero weight spaces (even over $\cO$). We define $W_1 = \oplus_{\mu \in \bbZ_{\geq 0} \cdot I} M_{\lambda - \mu}$, and $W_2$ to be the sum of the complementary weight spaces. Then there is a decomposition $\Res_{\underline{T}_\cO}^{\underline{G}_\cO} B(G; \lambda) = W_1 \oplus W_2$. 

We claim that this is a decomposition of $\underline{M}_{I, \cO}$-modules, that $\underline{N}_{I, \cO}$ acts trivially on $W_1$, and that there is an isomorphism $W_1 \cong B(M_I; \lambda)$ of $\underline{M}_{I, \cO}$-modules. Let us address each point in turn. By the main result of \cite{Cab84}, as well as \cite[4.1, Proposition]{Cab84}, we know that 
\[ \Res^{G_E}_{M_{I, E}} B(G; \lambda) = W_1 \otimes_\cO E \oplus W_2 \otimes_\cO E \]
as $M_{I, E}$-modules, that
\[ (B(G; \lambda) \otimes_\cO E)^{N_{I, E}} = W_1 \otimes_\cO E, \]
and that $W_1 \otimes_\cO E \cong B(M_I; \lambda) \otimes_\cO E$ as $M_{I, E}$-modules. The $\underline{M}_{I, \cO}$-invariance of the decomposition $B(G; \lambda) = W_1 \oplus W_2$ can be checked on $E$-points, so follows from what we have written above. The fact that $\underline{N}_{I, \cO}$ acts trivially on $W_1$ can also be checked on $E$-points.

It remains to check that there is an isomorphism $W_1 \cong B(M_I; \lambda)$ of $\underline{M}_{I, \cO}$-modules. By Frobenius reciprocity (i.e.\ \cite[I.3.4, Proposition]{Jan03}), we have for any $\underline{M}_{I, \cO}$-module $V$ an isomorphism
\[ \Hom_{\underline{M}_{I, \cO}}(V, B(M; \lambda)) = \Hom_{\underline{B}^-_\cO}(V, \cO(\lambda)). \]
The module $W_1$ has highest weight $\lambda$, so $W_1^\vee$ has lowest weight $- \lambda$, hence there is a non-zero $\underline{B}^-_\cO$-equivariant homomorphism $\cO(-\lambda) \to W_1^\vee$ (by \cite[II.1.19(7)]{Jan03}), hence a non-zero $\underline{B}^-_\cO$-equivariant homomorphism $W_1 \to \cO(\lambda)$, hence (by Frobenius reciprocity) a non-zero $\underline{M}_{I, \cO}$-equivariant homomorphism $f : W_1 \to B(M; \lambda)$. We can assume that $f \otimes_\cO k \neq 0$. We claim that $f$ is the desired isomorphism. We know that $f$ is an isomorphism after extending scalars to $E$, so it is enough to show that the map $f \otimes_\cO k$ is injective.

Suppose for contradiction that $\ker (f \otimes_\cO k) \neq 0$. Then the kernel of $f$ has a non-zero $\underline{B}$-socle, so contains the $\underline{B}$-socle of $W_1 \otimes_\cO k = B(G; \lambda)^{\underline{N}_{I, \cO}} \otimes_\cO k$, which equals $B(G; \lambda)^{\underline{N}_{\emptyset, \cO}} \otimes_\cO k = k(\lambda)$, by \cite[II.2.2, Proposition]{Jan03}. We deduce that $f$ determines a non-zero element of the group
\[ \Hom_{\underline{M}_{I, \cO}}(W_1 \otimes_\cO k / \ker (f \otimes_\cO k), B(G; \lambda) \otimes_\cO k) = \Hom_{\underline{B}^-_\cO}(W_1 \otimes_\cO k / \ker (f \otimes_\cO k), k(\lambda)) \]
which contradicts the fact that $W_1 \otimes_\cO k / \ker (f \otimes_\cO k)$ does not contain the weight $\lambda$, which occurs with multiplicity 1 in $W_1 \otimes_\cO k$. This contradiction shows that $f \otimes_\cO k$ is injective, and concludes the proof.
\end{proof}
We now change notation slightly, and suppose that $E \subset \overline{\bbQ}_p$ is a finite extension of $\bbQ_p$ which contains the image of all continuous embeddings $F \hookrightarrow \overline{\bbQ}_p$. If $\tau \in \Hom(F, E)$, then the above construction gives an $\cO[U]$-module $A(G; \lambda)$ for each $\lambda \in X^\ast(T_{E, \tau})^+$, where the subscript $\tau$ indicates that we extend scalars from $F$ to $E$ via the embedding $\tau$. 

If $\boldsymbol{\lambda} = (\lambda_{\tau})_{\tau \in \Hom(F, E)}$ is a tuple with $\lambda_{\tau} \in X^\ast(T_{E, \tau})^+$ for each $\tau \in \Hom(F, E)$, then we define $A(G; \boldsymbol{\lambda}) = \otimes_{\tau} A(G; \lambda_\tau)$, the tensor product being over $\cO$. Then $A(G; \boldsymbol{\lambda})$ is an $\cO[U]$-module, finite free over $\cO$, and $A(G; \boldsymbol{\lambda}) \otimes_\cO E$ has a natural structure of absolutely irreducible $E[G(F)]$-module. Proposition \ref{prop_direct_sum_decomposition_on_restriction_to_levi} now implies the following result:
\begin{corollary}\label{cor_local_tensor_product_of_weyl_modules}
Let $I \subset R^\text{rat}$ and let $\boldsymbol{\lambda} \in \prod_{\tau \in \Hom(F, E)} X^\ast(T_{E, \tau})^+$. Then there is a canonical decomposition $\Res^U_{U_{M_I}} A(G; \boldsymbol{\lambda}) = A(M; \boldsymbol{\lambda}) \oplus K$ of $\cO[U_{M_I}]$-modules, where $A(M; \boldsymbol{\lambda}) \subset A(G; \boldsymbol{\lambda})^{U_{N_I}}$.
\end{corollary}
\begin{proof}
Since tensor products respect direct sums, this is an immediate consequence of Proposition \ref{prop_direct_sum_decomposition_on_restriction_to_levi}.
\end{proof}

\subsection{Equivariant sheaves for abstract groups}\label{sec_equivariant_sheaves}

Let $X$ be a topological space, and let $G$ be a group that acts on the right on $X$ by homeomorphisms. (We call $X$ a $G$-space.) In this section, we consider (essentially following \cite[Ch. V]{Gro57}) the derived category of $G$-equivariant sheaves on $X$.
\begin{definition}
A $G$-equivariant sheaf on $X$ is a sheaf $\cF$ on $X$ equipped with isomorphisms $[g]_\cF : \cF \to g^\ast \cF$ for each $g \in G$, all satisfying the following conditions:
\begin{enumerate}
\item If $e \in G$ is the identity, then $[e]_\cF$ is the identity.
\item For each $g, h \in G$, we have $[g h]_\cF = g^\ast[h]_\cF \circ [g]_\cF$.
\end{enumerate}
We write $Sh_G(X)$ for the category of $G$-equivariant sheaves of abelian groups on $X$. If $R$ is a ring, then we write $\Sh_G(X, R)$ for the category of $G$-equivariant sheaves of $R$-modules on $X$.
\end{definition}
It is easy to see that $\text{Sh}_G(X,R)$ is an abelian category, and that the natural functor $\text{Sh}_G(X,R) \to \text{Sh}(X,R)$ (which forgets the $G$-action) commutes with the formation of kernels and cokernels.
\begin{lemma}
For any ring $R$, the category $\Sh_G(X,R)$ has enough injectives.
\end{lemma}
\begin{proof}
We just give the argument in the case $R = \bbZ$. Let $\cF \in \Sh_G(X)$. We must construct a monomorphism $\cF \hookrightarrow \cI$ for some injective object $\cI$. Choose for each orbit $y \in Y = X/G$ a representative $\xi(y) \in X$ and a monomorphism $\cF_{\xi(y)} \to A_y$, for some injective $G_{\xi(y)}$-module $A_y$. We then define 
\[ I_y = \Ind_{G_{\xi(y)}}^G A_y = \{ f : G \to A_y \mid \forall h \in H, g \in G, f(gh) = h^{-1} f(g) \}. \]
We interpret $I_y$ as a product of skyscraper sheaves supported on the orbit $y$, with stalk over $g \xi(y)$ given by the set of functions with support in $g G_{\xi(y)}$. It has a natural structure of $G$-equivariant sheaf. We define $\cI = \prod_{y \in Y} I_y$. Then there is a natural $G$-equivariant inclusion $\cF \hookrightarrow \cI$ and for any $\cG \in \Sh_G(X)$, we calculate
\[ \Hom_{\Sh_G(X)}(\cG, \cI) = \prod_{y \in Y} \Hom_{G_{\xi(y)}}(\cG_{\xi(y)}, A_y). \]
It follows that $\cI$ is injective, and this completes the proof of the lemma.
\end{proof}
To avoid a proliferation of notation, we now restrict to the case $R = \bbZ$. Everything we say has a clear analogue for the category $\Sh_G(X, R)$. If $H \subset G$ is a subgroup, then there is a natural restriction functor $\Res^G_H : \Sh_G(X) \to \Sh_H(X)$. We define a functor $\ind_H^G : \Sh_H(X) \to \Sh_G(X)$ as follows. Let $p : G \times X \to X$ denote projection to the second factor, and let $G \times H$ act on $G \times X$ by the formula $(g, h) \cdot (g', x) = (g g' h^{-1}, hx)$. Then the sheaf $p^\ast \cF$ admits a natural structure of $G \times H$-equivariant sheaf, and therefore descends naturally to a $G$-equivariant sheaf $\cF'$ on the quotient $G \times_H X$ (see Lemma \ref{lem_descent_along_free_actions} below). The induced map $f : G \times_H X \to X$ is a $G$-equivariant local homeomorphism and we define $\ind_H^G \cF$ to be the subsheaf of $f_\ast \cF'$ consisting of sections which stalkwise are supported in finitely many of copies of $X$ under the isomorphism $G \times_H X \cong \sqcup_{G / H} X$. (We use the notation $\ind$ instead of $\Ind$ as the functor $\ind_H^G$ plays the role of compact induction.)

If $Y$ is another space with $G$-action, and $f : X \to Y$ is a $G$-equivariant continuous map, then the usual pushforward and pullback of sheaves gives rise to functors $f_\ast : \Sh_G(X) \to \Sh_G(Y)$ and $f^\ast : \Sh_G(Y) \to \Sh_G(X)$. If $Y$ is a point, then we identify $\Sh_G(Y) = \Mod(G)$ and write $f_\ast = \Gamma_X$. (The left action of $G$ on $\Gamma_X \cF$ is given by the formula $g \cdot s = [g]_\cF^{-1} g^\ast s$.)
\begin{lemma}\label{lem_adjunctions} Let notation be as above.
\begin{enumerate}
\item The functors $(\ind_H^G, \Res_H^G)$ form an adjoint pair, and both $\ind_H^G$ and $\Res_H^G$ are exact.
\item The functors $(f^\ast, f_\ast)$ form an adjoint pair, and $f^\ast$ is exact.
\end{enumerate}
\end{lemma}
\begin{proof}
It is clear from the definition that there is a natural map $\cF \to \Res^G_H \ind_H^G \cF$ for any $\cF \in \Sh_H(X)$, and this gives rise to the desired adjunction. It is useful to note that the stalks of the induced sheaf can be calculated as 
\[ (\ind_H^G \cF)_x = \{ (s_g)_{g \in G} \mid s_g \in \cF_{g^{-1} x}, \forall h \in H, s_{gh} = h s_g, \text{ finitely supported modulo }H \}. \]
There is an isomorphism of underlying sheaves $\ind_H^G \cF \cong \oplus_{g \in G/H} \, g \otimes \cF$. This makes it clear that both $\Res_H^G$ and $\ind_H^G$ are exact, and proves the first part of the lemma. The second part follows easily from the corresponding result when $G$ is the trivial group.
\end{proof}
\begin{corollary}
The functors $\Res^G_H : \Sh_G(X) \to \Sh_H(X)$ and $f_\ast : \Sh_G(X) \to \Sh_G(Y)$ preserve injectives.
\end{corollary}
\begin{definition}\label{defn_free_g_space} Let $X$ be a $G$-space. We say that $X$ is free if the action of $G$ satisfies the following condition: every point $x \in X$ has a neighbourhood $U$ such that for all $g \in G - \{ e \}$, $g U \cap U = \emptyset$. This implies in particular that every point $x \in X$ has trivial stabilizer.
\end{definition}

If $\varphi : G \to H$ is a surjective homomorphism with kernel $K$, and $X$ is a $G$-space, and $Y$ is an $H$-space, then we say that a map $f : X \to Y$ is $\varphi$-equivariant if we have $f(xg) =  f(x) \varphi(g)$ for all $x \in X$, $g \in G$. In this case, we define a functor $f_\ast^K : \Sh_G(X) \to \Sh_H(Y)$ by the formula $f_\ast^K(\cF) = f_\ast(\cF)^K$ (i.e.\ $f_\ast^K(\cF) \subset f_\ast(\cF)$ is the subsheaf of $K$-invariants).
\begin{lemma}\label{lem_descent_along_free_actions} Let $\varphi : G \to H$ be a surjective homomorphism with kernel $K$.
\begin{enumerate}\item Suppose that $f : X \to Y$ is a $\varphi$-equivariant continuous map. Then the functors $f^\ast : \Sh_H(Y) \to \Sh_G(X)$, $f_\ast^K : \Sh_G(X) \to \Sh_H(Y)$ form an adjoint pair.
 \item Suppose instead that $X$ is a $G$-space on which $K$ acts freely, and $Y = X/K$, endowed with its quotient topology. Then the two functors $f_\ast^K : \Sh_G(X) \to \Sh_H(Y)$, $f^\ast : \Sh_H(Y) \to \Sh_G(X)$, are mutually inverse equivalences of categories.
\end{enumerate} 
\end{lemma}
\begin{proof}
The first part is \cite[Proposition 8.4.1]{Ber94}. The second part follows from \cite[Lemma 8.5.1]{Ber94}.
\end{proof}
Now suppose that $X$ is a $G$-space and that $G$ is a locally profinite group, and let $U \subset G$ be an open compact subgroup that acts freely on $X$. As we have seen, there is a left exact functor $\Gamma_U : \Mod(G) \to \Mod(\cH(G, U))$, $M \mapsto M^U$. We obtain a diagram of functors, commutative up to natural isomorphism:
\begin{equation}\label{eqn_hecke_functor_diagram}\begin{aligned} \xymatrix{ \Sh_G(X) \ar[r]^-{\Gamma_X} \ar[rd]_{\Res^G_U} & \Mod(G) \ar[r]^-{\Gamma_U} & \Mod(\cH(G, U)) \ar[r]^-{(\cdot)^\sim} & \Mod(\bbZ) \\
& \Sh_U(X) \ar[r]_{f_\ast^U} & \Sh(X/U) \ar[ur]_{\Gamma_{X/U}} } 
\end{aligned}
\end{equation}
 The functors $f_\ast^U$ and $\Res^G_U$ are exact and preserve injectives. As a formal consequence, we obtain:
\begin{proposition}\label{prop_hecke_action_on_cohomology}
With notation as above, there is a canonical isomorphism in $\mathbf{D}(\bbZ)$, for any $\cF \in \Sh_G(X)$:
\[ R (\Gamma_U \circ \Gamma_X)(\cF)^\sim \cong R \Gamma_{X/U}(f_\ast^U \cF). \]
\end{proposition}
We will often use the following slightly weaker consequence of the proposition: for any $\cF \in \Sh_G(X)$, there is a canonical homomorphism $T_G: \cH(G, U) \to \End_{\mathbf{D}(\bbZ)}(R \Gamma_{X/U} f_\ast^U \cF)$. (In the context of arithmetic locally symmetric spaces, such homomorphisms recover the usual action of Hecke operators on cohomology. We turn to this topic in \S \ref{sec_groups_and_symmetric_spaces}.) The above homomorphism can be given explicitly on basis elements as follows. We recall (cf. \S \ref{sec_hecke_algebras}) that the algebra $\cH(G,U)$ is free over $\bbZ$, a basis being given by the elements $[U \alpha U]$ with $\alpha \in G$. 

Let $V = U\cap\alpha U\alpha^{-1}$, let $p_1 : X/V \rightarrow X/U$ denote the natural projection, $p_2 : X/V \rightarrow X/U$ the map $X/V \to X/\alpha^{-1} V \alpha \to X/U$ given by acting by $\alpha$, then projecting. Both $p_1$ and $p_2$ are topological covering maps. If $\cF \in \Sh_G(X)$, then the isomorphism $\cF\cong\alpha^*\cF$ induces an isomorphism 
\[p_1^* f_*^U\cF = f_*^V\cF \cong \alpha^* f_*^{\alpha^{-1}V\alpha}\cF = p_2^* f_*^U\cF.\]  
We define an endomorphism $\theta(\alpha)$ of $R\Gamma_{X/U}(f_*^U\cF)$ as the composite
\[R\Gamma_{X/U}(f_*^U\cF) \overset{p_2^*}{\rightarrow} R\Gamma_{X/V}(p_2^*f_*^U\cF) \cong R\Gamma_{X/V}(p_1^*f_*^U\cF) \cong R\Gamma_{X/U}(p_{1,*}p_1^*f_*^U\cF) \rightarrow  R\Gamma_{X/U}(f_*^U\cF)\] where the final map is the trace, defined by the adjunction $(p_{1,*} = p_{1,!}, p_1^* = p_1^!)$. We refer to \cite[Theorem 3.1.5, Proposition 3.3.2]{Kas94} for the adjunction property and the identification $p_1^* = p_1^!$. 
\begin{lemma}\label{lem_hecke_action_as_expected}
Let $\cF \in \Sh_{G}(X)$. For $\alpha \in G$ the image of $[U\alpha U]$ in $\End_{\mathbf{D}(\bbZ)}(R\Gamma_{X/U} f_*^U\cF)$ equals $\theta(\alpha)$.
\end{lemma}
\begin{proof}
It suffices to check the same statement for $\Gamma_{X/U}(f_*^U\cF)$, since applying this to the sheaves appearing in an injective resolution $\cF \rightarrow \cI^\bullet$ gives the desired result. The lemma can then be proved by comparing the explicit descriptions of the Hecke action on $U$-invariants and trace map on global sections.
\end{proof}
We now present a kind of `Shapiro's lemma' for spaces. Let $G$ be a group, $H \subset G$ a subgroup, and $X$ an $H$-space.
\begin{proposition}\label{prop_shapiros_lemma_for_spaces}
There is a natural equivalence of categories $\Ind_H^G : \Sh_H(X) \cong \Sh_G(G \times_H X)$, and a natural isomorphism $\Gamma_{G \times_H X} \circ \Ind_H^G \cong \Ind_H^G \circ \Gamma_X$.
\end{proposition}
\begin{proof}
Let $\pi_1 : G \times X \to X$ and $\pi_2 : G \times X \to G \times_H X$ be the two projections. Let $G \times H$ act on $G \times X$ by the formula $(g, h) (g', x) = (h^{-1} g' g, x h)$. Then the subgroup $G \times \{ 1 \}$ acts freely on $G \times X$, and $\pi_1$ is the corresponding quotient map; and the subgroup $\{ 1 \} \times H$ acts freely on $G \times X$, and $\pi_2$ is the corresponding quotient map. We obtain a diagram of functors
\[ \xymatrix{ \Sh_H(X) \ar[r]^-{\pi_1^\ast} & \Sh_{G \times H}(G \times X) \ar[r]^-{\pi_{2, \ast}^H} & \Sh_G(G \times_H X).} \]
It follows from Lemma \ref{lem_descent_along_free_actions} that $\Ind_H^G = \pi_{2, \ast}^H \circ \pi_1^\ast$ is an equivalence of categories, with inverse given by $\pi_{1, \ast}^G \circ \pi_2^\ast$. The natural isomorphism $\Gamma_{G \times_H X} \circ \Ind_{H}^G \cong \Ind_H^G \circ \Gamma_X$ is then an easy consequence of the definitions. 
\end{proof}
\begin{corollary}\label{cor_derived_induction_of_equivariant_sheaves}
With notation as in the proposition, we have a natural isomorphism of derived functors $R \Gamma_{G \times_H X}  \circ \Ind_{H}^G \cong \Ind_H^G \circ R \Gamma_X$.
\end{corollary}
\begin{proof}
This follows from Proposition \ref{prop_shapiros_lemma_for_spaces} and the formula for the composition of derived functors.
\end{proof}
\subsection{Equivariant sheaves for topological groups}\label{sec_top_eq_sheaves}
We will also consider $G$-equivariant sheaves where $G$ is a topological group acting continuously on a topological space $X$, following \cite{Ber94}. 
\begin{definition}
Let $G$ be a topological group and $X$ be a topological space. 
\begin{enumerate}\item We say that $X$ is a $G$-space if it is equipped with a continuous right action of $G$, i.e.\ the multiplication map $m: X \times G \rightarrow X$ is continuous. Write $p: X \times G \rightarrow X$ for the projection map. 

\item If $X$ is a $G$-space, a $G$-equivariant sheaf on $X$ is a sheaf $\cF$ on $X$ equipped with an isomorphism $\theta: p^*\cF \cong m^*\cF$ satisfying the usual cocycle condition (see \cite[\S0.2]{Ber94} for the analogous formula in the case of a left action).
\end{enumerate}
We write $\Sh_G(X)$ for the abelian category of $G$-equivariant sheaves of abelian groups on $X$. For a ring $R$, we write $\Sh_G(X,R)$ for the abelian category of $G$-equivariant sheaves of $R$-modules on $X$.
\end{definition}
If $G$ is endowed with the discrete topology, then the above definition coincides with the one given in the previous section. We will usually restrict ourselves to the simplest situation, where the action of $G$ on $X$ is free.
\begin{definition}
Let $G$ be a topological group and $X$ be a $G$-space. We say that $X$ is free if the quotient map $q: X \rightarrow X/G$ is a locally trivial $G$-torsor. In other words, there exists an open cover $\{U_i\}_{i\in I}$ of $X/G$ and $G$-equivariant isomorphisms $U_i \times G \cong q^{-1}(U_i)$.
\end{definition}
\begin{lemma}\label{lem_descent_torsor}
Let $X$ be a free $G$-space. Then the functor $q^*: \Sh(X/G) \rightarrow \Sh_G(X)$ is an equivalence of categories. An inverse is given by $q_*^G$ (defined by the same formula as in the case where $G$ is discrete, see before Lemma \ref{lem_descent_along_free_actions}). 
\end{lemma}
\begin{proof}
This is well known. It is a special case of descent along a torsor \cite[Theorem 4.46]{Vistoli}.
\end{proof}
It follows that if $X$ is a free $G$-space, then $\Sh_G(X)$ has enough injectives.
\subsubsection{Equivariant sheaves and smooth representations}\label{sec_smooth_eqvt}
\begin{lemma}\label{lem_smooth_compact}
Let $G$ be a topological group and $X$ be a $G$-space. Suppose $X$ is compact. Then for $\cF\in \Sh_G(X)$ the global sections of $\cF$ form a smooth $G$-representation. 
\end{lemma}
\begin{proof}
Let $s$ be a global section of $\cF$. We consider the two sections $\theta(p^*s)$ and $m^*s$ of $m^*\cF$ over $X\times G$. For $x \in X$, the stalks of $p^*s$ and $m^*\cF$ at $(x,e)$ are both given by $\cF_x$ and $\theta$ induces the identity map $p^*\cF_{(x,e)} \cong m^*\cF_{(x,e)}$. In particular, $\theta(p^*s)$ and $m^*s$ have the same image in the stalk at $(x,e)$, and hence coincide on some open neighbourhood $W_x \subset X\times G$ of $(x,e)$. We have $U_x \times G_x \subset W_x$ for $U_x$ an open neighbourhood of $x$ in $X$ and $G_x$ an open neighbourhood of $e$ in $G$.

Since $X$ is compact, we obtain a finite open cover $U_i$ of $X$ and open neighbourhoods $G_i$ of $e$ such that for all $(x,g) \in U_i \times G_i$, $\theta(p^*s)$ and $m^*s$ have the image in $m^\ast \cF(U_i \times G_i)$. We conclude that there is an open neighbourhood $H$ of $e$ in $G$ such that $\theta(p^*s)$ and $m^*s$ have the same restriction to $X\times H$. In other words, $s$ is fixed by an open neighbourhood of $e$ in $G$, and hence its stabilizer is an open subgroup of $G$.
\end{proof}
If $X = \text{pt} = \{ x \}$, then an object $\cF$ in $\Sh_G(X,R)$ gives rise to an $R\text{-module}$ $\cF_x$ equipped with an action of $G$. Lemma \ref{lem_smooth_compact} shows that this gives a functor $\Sh_G(X,R)\rightarrow \Modsm(G,R)$.
\begin{lemma}
The functor $\cF \rightarrow \cF_x$ induces an an equivalence of categories between $\Sh_G(\text{pt},R)$ and $\Modsm(G,R)$. 
\end{lemma}
\begin{proof}
The functor is clearly fully faithful, so we need to check essential surjectivity. For $M \in \Modsm(G,R)$ we set $\cF_M$ to be the sheaf on $\{x\}$ with sections $M$. We have $p^*\cF_M = m^*\cF_M$ and this is the sheaf of locally constant functions from $G$ to $M$. For $U \subset G$ an open subset we define $\theta: p^*\cF_M(U) \cong m^*\cF_M(U)$ by $\theta(f)(g) = gf(g)$, for $f$ a locally constant function $U \rightarrow M$ and $g \in U$. Since the action of $G$ on $M$ is smooth, $\theta(f)$ is again a locally constant function from $U$ to $M$. The cocycle condition for $\theta$ can be checked on stalks, where it amounts to the action of $G$ on $M$ being a group action.
\end{proof}
\begin{definition}\label{def_derived_free_action}
Let $X$ be a compact $G$-space. Denote the left exact functor obtained by taking global sections by \[\Gamma_X: \Sh_{G}(X,R) \rightarrow \Modsm(G,R).\]
If $X$ is a compact free $G$-space, we denote by $R\Gamma_X$ the right derived functor
\[R\Gamma_X: \mathbf{D}^+(\Sh_{G}(X,R)) \rightarrow \mathbf{D}^+_\text{sm}(G, R).\]
\end{definition}
\begin{lemma}\label{presinj}
Let $X$ be a compact $G$-space. The functor $\Gamma_X: \Sh_{G}(X,R) \rightarrow \Modsm(G,R)$ preserves injectives. 
\end{lemma}
\begin{proof}
The functor $\Gamma_X$ can be viewed as the direct image functor $\Sh_{G}(X,R) \rightarrow \Sh_{G}(pt,R)$. This has an exact left adjoint given by the inverse image functor, so $\Gamma_X$ preserves injectives.
\end{proof}

Now suppose that $G$ is a topological group, and $H$ is a locally profinite group. We suppose that $X$ is a compact $G \times H_\delta$-space, where $H_\delta$ indicates $H$ with the discrete topology. Let $U \subset H$ be an open compact subgroup such that $G\times U_\delta$ acts freely on $X$. We obtain a diagram of functors, commutative up to natural isomorphism, analogous to the diagram (\ref{eqn_hecke_functor_diagram}): 
\begin{equation}\label{eqn_completed_coh_hecke_functor_diagram}\begin{aligned}
\xymatrix{ \Sh_{G\times H_\delta}(X,R) \ar[r]^-{\Gamma_X} \ar[rd]_{\Res^{G\times H_\delta}_{G\times U_\delta}} & \Modsm(G\times H_\delta,R) \ar[r]^-{\Gamma_U} & \Modsm(G,\cH(H, U)\otimes_\bbZ R) \ar[r]^-{(\cdot)^\sim} & \Modsm(G,R) \\
& \Sh_{G\times U_\delta}(X,R) \ar[r]_{f_\ast^U} & \Sh_G(X/U,R) \ar[ur]_{\Gamma_{X/U}} } 
\end{aligned}
\end{equation}
Note that $X/U$ is a free $G$-space. We also have an equivalence $\Sh_{G\times H_\delta}(X,R) \cong \Sh_{H_\delta}(X/G,R)$, so this category has enough injectives. We can therefore define a right derived functor $R(\Gamma_U \circ \Gamma_X) : \mathbf{D}^+(\Sh_{G\times H_\delta}(X,R))\rightarrow \mathbf{D}^+_\text{sm}(G,\cH(H, U) \otimes_\bbZ R)$.
We obtain:
\begin{proposition}\label{prop_top_heck_action_on_cohomology}There is a canonical isomorphism in $\mathbf{D}^+_\text{sm}(G, R)$, for any $\cF \in \Sh_{G\times H_\delta}(X,R)$:
\[ R (\Gamma_U \circ \Gamma_X)(\cF)^\sim \cong R \Gamma_{X/U}(f_\ast^U \cF). \]
\end{proposition}
As in the discrete case, we will use the following consequence of the proposition: for any $\cF \in \Sh_{G\times H_\delta}(X,R)$, there is a canonical homomorphism 
\begin{equation}\label{eqn_hecke_homomorphism_topological_case} \cH(H,U) \rightarrow \End_{\mathbf{D}^+_\text{sm}(G, R)}(R\Gamma_{X/U}f_\ast^U \cF). 
\end{equation}
 For $\alpha \in H$, define an endomorphism $\theta(\alpha)$ of $R\Gamma_{X/U} f_*^U\cF$ as in the discrete case by pullback and pushforward. The same proof as before now yields the analogue of Lemma \ref{lem_hecke_action_as_expected}:
\begin{lemma}\label{lem_hecke_action_as_expected_top}
Let $\cF \in \Sh_{G\times H_\delta}(X,R)$. For $\alpha \in H$, the image of  $[U\alpha U]$ under the homomorphism (\ref{eqn_hecke_homomorphism_topological_case}) equals $\theta(\alpha)$.
\end{lemma}

\subsection{Completed cohomology}\label{cc}
We now recall some elements of the theory of completed cohomology. We begin by working in a general context as in \cite[1.1]{emcalsumm} and \cite[2.2]{Hill}. Let $G_0$ be a profinite group with a countable basis of neighbourhoods of the identity given by normal open subgroups \[\cdots \subset G_n \subset \cdots \subset G_1 \subset G_0.\] 
Suppose given a tower of compact topological spaces \[\cdots \rightarrow X_n \rightarrow \cdots \rightarrow X_1 \rightarrow X_0,\] each equipped with an action of $G_0$. We moreover suppose that:
\begin{enumerate}
\item The maps $X_{n+1} \rightarrow X_n$ are $G_0$-equivariant.
\item $G_n$ acts trivially on $X_n$ and $X_n$ is a (locally trivial) $G_0/G_n$-torsor over $X_0$.
\end{enumerate}
Finally, we assume that $X_0$ admits an open covering by contractible subsets (for example, $X_0$ is locally contractible). In the above situation, we define a topological space 
\[X = \plim_n X_n,\] endowed with the projective limit topology. $X$ is a compact $G_0$-space. We write $\pi_n$ for the maps $X_n \rightarrow X_0$ and $\pi$ for the map $X \rightarrow X_0$.
\begin{lemma} The space $X$ is a free $G_0$-space and the natural map $X/G_0 \rightarrow X_0$ is an isomorphism. 
\end{lemma}
\begin{proof}It is clear that the canonical map $X \rightarrow X_0$ identifies $X_0$ with the quotient $X/G_0$. To show that the $G_0$ action is free, we must show that the quotient map $X \rightarrow X_0$ is a locally trivial $G_0$-torsor. Let $U$ be a contractible open subset of $X_0$. For each $n$, the fibre product $X_{n}|_U:=X_n \times_{X_0} U$ is a torsor over $U$ for the finite group $G_0/G_n$. We therefore have an isomorphism of $G_0/G_n$-torsors over $U$: \[\tau_n: X_{n}|_U \cong U \times (G_0/G_n).\] We are going to construct an isomorphism of $G_0$-spaces \[X|_U := \varprojlim_n X_{n}|_U \cong U \times G_0\]
 by modifying the isomorphisms $\tau_n$. Suppose we have isomorphisms \[\tau'_i: X_{i}|_U \cong U \times (G_0/G_i)\] for $0 \le i \le n-1$, which, for $1\le i \le n-1$, send the transition maps $X_{i}\rightarrow X_{i-1}$ to the obvious projection $U \times (G_0/G_i) \rightarrow U \times (G_0/G_{i-1})$. 

We consider the surjective map of $G_0$-spaces $U \times (G_0/G_n) \rightarrow U\times (G_0/G_{n-1})$ induced by the isomorphisms $\tau_{n-1}', \tau_n$ and the transition map $X_n \rightarrow X_{n-1}$. This map sends $(u,x)$ to $(u,\alpha(x))$, where $\alpha: G_0/G_n \rightarrow G_0/G_{n-1}$ is a map of $G_0$-sets. This map is therefore determined by $\alpha(1) \in G_0/G_{n-1}$. We define $\tau_n'$ by $\tau_n'(x) = \tau_n(x)g$, where $g$ is any representative of $\alpha(1)$ in $G_0$. 

We set $\tau_0'$ equal to the identity, and by induction we have constructed $\tau_n'$ as above for all $n$. Now taking the projective limit gives the desired trivialisation of $X|_U$. 
\end{proof}
\begin{lemma}\label{lemmalocfin}Let $R$ be a ring. 
The category $\Modsm(G_0,R)$ has a generator and exact inductive limits. In particular, $\Modsm(G_0,R)$ has enough injectives. 
\end{lemma}
\begin{proof}
A generator is given by $X = \bigoplus_{n \ge 0}\Ind^{G_0}_{G_n}R$, since $\Hom(X,M) = \prod_{n \ge 0}M^{G_n}$ which is non-zero for all $M \in \Modsm(G_0,R)$. It is clear that inductive limits exist in $\Modsm(G_0,R)$, and they are exact by \cite[Proposition I.6b]{Gabriel-AbCat}. 
\end{proof}
Given $\cF_0 \in \Sh(X_0)$,  we set $\cF_n = \pi_n^*\cF_0 \in \Sh_{G_0/G_n}(X_n)$ and set $\cF = \pi^*\cF_0 \in \Sh_{G_0}(X)$.
\begin{lemma}\label{lem_smooth_cc}
The natural maps $\Gamma_{X_n}(\mathcal{F}_n) \rightarrow \Gamma_X(\mathcal{F})$ induce an isomorphism \[\varinjlim_n \Gamma_{X_n}(\mathcal{F}_n) \cong \Gamma_X(\mathcal{F}).\]
\end{lemma}
\begin{proof}
The natural maps $\Gamma_{X_n}(\mathcal{F}_n) \rightarrow \Gamma_X(\mathcal{F})$ identify $\Gamma_{X_n}(\mathcal{F}_n)$ with $\Gamma_X(\mathcal{F})^{G_n}$, by Lemma \ref{lem_descent_torsor}. By Lemma \ref{lem_smooth_compact}, $\Gamma_X(\mathcal{F})$ is smooth, which gives the desired result.
\end{proof}

\begin{lemma}\label{altdes}The functor $R\Gamma_X : \mathbf{D}^+(\Sh_{G_0}(X,R)) \rightarrow \mathbf{D}^+_\text{sm}(G_0,R)$, when composed with the equivalence of triangulated categories $\mathbf{D}^+(\Sh(X_0,R))\cong \mathbf{D}^+(\Sh_{G_0}(X,R))$, is the right derived functor of the functor $\Sh(X_0,R) \rightarrow \Modsm(G_0,R)$ given by 
\[\mathcal{F}_0 \mapsto \varinjlim_n \Gamma_{X_n}(\mathcal{F}_n).\]
\end{lemma}
\begin{proof}
This follows from Lemma \ref{lem_smooth_cc}.
\end{proof}
\begin{lemma}\label{identcoho}
There are canonical isomorphisms 
\[H^i(R\Gamma_X(\mathcal{F})) \cong \varinjlim_n H^i(X_n,\mathcal{F}_n).\]
\end{lemma}
\begin{proof}
The previous lemma identifies $R\Gamma_X$ with the derived functor of \[\mathcal{F}_0 \mapsto \varinjlim_n \Gamma_{X_n}(\mathcal{F}_n)\] from $\Sh(X_0)$ to $\Modsm(G_0)$. Taking an injective resolution $\mathcal{I}_0^\bullet$ of $\mathcal{F}_0$ we get injective resolutions $\mathcal{I}_n^\bullet=\pi_n^*(\mathcal{I}_0^\bullet)$ of $\mathcal{F}_n$ for each $n$ ($\pi_{n,!}$ is exact, so $\pi_n^*$ preserves injectives). Now $H^i(R\Gamma_X(\mathcal{F}))$ is (by definition) given by $H^i(\varinjlim_n\Gamma_{X_n}(\mathcal{I}_n^\bullet)) = \varinjlim_nH^i(\Gamma_{X_n}(\mathcal{I}_n^\bullet)) = \varinjlim_n H^i(X_n,\mathcal{F}_n).$ 
\end{proof}

\begin{definition}
For a ring $R$, denote by \[R\Gamma_{G_n} : \mathbf{D}^+_\text{sm}(G_0,R) \rightarrow \mathbf{D}^+(G_0/G_n, R) \] the right derived functor of taking $G_n$-invariants.
\end{definition}
Since $G_0$ is compact, the derived functors $R\Gamma_{G_n}$ may be computed using standard resolutions. For $M \in \Modsm(G_0,R)$ and $r \ge 0$ we denote by $X^r(M)$ the object of $\Modsm(G_0,R)$ given by locally constant maps from $G_0^{r+1}$ to $M$. The action of $G_0$ is given by $(\sigma f)(\sigma_0,\ldots,\sigma_r) = \sigma f(\sigma^{-1}\sigma_0,\ldots,\sigma^{-1}\sigma_r)$. As in \cite[I.2]{NSW}, we define a complex $X^\bullet(M)$: there are maps $d_i: G^{r}\rightarrow G^{r-1}$ given by omitting the $i$th term, which induce maps $d_i^*:X^{r-1}\rightarrow X^{r}$, and the maps in the complex are given by 
\[\delta = \sum_{i=0}^r (-1)^id_i^*:X^{r-1}\rightarrow X^r.\] 
The map $M \rightarrow X^0(M)$ given by sending $m$ to the constant function with value $m$ induces a quasi-isomorphism $M \rightarrow X^\bullet(M)$ in $\Modsm(G_0,R)$ (\cite[Proposition 1.2.1]{NSW}).

\begin{lemma}
For $i > 0$ and $r \ge 0$ we have $R^i\Gamma_{G_n} X^r(M) = 0$. As a consequence, the natural map $X^\bullet(M)^{G_n} \rightarrow R\Gamma_{G_n} X^\bullet(M)$ is a quasi-isomorphism and we have an isomorphism $X^\bullet(M)^{G_n} \cong R\Gamma_{G_n} M$ in $\mathbf{D}^+(G_0/G_n,R)$.
\end{lemma}
\begin{proof}
First we recall some standard results in continuous group cohomology, see for example \cite[Chapter I]{NSW}. The functors $M \mapsto \Gamma^i(M):= H^i(X^\bullet(M)^{G_n})$ define a cohomological $\delta$-functor from $\Modsm(G_0,R)$ to $\Mod(G_0/G_n, R)$, and we have $\Gamma^i(X^r(M)) = 0$ for all $i > 0$ and $r \ge 0$ (the modules $X^r(M)$ are induced). The map $M \rightarrow X^0(M)$ induces an isomorphism $\Gamma_{G_n} M \cong \Gamma^0(M)$. Therefore the functors $\Gamma^i$ form a universal $\delta$-functor, and hence $R^i\Gamma_{G_n} X^r(M) = \Gamma^i(X^r(M)) = 0$.
\end{proof}

\begin{lemma}\label{lem_flat_coeff_change}
Let $R \rightarrow R'$ be a flat ring map and denote the extension of scalars functor $\Modsm(G_0,R) \rightarrow \Modsm(G_0,R')$ by $(-)\otimes_R {R'}$. Then there is a natural isomorphism of functors $\mathbf{D}^+_\text{sm}(G_0,R)\rightarrow \mathbf{D}^+(G_0/G_n, R')$:
\[R\Gamma_{G_n}(-) \otimes_R R' \cong R\Gamma_{G_n}((-)\otimes_R R')\]
\end{lemma}
\begin{proof}
For $M \in \Modsm(G_0,R)$ we have a natural isomorphism $X^\bullet(M)\otimes_R R' \cong X^\bullet(M\otimes_R R')$ (each function in $X^r(M\otimes_R R')$ has a finite set of values, so it is a finite $R'$-linear combination of $M$-valued functions). Since we can compute $R\Gamma(G_n,-)$ using resolutions by the acyclic objects $X^\bullet(M)$, and $(-)\otimes_R R'$ sends these acyclic resolutions to acyclic resolutions, we obtain a natural isomorphism between $R\Gamma_{G_n}((-)\otimes_R R')$ and the derived functor of \[M\mapsto (M\otimes_R R')^{G_n}.\] 
On the other hand, since $(-)\otimes_R R'$ is exact, $R\Gamma_{G_n}(-)\otimes_R R'$ is naturally isomorphic to the derived functor of \[M \mapsto M^{G_n}\otimes_R R'.\]
To prove the lemma, it suffices to show that the natural map $M^{G_n}\otimes_R R' \rightarrow (M\otimes_R R')^{G_n}$ is an isomorphism. Since the action of $G_n$ on $M$ is smooth, it suffices to check that this map is an isomorphism for $G_n$ a finite group. This follows from the isomorphism (\ref{homtensor}) in the proof of Lemma \ref{lem_derived_homtensor}, as $M^{G_n} = \Hom_{R[G_n]}(R,M)$ and $R$ is a finitely presented $R[G_n]$-module.
\end{proof}

\begin{lemma}\label{lem_invts_quo}
Let $R$ be a ring, and let $\mathcal{F}\in \Sh_{G_0}(X,R)$ and $n \ge 0$. Denote by $\mathcal{F}_n \in \Sh_{G_0/G_n}(X_n,R)$ the sheaf on $X_n$ obtained by descent from $\cF$. There are natural isomorphisms in $\mathbf{D}^+(G_0/G_n,R)$: \[R\Gamma_{G_n} R\Gamma_X\mathcal{F} \cong R\Gamma_{X_n}\mathcal{F}_n \] extending the natural isomorphism $\Gamma_{G_n} \Gamma_X\mathcal{F} \cong \Gamma_{X_n}\mathcal{F}_n$.
\end{lemma}
\begin{proof}
This follows from Lemma \ref{presinj}.
\end{proof}

\subsubsection{Completed cohomology without taking a limit of spaces}\label{ccnoprojlim}
We now present a variant of the constructions of \S \ref{cc} which works just with sheaves at `finite levels' $X_n$ instead of passing to the limit $X$. This variant will then be applicable in a more general situation: for example, when the group actions are not free (we will later work with minimal compactifications as well as Borel--Serre compactifications), and the spaces $X_n$ are algebraic varieties or even adic spaces.

In this section a `space' means one of the following: a topological space, an adic space over a complete and algebraically closed extension of $\bbQ_p$, or an algebraic variety over an algebraically closed field of characteristic $0$. We again let $G_0$ be a profinite group with a countable basis of neighbourhoods of the identity given by normal open subgroups 
\[\cdots \subset G_n \subset \cdots \subset G_1 \subset G_0.\] 
Suppose given a tower of spaces 
\[\cdots \rightarrow X_n \rightarrow \cdots \rightarrow X_1 \rightarrow X_0,\] 
with each $X_n$ equipped with an action of the finite group $G_0/G_n$ and the transition maps equivariant with respect to these actions. Contrary to the last section, we do not assume that these group actions are free.

We now consider categories $S_n = \Sh_{G_0/G_n}(X_n,R)$, where we take equivariant sheaves on the topological space $X_n$ or the \'{e}tale site of the algebraic variety or adic space $X_n$ as appropriate. (See \cite{Huber} for the definition of the \'etale site of an adic space.) Note that the transition maps $X_m \rightarrow X_n$ for $m \ge n$ give pairs of functors $\pi_{m,n *}^{G_n/G_m}: S_m \rightarrow S_n$ and $\pi_{m,n}^*: S_n \rightarrow S_m$. By the first part of Lemma \ref{lem_descent_along_free_actions} and its analogue for the other sites, $\pi_{m,n *}^{G_n/G_m}$ is a right adjoint to $\pi_{m,n}^*$. 
\begin{definition}
We define $S$ to be the category whose objects are collections of sheaves $(\cF_n)_{n \geq 0}$ with $\cF_n \in S_n$, equipped with morphisms $\theta_{m,n}:\pi_{m,n}^*\cF_n \rightarrow \cF_m$ for $m \ge n$, which satisfy a cocycle relation: for $i\ge j \ge k$ we have \[\theta_{i,k} = \theta_{i,j}\circ \pi_{i,j}^*\theta_{j,k}.\]
The morphisms in $S$ are given by families of morphisms $f_n: \cF_n \rightarrow \cG_n$, such that for $m \ge n$ \[f_m \circ\theta_{m,n} = \theta_{m,n}\circ \pi_{m,n}^*(f_n).\]
\end{definition}
The category $S$ depends on the ring $R$ of coefficients, although we do not include this in the notation. We observe that $S$ is the category of $R$-module objects in the total topos of a fibred topos (see \cite[D\'{e}finition 7.4.1]{SGA42}), although this is not made use of in the sequel. In particular, we can deduce immediately that $S$ has enough injectives, but we will show this directly (Lemma \ref{lem_S_enough_inj}).
\begin{definition}
Denote by $\iota_n^*$ the functor $S \rightarrow S_n$ given by $(\cF_m)_{m\ge 0}\mapsto \cF_n$. Denote by $\iota_{n,!}$ the functor $S_n \rightarrow S$ given by \[\cF \mapsto (0,\ldots,0,\cF,\ldots,\pi_{m,n}^*\cF,\ldots).\]
Denote by $\iota_{n,*}$ the functor $S_n \rightarrow S$ given by \[\cF \mapsto (\pi_{n,0,*}^{G_0/G_n}\cF,\ldots,\pi_{n,m,*}^{G_m/G_n}\cF,\ldots,\cF,0,0,\ldots).\]
\end{definition}
\begin{lemma}
The functor $\iota_{n,!}$ is a left adjoint to $\iota_n^*$, and $\iota_{n,*}$ is a right adjoint to $\iota_n^*$. The functor $\iota_{n,!}$ is exact.
\end{lemma}
\begin{proof}
First we check the adjointness properties. Suppose we have $\cF \in S_n$ and $\cG_\bullet \in S$. Then an element of $\Hom_S(\iota_{n,!}\cF , \cG)$ is given by a collection of maps $f_m$ in $\Hom_{S_m}(\pi_{m,n}^*\cF_n,\cG_m)$ for $m \ge n$ such that \[f_{m} = \theta_{m,n}\circ \pi_{m,n}^*(f_n).\] So we see immediately that everything is determined by $f_n$, and $\Hom_S(\iota_{n,!}\cF , \cG) = \Hom_{S_n}(\cF,\cG_n) = \Hom_{S_n}(\cF,\iota_n^*\cG)$ as required.

Next, suppose we have $\cF_\bullet \in S$ and $\cG \in S_n$. Given a map $f_n: \cF_n \rightarrow \cG$ there is a natural way to produce a map $\cF_\bullet \rightarrow \iota_{n,*}\cG$: for $m \le n$ we let $f_m: \cF_m \rightarrow \pi_{n,m,*}^{G_m/G_n}\cG$ be the map corresponding by adjunction to $\pi_{n,m}^*\cF_m \overset{\theta_{n,m}}{\rightarrow} \cF_n \overset{f_n}{\rightarrow}\cG$. This gives a natural map $\Hom_{S_n}(\iota_n^*\cF, \cG) \rightarrow \Hom_S(\cF,\iota_{n,*}\cG)$. To prove that this is a  bijection, we must show that this choice of $f_m$ is the unique map making the following diagram commute:
\[ \xymatrix{ \pi_{n,m}^*\cF_m \ar[r]^{\pi_{n,m}^*f_m} \ar[d]^{\theta_{n,m}} & \pi_{n,m}^*\pi^{G_m/G_n}_{n,m *}\cG \ar[d] \\
\cF_n \ar[r]^{f_n}  & \cG } \]
The right hand vertical arrow here is the counit of the adjunction. So the composition of the top horizontal and right vertical maps is identified with $f_m$ under the bijection $\Hom_{S_n}(\pi_{n,m}^*\cF_m,\cG) = \Hom_{S_m}(\cF_m,\pi_{n,m,*}^{G_m/G_n}\cG)$. This forces $f_m$ to be the map we have defined above.

This implies that $\iota_{n}^*$ is exact, so it preserves kernels and images. In particular, kernels and images of maps in $S$ are given by componentwise kernels and images, so one can check exactness of complexes in $S$ componentwise. Since $\pi_{m,n}^*$ is an exact functor, $\iota_{n,!}$ is exact.
\end{proof}
\begin{lemma}
For $\cF\in S$ an injective object the equivariant sheaf $\iota_n^*\cF \in S_n$ and the underlying sheaf in $\Sh(X_n,R)$ are injective. The functor $\iota_{n,*}$ also preserves injectives.
\end{lemma}
\begin{proof}
The functors $\iota_n^*$ and $\iota_{n,*}$ have exact left adjoints, as does the forgetful functor from $S_n$ to $\Sh(X_n,R)$ (by the same construction as in Lemma \ref{lem_adjunctions}). It follows that all of these functors preserve injectives.
\end{proof}
\begin{lemma}\label{lem_S_enough_inj}
The category $S$ has enough injectives.
\end{lemma}
\begin{proof}
Let $\cF \in S$. For each $n$, we pick a monomorphism $\iota_n^*\cF \hookrightarrow \cI_n$ to an injective object of $S_n$. By adjointness we have maps $\cF \rightarrow \iota_{n,*}\cI_n$ for each $n$, so we obtain a map $\cF\rightarrow \cI := \prod_{n\ge 0}\iota_{n,*}\cI_n$. Since products of injectives are injective, $\cI$ is injective. The map $\cF \rightarrow \cI$ is monic, as this can be checked componentwise. So $S$ has enough injectives.
\end{proof}

\begin{definition}\label{def_fibered_cc}
We denote by $\widetilde{\Gamma}$ the functor $S \rightarrow \Modsm(G_0,R)$ given by \[\cF \mapsto \varinjlim_n \Gamma_{X_n}(\cF_n).\]
We denote by $R\widetilde{\Gamma} : \mathbf{D}^+(S) \to \mathbf{D}_\text{sm}^+(G_0, R)$ the right derived functor of $\widetilde{\Gamma}$.
\end{definition}
\begin{lemma}\label{lem_completed_coho}
There are natural isomorphisms \[R^i\widetilde{\Gamma}(\cF) \cong \varinjlim_n H^i(X_n,\cF_n).\]
\end{lemma}
\begin{proof}
We take an injective resolution $\cF\rightarrow\cI^\bullet$. For each $n$, $\cF_n \rightarrow \cI^\bullet_n$ is an injective resolution.
Since direct limits are exact in $\Modsm(G_0,R)$, we have \[R^i\widetilde{\Gamma}(\cF) = H^i(\varinjlim_{n}\Gamma_{X_n}(\cI^\bullet_n)) = \varinjlim_nH^i(\Gamma_{X_n}(\cI^\bullet_n)) = \varinjlim_nH^i(X_n,\cF_n).\]
\end{proof}

\begin{lemma}\label{lem_fibered_free_equiv}
Suppose the $X_n$ are compact locally contractible topological spaces and the spaces $X_n \rightarrow X_0$ are $G_0/G_n$-torsors. In other words, we suppose that the formalism of the previous section applies to $X_\bullet$. We regard $R\Gamma_X$ as a functor on $S_0$, using Lemma \ref{altdes}. Then there is a natural isomorphism
\[R\Gamma_X \cong R\widetilde{\Gamma}\circ\iota_{0,!}.\]
\end{lemma}
\begin{proof}
We have a natural isomorphism of functors $\Gamma_X \cong \widetilde{\Gamma}\circ\iota_{0,!}$. By Lemma \ref{lem_universal_property_of_derived_functor} we obtain a natural transformation \[R\Gamma_X \rightarrow R\widetilde{\Gamma}\circ\iota_{0,!}.\] Tracing through the constructions and applying Lemma \ref{lem_completed_coho}, we see that this is an isomorphism.
\end{proof}
\begin{lemma}
Suppose that the $X_n$ are compact Hausdorff topological spaces, let $j_0 : Y_0 \hookrightarrow X_0$ be an open subspace, and let $j_n : Y_n \hookrightarrow X_n$ be the pullback of $j_0$ for each $n \geq 0$. Let $S_Y$ denote the analogue of the category $S$ for the tower $(Y_n)_{n \geq 0}$, and let $j_! : S_Y \to S$ denote the exact functor induced by the functors $j_{n, !}$. Let $\widetilde{\Gamma}_c: S_Y \rightarrow \Modsm(G_0,R)$ be the functor
\[ \widetilde{\Gamma}_c(\cF) = \varinjlim_n \Gamma_c(Y_n,\cF_n) = \widetilde{\Gamma}\circ j_!\cF. \]
 Then there is a natural isomorphism of functors
\[ R\widetilde{\Gamma}_c \cong R\widetilde{\Gamma}\circ j_! : \mathbf{D}^+(S_Y) \to \mathbf{D}_\text{sm}^+(G_0, R). \]
\end{lemma}
\begin{proof}
It suffices to check that for $\cI$ an injective object of $S_Y$ we have $R^i\widetilde{\Gamma} j_!\cI = 0$ for all $i > 0$. We have $R^i\widetilde{\Gamma} j_!\cI = \varinjlim_n H^i(X_n,j_{n,!}\cI_n)$, so the vanishing follows the fact that $j_{n,!}$ takes injective sheaves to soft sheaves, which are $\Gamma$-acyclic (see \cite[Proposition III.7.2]{Ive86}).
\end{proof}
\subsubsection{Comparing topologies}\label{sec_compare_fibered_cc}
Suppose given a complete algebraically closed extension $C$ of $\bbQ_p$ with compatible embeddings $\overline{\bbQ}\subset \bbC$ and $\overline{\bbQ}\subset C$, together with a tower $X_n^{\mathrm{alg}}$ of proper schemes over $\overline{\bbQ}$. We set $X_n^{\mathrm{top}} = X_n^{\mathrm{alg}}(\bbC)$ with the usual topology, set $X_n^{\mathrm{cl}}$ to be the `local isomorphisms' site on $X_n^{\mathrm{alg}}(\bbC)$ as defined in \cite[XI 4]{SGA43}, and set $X_n^{\mathrm{ad}}$ to be the adic space associated to $X_{n,C}^{\mathrm{alg}}$.

If $? \in \{ \text{top}, \text{cl}, \text{alg}, \text{ad} \}$, then we denote by $\widetilde{\Gamma}^{?}$ and $R\widetilde{\Gamma}^{?}$ the functors obtained by applying the formalism of \S \ref{ccnoprojlim} to the tower of spaces $(X_n^{?})_{n \geq 0}$ (with the appropriate site). Recall that there is a morphism of sites 
\[\epsilon_n: X_{n}^\mathrm{cl} \rightarrow (X_{n,\bbC}^\mathrm{alg})_{\text{\'et}}\rightarrow (X_{n}^\mathrm{alg})_{\text{\'et}}\] 
together with an inclusion of sites $X_{n}^\mathrm{cl} \rightarrow X_n^{\mathrm{top}}$ inducing an equivalence of topoi. This induces an exact functor  $\epsilon^*$ from $S^\mathrm{top}$ to $S^\mathrm{alg}$. As a result we obtain a base change natural transformation $R\widetilde{\Gamma}^\mathrm{alg} \rightarrow R\widetilde{\Gamma}^\mathrm{top}\epsilon^*$, for example by applying Lemma \ref{lem_universal_property_of_derived_functor}. This natural transformation is an isomorphism by Lemma \ref{lem_completed_coho} and the usual comparison theorem for cohomology \cite[XVII Corollaire 5.3.5]{SGA43}.

Similarly, we obtain a natural isomorphism $R\widetilde{\Gamma}^\mathrm{alg} \rightarrow R\widetilde{\Gamma}^\mathrm{ad}\epsilon^*$, where $\epsilon^*$ is induced by the morphisms of sites $\epsilon_n: (X_{n}^\mathrm{ad})_{\text{\'et}} \rightarrow (X_{n,C}^\mathrm{alg})_{\text{\'et}}\rightarrow (X_{n}^\mathrm{alg})_{\text{\'et}}$ (we use (3.2.8) and Theorem 3.7.2 of \cite{Huber}).

\subsection{Almost smooth representations}\label{sec_almost_smooth}Fix a complete and algebraically closed extension $C$ of $\bbQ_p$. Let $\cO_C \subset C$ denote the ring of integers and $\ffrm \subset \cO_C$ the maximal ideal. Let $\pi \in \ffrm - \{ 0 \}$. Fix $N \geq 1$ and let $\bbV = \cO_C/(\pi^N)$.   Recall \cite{GabberRamero}, \cite[\S4]{perfectoid} that the almost context $(\cO_C,\ffrm)$ allows us to define a category of almost $\cO_C$-modules, or $\cO_C^a$-modules, that we denote $\Mod(\cO_C^a)$. This is obtained by localising the category $\Mod(\cO_C)$ of $\cO_C$-modules at the Serre subcategory comprising modules which are killed by $\ffrm$.

One then defines $\cO_C^a$-algebras, and in particular we have an $\cO_C^a$-algebra $\bbV^a$, which we denote by $A$. There is an exact localisation functor $(-)^a$ from $\Mod(\bbV)$ to $\Mod(A)$. This functor has a right adjoint $M \mapsto M_\ast = \Hom_A(A, M)$, the functor of almost elements; for any $M \in \Mod(A)$, the adjunction morphism $(M_\ast)^a \to M$ is an isomorphism.  The functor  $(-)^a$ also has an exact left adjoint $M \mapsto M_! = \ffrm \otimes_{\cO_C}M_\ast$ (see \cite[2.2.21, 2.2.23]{GabberRamero}), and the adjunction morphism $M\mapsto (M_!)^a$ is again an isomorphism.

Let $G_0$ be a profinite group equipped with a nested sequence $G_0 \supset G_1 \supset \dots$ of open compact subgroups.
\begin{definition}\label{def_almost_smooth}
We say $X \in \Modsm(G_0,\bbV)$ is almost zero if $\ffrm X = 0$. The full subcategory of almost zero objects in $\Modsm(G_0,\bbV)$ is a Serre subcategory and we denote by $\Modsm(G_0,A)$ the abelian category obtained as the quotient of $\Modsm(G_0,\bbV)$ by this Serre subcategory.
\end{definition}
\begin{lemma}
The functor $(-)_!$ induces an exact left adjoint to the localisation functor $\Modsm(G_0,\bbV)\rightarrow \Modsm(G_0,A)$, which we also denote by $(-)_!$. The unit of the adjunction $M \rightarrow (M_!)^a$ is a natural isomorphism from the identity functor to the composition $((-)_!)^a$.
\end{lemma}
\begin{proof}
We note that $X \mapsto \ffrm\otimes_{\cO_C} X$ defines an exact functor $\Modsm(G_0,\bbV) \rightarrow \Modsm(G_0,\bbV)$ which is zero on almost zero objects. Therefore we obtain an exact functor $(-)_!$ from $\Modsm(G_0,A)$ to $\Modsm(G_0,\bbV)$, by the universal property of a quotient by a Serre subcategory. This is seen to be left adjoint to the localisation functor by the same argument used to deduce \cite[(2.2.4)]{GabberRamero}: namely, describe $\Hom$ groups in the localisation $\Modsm(G_0,A)$ by calculus of fractions and observe that $\ffrm\otimes_{\cO_C} X \rightarrow X$ is initial in the category of almost isomorphisms to $X \in \Modsm(G_0,\bbV)$. In particular, for $X, Y \in \Modsm(G_0,\bbV)$ we obtain a natural isomorphism of $\bbV$-modules 
\[\Hom_{\Modsm(G_0,A)}(X^a,Y^a) = \Hom_{\Modsm(G_0,\bbV)}(\ffrm\otimes_{\cO_C}X,Y).\]
Finally, for $M = X^a$, the argument mentioned above shows that the unit of the adjunction is given by the map $X^a \rightarrow (\ffrm\otimes_{\cO_C}X)^a$ which is the inverse of the almost isomorphism $\ffrm\otimes_{\cO_C}X \rightarrow X$. This implies that the unit of the adjunction is a natural isomorphism, as required.
\end{proof}
\begin{remark}
The functor $(-)_\ast$ induces a functor from $\Modsm(G_0,A)$ to $\Mod(G_0,\bbV)$ but the resulting $\bbV[G_0]$-module may not always be smooth. To define the functor, we note that for $X \in \Modsm(G_0,\bbV)$, $\Hom_{A}(A,X^a)=\Hom_\bbV(\ffrm\otimes_{\cO_C}\bbV,X)$ is naturally equipped with an action of $G_0$ (which may not be smooth). Taking smooth vectors gives a functor $(-)_\ast^\text{sm}: \Modsm(G_0,A) \rightarrow \Modsm(G_0,\bbV)$, which can be checked to be a right adjoint to the localisation functor. However, we will not use this right adjoint in the sequel. 
\end{remark}

\begin{definition} We denote by $\mathbf{D}^+_\text{sm}(G_0,\bbV)$ (respectively $\mathbf{D}^+_\text{sm}(G_0,A)$) the bounded-below derived categories of $\Modsm(G_0,\bbV)$ (resp. $\Modsm(G_0,A)$).
\end{definition}
The (exact) localisation functor $(-)^a$ induces a functor $\mathbf{D}^+_\text{sm}(G_0,\bbV)\rightarrow \mathbf{D}^+_\text{sm}(G_0,A)$.

\begin{lemma}\label{lem_almost_presinj}
The localisation functor $\Modsm(G_0,\bbV)\rightarrow \Modsm(G_0,A)$ preserves injectives.
\end{lemma}
\begin{proof}
This follows from exactness of the left adjoint $(-)_!$.
\end{proof}
\begin{lemma}
The category $\Modsm(G_0,A)$ has enough injectives.
\end{lemma}
\begin{proof}
For $M \in \Modsm(G_0,A)$ we have a monomorphism $M_! \hookrightarrow I$, with $I$ an injective object of $\Modsm(G_0,\bbV)$, since $\Modsm(G_0,\bbV)$ has enough injectives (Lemma \ref{lemmalocfin}). Then applying the localisation functor gives a monomorphism $M \cong M_!^a \hookrightarrow I^a$. By Lemma \ref{lem_almost_presinj}, $I^a$ is injective. 
\end{proof}

\begin{definition}
For $n \ge 0$ we denote by \[\Gamma_{G_n} : \Modsm(G_0,A) \rightarrow \Modsm(G_0/G_n,A)\] the functor given by
\[M \mapsto \Hom_{\Modsm(G_n, A)}((\mathbf{1}_\bbV)^a,M)^a.\]
Here, we note that if $M= X^a$, $\Hom_{\Modsm(G_n, A)}((\mathbf{1}_\bbV)^a,M) = \Hom_{\Modsm(G_n, \bbV)}(\ffrm\otimes_{\cO_C}\bbV,X)$ which we view as an object of $\Mod(G_0/G_n,\bbV)=\Modsm(G_0/G_n,\bbV)$, and then apply $(-)^a$ to get something in $\Modsm(G_0/G_n,A)$.
\end{definition}
\begin{lemma}\label{lem_almost_hom}
For $M, N \in \Modsm(G_0,\bbV)$, the natural map of $\bbV$-modules \[\Hom_{\Modsm(G_0,\bbV)}(M,N) \rightarrow \Hom_{\Modsm(G_0,A)}(M^a,N^a)\] induces an isomorphism of $A$-modules \[\Hom_{\Modsm(G_0,\bbV)}(M,N)^a \rightarrow \Hom_{\Modsm(G_0,A)}(M^a,N^a)^a.\]
In particular, for $X \in \Modsm(G_0,\bbV)$ the map \[\Hom_{\Modsm(G_n,\bbV)}(\mathbf{1}_\bbV,X|_{G_n}) \rightarrow \Hom_{\Modsm(G_n,A)}((\mathbf{1}_\bbV)^a,X|_{G_n}^a)\] induces an isomorphism $(\Gamma_{G_n} X)^a\cong\Gamma_{G_n}( X^a)$.
Moreover, for $M, N \in \mathbf{D}^+_\text{sm}(G_0,\bbV)$ we similarly have a natural isomorphism of $A$-modules \[\Hom_{\mathbf{D}^+_\text{sm}(G_0,\bbV)}(M,N)^a \rightarrow \Hom_{\mathbf{D}^+_\text{sm}(G_0,A)}(M^a,N^a)^a.\]
\end{lemma}
\begin{proof}
We have $\Hom_{\Modsm(G_0,A)}(M^a,N^a) = \Hom_{\Modsm(G_0,\bbV)}(\ffrm \otimes_{\cO_C}M,N)$, by Lemma \ref{lem_almost_hom}, and the natural (multiplication) map $\ffrm \otimes_{\cO_C}M \rightarrow M$ is an almost isomorphism (the kernel and cokernel are killed by $\ffrm$). The induced map \[\Hom_{\Modsm(G_0,\bbV)}(M,N) \rightarrow \Hom_{\Modsm(G_0,\bbV)}(\ffrm \otimes_{\cO_C}M,N)\] is therefore an almost isomorphism. 

To check the `moreover' statement, we work with the homotopy categories $\mathbf{K}^+ = \mathbf{K}^+(\text{Inj}_\text{sm}(G_0,\bbV))$ and $\mathbf{K}^{+,a}=\mathbf{K}^+(\text{Inj}_\text{sm}(G_0,A))$, where $\text{Inj}$ denotes the full subcategory of injective objects in $\Modsm(G_0,\bbV)$ and $\Modsm(G_0,A)$. These categories are equivalent to the bounded below derived categories $\mathbf{D}^+_\text{sm}(G_0,\bbV)$ and $\mathbf{D}^+_\text{sm}(G_0,A)$ (\cite[Theorem 10.4.8]{Wei94}). We denote by $\mathbf{Kom}^+$ and $\mathbf{Kom}^{+,a}$ the categories of bounded below complexes of injectives. For $M, N \in \text{ob}(\mathbf{K}^+) = \text{ob}(\mathbf{Kom}^{+})$ we set $H^+ = \prod_{i \in \bbZ}\Hom_{\Modsm(G_0,\bbV)}(M^i,N^{i-1})$ and set $H^{+,a} = \prod_{i \in \bbZ}\Hom_{\Modsm(G_0,A)}((M^i)^a,(N^{i-1})^a)$. Since $(-)^a$ has a left adjoint, it commutes with direct products, and so the natural map of $\bbV$-modules $H^+ \rightarrow H^{+,a}$ induces an isomorphism of $A$-modules $(H^+)^a \cong (H^{+,a})^a$ (by the first part of the lemma). Now consider the commutative diagram of $\bbV$-modules, with exact rows \[ \xymatrix{ H^+ \ar[d] \ar[r] & \Hom_{\mathbf{Kom}^+}(M,N) \ar[d] \ar[r] &  \Hom_{\mathbf{K}^+}(M,N) \ar[d]\ar[r] & 0\\
H^{+,a} \ar[r] & \Hom_{\mathbf{Kom}^{+,a}}(M^a,N^a) \ar[r] &  \Hom_{\mathbf{K}^{+,a}}(M^a,N^a)\ar[r] & 0} \] where the left hand horizontal maps are given by sending $(s^i)_{i \in \bbZ}$ to $ds + sd$. The first part of the lemma shows that the first two vertical maps are almost isomorphisms. Therefore the third vertical map is an almost isomorphism, as required.
\end{proof}
Note that $\Gamma_{G_n}$ is left exact, as it is a composition of the left exact $\Hom$-functor and the exact localisation functor.
\begin{definition}
Denote by $R\Gamma_{G_n} : \mathbf{D}^+_\text{sm}(G_0,A) \rightarrow \mathbf{D}^+_\text{sm}(G_0/G_n,A)$ the right derived functor of $\Gamma_{G_n}$.
\end{definition}
\begin{lemma}\label{lem_commute_a_and_invts}
There is a natural isomorphism of functors $ R\Gamma_{G_n} (-)^a \cong (R\Gamma_{G_n} (- ))^a : \mathbf{D}^+_\text{sm}(G_0,\bbV) \rightarrow \mathbf{D}^+_\text{sm}(G_0/G_n,A)$.
\end{lemma}
\begin{proof}
Both functors are the right derived functor of $X \mapsto \Gamma_{G_n}(X)^a = \Gamma_{G_n}(X^a)$, since $(-)^a$ is exact and preserves injectives.
\end{proof}

\section{Arithmetic locally symmetric spaces}\label{sec_groups_and_symmetric_spaces}
In this section, we will describe the spaces associated to linear algebraic groups over number fields, and use them to define our derived Hecke algebras.
\subsection{Symmetric spaces}
Let $G$ be a connected linear algebraic group over $\bbQ$, and let $R_d G$ denote the $\bbQ$-split part of the radical of $G$. Following Borel--Serre \cite{Bor73}, we make the following definition.
\begin{definition}
A space of type $S-\bbQ$ for $G$ is a pair consisting of a \(left\) homogeneous space $X_G$ under $G(\bbR)$ and a family $(L_x)_{x \in X}$ of Levi subgroups of $G_\bbR$ satisfying the following two conditions:
\begin{enumerate}
\item The isotropy groups $G_x = \Stab_{G(\bbR)}(x)$ are of the form $G_x = K \cdot S(\bbR)$, where $S \subset R_d G$ is a maximal torus and $K \subset G(\bbR)$ is a maximal compact subgroup normalizing $S$.
\item For each $x \in X$, we have $G_x \subset L_x$ and $L_{g \cdot x} = g L_x g^{-1}$ for all $g \in G(\bbR)$. 
\end{enumerate}
\end{definition}
It follows from \cite[Lemma 2.1]{Bor73} that there is a unique $G(\bbR)$-conjugacy class of such subgroups $G_x = S(\bbR) \cdot K$; the homogeneous space $X_G$ is therefore determined up to isomorphism. It is connected, because a maximal compact subgroup $K$ meets every connected component of $G(\bbR)$. On the other hand, the family of Levi subgroups $(L_x)_{x \in X}$ involves a choice. Henceforth, we write $X_G$ for a fixed choice of space of type $S-\bbQ$. The space $X_G$ is orientable: in fact, it is diffeomorphic to Euclidean space (\cite[Remark 2.4]{Bor73}). 

These spaces are studied in great generality in \cite{Bor73}. For us, examples will arise as follows:
\begin{itemize}
\item If $G$ is reductive, then there is a unique isomorphism class of space of type $S-\bbQ$ for $G$, as is clear from the definition.
\item If $G$ is reductive and $P \subset G$ is a rational parabolic subgroup, then $P(\bbR)$ acts transitively on $X_G$. For any $x \in X_G$, there is a unique Levi subgroup $L_x' \subset P_\bbR$ which is stable under the Cartan involution of $G_\bbR$ associated to $K_x$, the maximal compact subgroup of $G_x$ (hence of $G(\bbR)$); see \cite[(1.9), Corollary]{Bor73}. 

Let $S_P = (R_d P / (R_u P \cdot R_d G))$, a $\bbQ$-split torus, and let $A_P = S_P(\bbR)^\circ$. There is a canonical action of $A_P$ on $X_G$, called the geodesic action, and given by the formula (for $a \in A_P$, $x \in X_G$) $a \bullet x = a_x \cdot x$, where $a_x \in L'_x(\bbR)$ is any lifting of $a \in A_P$, cf. \cite[(3.2)]{Bor73}. This action of $A_P$ commutes with the action of $P(\bbR)$ on $X_G$, which therefore descends to the quotient $X_P = A_P \backslash X_G$. For any $a \in A_P,$ $x \in X_G$, we have $L'_x = L'_{a \bullet x}$, and the quotient $X_P = A_P \backslash X_G$ becomes a space of type $S - \bbQ$ for $P$ when equipped with the family of Levi subgroups $(L'_x)_{x \in X_P}$. 
\end{itemize}
For notational purposes it is convenient to allow groups over arbitrary number fields, so now suppose that $F$ is a number field and that $G$ is a connected linear algebraic group over $F$. We will write $X_G$ for a fixed choice of space of type $S-\bbQ$ for the restriction of scalars $\Res^F_\bbQ G$.

We now consider certain adelic arithmetic quotients of $X_G$. Choose an element $g = (g_v)_v \in G(\bbA_F^\infty)$, and consider for each finite place $v$ the subgroup $\Gamma_v \subset \overline{F}_v^\times$, defined as the torsion subgroup of the subgroup of $\overline{F}_v^\times$ generated by the eigenvalues of $g_v$ in any faithful representation of $G$. The element $g$ is said to be neat if $\cap_v \Gamma_v$ is the trivial group. (This intersection has a sense since finite subgroups of $\overline{\bbQ}^\times$ are invariant under Galois automorphisms.) An open compact subgroup $U \subset G(\bbA_F^\infty)$ is said to be neat if all of its elements are neat. (This definition of neatness is the one used by Pink \cite{Pin90}.)

If $U$ is neat, then for all $g \in G(\bbA_F^\infty)$, the group $\Gamma_{g, U} = G(F) \cap g U g^{-1}$ is neat as an arithmetic subgroup of $G(F \otimes_\bbQ \bbR)$. In particular, it is torsion-free and acts freely and properly discontinuously on $X_G$, preserving orientations. Moreover, if $H \subset G$ is a subgroup then $U \cap H(\bbA_F^\infty)$ is neat, and if $G \to H$ is a surjective homomorphism then the image of $U$ in $H(\bbA_F^\infty)$ is again neat. 

We write $\mathfrak{X}_G = G(F) \backslash \left[ G(\bbA_F^\infty) \times X_G \right]$, where before forming the quotient $G(\bbA_F^\infty)$ is \emph{endowed with the discrete topology}. Then $\mathfrak{X}_G$ is a $G(\bbA_F^\infty)$-space, in the sense of \S \ref{sec_equivariant_sheaves}. It follows that $\mathfrak{X}_G$ is isomorphic to an uncountable disjoint union of connected smooth manifolds, and for any neat open compact subgroup $U \subset G(\bbA_F^\infty)$, $\mathfrak{X}_G$ is a free $U$-space (in the sense of Definition \ref{defn_free_g_space}). If $U$ is such a subgroup, then we write $X_G^U$ for the quotient
\[ X_G^U = \mathfrak{X}_G/U = G(F) \backslash \left[ G(\bbA_F^\infty)/U \times X_G \right]. \]
If $S$ is a finite set of finite places of $F$ then we will write $G^S = G(\bbA_F^{\infty, S})$, where $\bbA_F^{\infty, S}$ is the ring of finite adeles, deprived of its $S$-components, and $G_S = \prod_{v \in S} G(F_v)$. Thus $G(\bbA_F^\infty) = G^S \times G_S$. In a slight abuse of notation, we will also write $G^\infty = G(\bbA_F^\infty)$ and $G_\infty = G(F \otimes_\bbQ \bbR)$. 

In order to describe a reasonable class of level subgroups, we will fix an integral model $\underline{G}$ of $G$, i.e.\ a flat affine group scheme over $\cO_F$ with generic fibre $G$. Such a structure having been fixed, we will write $\cJ_G$ for the set of neat open compact subgroups of $G(\bbA_F^\infty)$ of the form $U = \prod_v U_v$ with $U_v \subset \underline{G}(\cO_{F_v})$ for all $v$. When $G = \GL_{n, F}$ we will always choose the natural integral structure $\underline{G} = \GL_{n, \cO_F}$, in which case $\cJ_G$ is the set of neat open compact subgroups of $\GL_n(\bbA_F^\infty)$ of the form $U = \prod_v U_v$, with $U_v \subset \GL_n(\cO_{F_v})$ for all $v$.
\begin{lemma}\label{lem_smooth_structure_on_X_G} Let $U \in \cJ_G$.
\begin{enumerate}
\item The quotient $G(F) \backslash G^\infty / U$ is finite. Writing $g_1, \dots, g_s \in G^\infty$ for a set of representatives and $\Gamma_{g_i, U} = G(F) \cap g_i U g_i^{-1}$, we have a homeomorphism
\[ X_G^U \cong \coprod_{i=1}^s \Gamma_{g_i, U} \backslash X_G, \]
that we use to endow $X_G^U$ with the structure of orientable smooth manifold.
\item There is an equivalence of categories $\Sh_U(\mathfrak{X}_G) \cong \Sh(X_G^U)$.
\end{enumerate}
\end{lemma}
\begin{proof}
The first part is finiteness of the class number for $G$, which follows from \cite[Theorem 5.1]{Pla94}. The second part follows from Lemma \ref{lem_descent_along_free_actions}, since $U$ acts freely on $\mathfrak{X}_G$.
\end{proof}
We will need to consider some naturally arising families of sheaves on the spaces $X_G^U$. Let $S$ be a set of finite places of $F$ and let $U_S \subset \prod_{v \in S} \underline{G}(\cO_{F_v})$ be an open compact subgroup. We will write $\cJ_{G, U_S} \subset \cJ_G$ for the set of $U \in \cJ_G$ of the form $U = U_S U^S$. If $S$ is finite and $M$ is a $\bbZ[U_S]$-module, viewed as an object of $\Sh_{G^S \times U_S}(\text{pt})$, then we write $\underline{M}_G$ for its pullback to $\Sh_{G^S \times U_S}(\mathfrak{X}_G)$, and $\underline{M}_G^U$ for its image in $\Sh(X_G^U)$. Lemma \ref{lem_smooth_structure_on_X_G} and the diagram (\ref{eqn_hecke_functor_diagram}) then imply:
\begin{corollary}\label{cor_correct_cohomology_groups}
There is a natural isomorphism for any $\bbZ[U_S]$-module $M$:
\[ R \Gamma_U^\sim R \Gamma_{\mathfrak{X}_G} \underline{M}_G \cong R \Gamma_{X_G^U} \underline{M}_G^U. \]
\end{corollary}
This shows that our use of the discrete topology on $G^\infty$ does not cause pathologies.
\subsubsection{Quotient by unipotent radical}

We continue to denote by $G$ a connected linear algebraic group over a number field $F$. As a warm-up for later, we now discuss what happens when when we consider the morphism $G \to H = G / N$, with $N = R_u G$ the unipotent radical of $G$. In this case, the group $N(F \otimes_\bbQ \bbR)$ acts freely on $X_G$ and we can take $X_H = N(F \otimes_\bbQ \bbR) \backslash X_G$ (see \cite[(2.8)]{Bor73}). Let $S$ be a finite set of finite places of $F$, and let $U_S \subset G_S$ be a fixed open compact subgroup. We will freely use the identification $\cH(G^S \times U_S, U) \cong \cH(G^S, U^S)$, and similarly for the groups $H$ and $N$. For any $U \in \cJ_{G, U_S}$, we write $U_H \in \cJ_{H, U_{H, S}}$ for its image in $H^\infty$ and $U_N \in \cJ_{N, U_{N, S}}$ for its intersection with $N^\infty$. There is a natural projection $\pi_{G, H} : \mathfrak{X}_G \to \mathfrak{X}_H$, and for any $U \in \cJ_{G, U_S}$ a quotient projection $\pi_{U, U_H} : X_G^U \to X_H^{U_H}$. The map $\pi_{U, U_H}$ is a submersion with compact nilmanifold fibres.

Now fix a Levi decomposition $G = H \ltimes N$, and fix a subgroup $U \in \cJ_{G, U_S}$ which is decomposed, i.e. such that $U = U_H \ltimes U_N$. In this case we have constructed in Lemma \ref{lem_hecke_integration_along_fibres} a homomorphism $r_H : \cH(G^S, U^S) \to \cH(H^S, U_H^S)$ and a corresponding functor $r_H^\ast : \Mod(\cH(H^S, U_H^S)) \to \Mod(\cH(G^S, U^S))$. In this situation, we want to construct for any $U_{H, S}$-module $A$ a homomorphism in $\mathbf{D}(\cH(G^S, U^S))$:
\[ i : r_H^\ast R \Gamma_{U_H} R \Gamma_{\mathfrak{X}_H} \underline{A}_H \to R \Gamma_U R \Gamma_{\mathfrak{X}_G} \underline{A}_G. \]
To this end, we consider the following diagram of functors:
\[ \xymatrix{ \Sh_{G^S \times U_S}(\mathfrak{X}_{G})  \ar[r]^{\Gamma_{\mathfrak{X}_G}} & \Sh_{G^S \times U_S}(\text{pt}) \\
\Sh_{H^S \times U_{H, S}}(\mathfrak{X}_H) \ar[u]^{\pi^\ast}\ar[r]_{\Gamma_{\mathfrak{X}_H}} & \Sh_{H^S \times U_{H, S}}(\text{pt}) \ar[u]_{\Inf_{H^S \times U_{H, S}}^{G^S \times U_S}}. } \]
Pullback of global sections gives a natural transformation $ \Inf_{H^S \times U_{H, S}}^{G^S \times U_S} \circ  \, \Gamma_{\mathfrak{X}_H} \to \Gamma_{\mathfrak{X}_G} \circ \pi^\ast$; by Lemma \ref{lem_universal_property_of_derived_functor}, we obtain a canonical morphism for any $A \in \Mod(U_{H, S})$: 
\[ \Inf_{H^S \times U_{H, S}}^{G^S \times U_S} R \Gamma_{\mathfrak{X}_H} \underline{A}_H \to R \Gamma_{\mathfrak{X}_G} \pi^\ast \underline{A}_H \cong R \Gamma_{\mathfrak{X}_G} \underline{A}_G. \]
(Note that there is a canonical natural isomorphism $\Inf_{H^S \times U_{H, S}}^{G^S \times U_S} R \Gamma_{\mathfrak{X}_H} \cong R(\Inf_{H^S \times U_{H, S}}^{G^S \times U_S}\Gamma_{\mathfrak{X}_H})$.) Combining this with Corollary \ref{cor_universality_of_derived_functor}, we obtain our desired morphism $i$ as the composite
\[ r_H^\ast R \Gamma_{U_H}  R \Gamma_{\mathfrak{X}_H} \underline{A}_H  \to R \Gamma_U \Inf_{H^S \times U_{H, S}}^{G^S \times U_S} R \Gamma_{\mathfrak{X}_H} \underline{A}_H \to R \Gamma_U R \Gamma_{\mathfrak{X}_G} \underline{A}_G. \]
We now want to construct for each $A \in \Mod(U_{H, S})$ a splitting 
\[ s : (R \Gamma_U R \Gamma_{\mathfrak{X}_G} \pi^\ast \underline{A}_H)^\sim \to (R \Gamma_{U_H} R \Gamma_{\mathfrak{X}_H} \underline{A}_H)^\sim \]
of $i^\sim$. To this end, we introduce a new space $\mathfrak{Y}_H = H(F) \backslash \left[ H(\bbA_F^\infty) \times X_G \right]$, with the action of $H(F)$ on $X_G$ induced by our fixed Levi decomposition $G = H \ltimes N$. There is a natural $H^\infty$-equivariant map $\theta : \mathfrak{Y}_H \to \mathfrak{X}_G$, and the composite $\sigma : \mathfrak{Y}_H \to \mathfrak{X}_G \to \mathfrak{X}_H$ is a fibre bundle with fibre $N(F \otimes_\bbQ \bbR)$. In particular, the endofunctor $R \sigma_\ast \sigma^\ast$ of $\mathbf{D}^+(\Sh_{H^S\times U_{H, S}}\mathfrak{X}_H)$ is naturally isomorphic to the identity functor, by adjunction (cf. \cite[Proposition 2.7.8]{Kas94}).

There is a natural transformation $R \Gamma^\sim_U R \Gamma_{\mathfrak{X}_G} \to R \Gamma_{U_H}^\sim R \Gamma_{\mathfrak{Y}_H} \theta^\ast$. Applying this to a sheaf $\pi^\ast \underline{A}_H$ we obtain our desired morphism $s$ as the composite
\[ R \Gamma^\sim_U R \Gamma_{\mathfrak{X}_G} \pi^\ast \underline{A}_H \to \Gamma_{U_H}^\sim R \Gamma_{\mathfrak{Y}_H} \theta^\ast \pi^\ast \underline{A}_H \cong \Gamma_{U_H}^\sim R \Gamma_{\mathfrak{Y}_H} \sigma^\ast \underline{A}_H \cong \Gamma_{U_H}^\sim R \Gamma_{\mathfrak{X}_H} R\sigma_\ast \sigma^\ast \underline{A}_H \cong \Gamma_{U_H}^\sim R \Gamma_{\mathfrak{X}_H} \underline{A}_H. \]
It is easy to see that $s$ is a splitting of the morphism $i^\sim : R \Gamma^\sim_{U_H} R \Gamma_{\mathfrak{X}_H} \underline{A}_H \to R \Gamma^\sim_U R \Gamma_{\mathfrak{X}_G} \pi^\ast \underline{A}_H$ in $\mathbf{D}(\bbZ)$. Putting all of this together, we have proved the following:
\begin{proposition}\label{prop_pullback_hecke_algebra}
For any $A \in \Mod(U_{H, S})$, there are natural morphisms
\[ i : r_H^\ast R \Gamma_{U_H} R \Gamma_{\mathfrak{X}_H} \underline{A}_H \to R \Gamma_U R \Gamma_{\mathfrak{X}_G} \pi^\ast \underline{A}_H \]
in $\mathbf{D}(\cH(G^S, U^S))$ and
\[ s : R \Gamma^\sim_U R \Gamma_{\mathfrak{X}_G} \pi^\ast \underline{A}_H \to R \Gamma^\sim_{U_H} R \Gamma_{\mathfrak{X}_H} \underline{A}_H \]
in $\mathbf{D}(\bbZ)$, satisfying $s i^\sim = 1$. In particular, there is a commutative diagram of $\bbZ$-algebras:
\[ \xymatrix{ \cH(G^S, U^S)\ar[d]_{r_H} \ar[r]^-{T_G} & \End_{\mathbf{D}(\bbZ)}(R \Gamma^\sim_U R \Gamma_{\mathfrak{X}_G} \pi^\ast \underline{A}_H) \ar[d]^-{t \mapsto s t i^\sim}
\\ \cH(H^S, U_H^S) \ar[r]^-{T_H} & \End_{\mathbf{D}(\bbZ)}(R \Gamma^\sim_{U_H} R \Gamma_{\mathfrak{X}_H} \underline{A}_H). } \]
\end{proposition}
\begin{proof}
It remains to check that the diagram is commutative. For this, it is enough to show that for any $t \in \cH(H^S, U_H^S)$, the equality $T_G(t) i^\sim = i^\sim T_H(r_H(t))$ holds inside 
\[ \Hom_{\mathbf{D}(\bbZ)}(R \Gamma^\sim_{U_H} R \Gamma_{\mathfrak{X}_H}, R \Gamma^\sim_U R \Gamma_{\mathfrak{X}_G} \pi^\ast \underline{A}_H). \]
This follows immediately from the fact that $i^\sim$ arises from a map in $\mathbf{D}(\cH(G^S, U^S))$.
\end{proof}
We now generalize this slightly.
\begin{proposition}\label{prop_push_pull_hecke_algbera}
Let $B$ be a $\bbZ[U_S]$-module, and let $C = \Res^{U_S}_{U_{H, S}} B$. Suppose that $C$ admits a decomposition $C = A \oplus K$ of $\bbZ[U_{H, S}]$-modules, where $A \subset B^{U_{N, S}}$. Then:
\begin{enumerate}
\item We have $\pi^\ast \underline{A}_H = \underline{A}_G$.
\item There are natural morphisms
\[ i : r_H^\ast R \Gamma_{U_H} R \Gamma_{\mathfrak{X}_H} \underline{A}_H \to R \Gamma_U R \Gamma_{\mathfrak{X}_G} \underline{B}_G \]
in $\mathbf{D}(\cH(G^S, U^S))$ and
\[ s : R \Gamma^\sim_U R \Gamma_{\mathfrak{X}_G} \underline{B}_G \to R \Gamma^\sim_{U_H} R \Gamma_{\mathfrak{X}_H} \underline{A}_H \]
in $\mathbf{D}(\bbZ)$, satisfying $s i^\sim = 1$. In particular, there is a commutative diagram of $\bbZ$-algebras:
\[ \xymatrix{ \cH(G^S, U^S)\ar[r]^-{T_G} \ar[d]_{r_H} & \End_{\mathbf{D}(\bbZ)}(R \Gamma^\sim_U R \Gamma_{\mathfrak{X}_G}  \underline{B}_G) \ar[d]^{t \mapsto s t i^\sim} 
\\\cH(H^S, U_H^S) \ar[r]_-{T_H} & \End_{\mathbf{D}(\bbZ)}(R \Gamma^\sim_{U_H} R \Gamma_{\mathfrak{X}_H} \underline{A}_H). } \]
\end{enumerate}
\end{proposition}
\begin{proof}
The isomorphism $\pi^\ast \underline{A}_H \cong \underline{A}_G$ is clear from the definitions. The inclusion $A \subset B^{U_{N_S}}$ implies the existence of a $U_S$-equivariant map $A \to B$, hence $\pi^\ast \underline{A}_H \to \underline{B}_G$. We define $i$ as the composite
\[ r_H^\ast R\Gamma_{U_H} R \Gamma_{\mathfrak{X}_H} \underline{A}_H \to R \Gamma_U R \Gamma_{\mathfrak{X}_G} \pi^\ast \underline{A}_H \to R \Gamma_U R \Gamma_{\mathfrak{X}_G} \underline{B}_G, \]
where the first arrow is the one constructed in Proposition \ref{prop_pullback_hecke_algebra} and the second arrow arises from the map $\pi^\ast \underline{A}_H \to \underline{B}_G$. We define $s$ as the composite
\[ \begin{split} R \Gamma_U^\sim R \Gamma_{\mathfrak{X}_G} \underline{B}_G \to R \Gamma_{U_H}^\sim R \Gamma_{\mathfrak{Y}_H} \theta^\ast \underline{B}_G \cong R \Gamma_{U_H}^\sim R \Gamma_{\mathfrak{Y}_H} \sigma^\ast \underline{C}_H \cong R \Gamma_{U_H}^\sim R \Gamma_{\mathfrak{X}_H} \underline{C}_H\\  \cong R \Gamma_{U_H}^\sim R \Gamma_{\mathfrak{X}_H} (\underline{A}_H \oplus \underline{K}_H) \to R \Gamma_{U_H}^\sim R \Gamma_{\mathfrak{X}_H} \underline{A}_H. \end{split}\]
The first arrow is constructed as in the proof of Proposition \ref{prop_pullback_hecke_algebra}, the isomorphism $\theta^\ast \underline{B}_G \cong \sigma^\ast \underline{C}_H$ follows from the definitions, the second isomorphism follows as in the proof of Proposition \ref{prop_pullback_hecke_algebra}, the third isomorphism follows from the isomorphism $C \cong A \oplus K$, and the final arrow is projection onto the first factor of $C = A \oplus K$. 

To complete the proof of the proposition, it remains to show that the equality $s i^\sim = 1$ holds inside $\End_{\mathbf{D}(\bbZ)}(R \Gamma_{U_H}^\sim R \Gamma_{\mathfrak{X}_H} \underline{A}_H)$. The composite $s i^\sim$ is equal to the composite
\[ R \Gamma^\sim_{U_H} R \Gamma_{\mathfrak{X}_H} \underline{A}_H \to R \Gamma^\sim_{U_H} R \Gamma_{\mathfrak{X}_H} \underline{C}_H \to R \Gamma^\sim_{U_H} R \Gamma_{\mathfrak{X}_H} \underline{A}_H, \]
where the first arrow is induced by the inclusion $A \subset C$ of $\bbZ[U_{H, S}]$-modules and the second by projection along the direct sum decomposition $C = A \oplus K$. It follows that $s i^\sim$ is induced by the identity morphism of $A$, hence is equal to the identity. This completes the proof.
\end{proof}
In the applications, we will need this result in a slightly different form: 
\begin{corollary}\label{cor_parabolic_hecke_action_on_compactly_supported_cohomology}
With assumptions as in Proposition \ref{prop_push_pull_hecke_algbera}, there exists a commutative diagram
\[ \xymatrix{ \cH(G^S, U^S) \ar[r]^-{T_G} \ar[d]^{r_H}  & \End_{\mathbf{D}(\bbZ)}(R \Gamma_{X_G^U} \underline{B}^U_G)\ar[d]^{t \mapsto s t i^\sim} \ar[r]^-{t \mapsto t F} & \Hom_{\mathbf{D}(\bbZ)}(R \Gamma_{X_G^U, c} \underline{B}_G^U, R \Gamma_{X_G^U} \underline{B}_G^U) \ar[d]^{t \mapsto s t i_c} \\
\cH(H^S, U_H^S)\ar[r]^-{T_H} & \End_{\mathbf{D}(\bbZ)}(R \Gamma_{X_H^{U_H}} \underline{A}^{U_H}_H) \ar[r]^-{t \mapsto t F} & \Hom_{\mathbf{D}(\bbZ)} (R \Gamma_{X_H^{U_H}, c} \underline{A}^{U_H}_H, R \Gamma_{X_H^{U_H}} \underline{A}_H^{U_H}). } \]
\end{corollary}
\begin{proof}
We first define the relevant objects and morphisms. For a space $X$, $R \Gamma_{X, c}$ denotes cohomology with compact support. The maps $s$ and $i$ are as in the statement of the proposition, we write $F$ for the natural `forget supports' maps $R \Gamma_c \to R \Gamma$, and $i_c$ is the natural pullback
\[ R \Gamma_{X_H^{U_H}, c} \underline{A}_H^{U_H} \to R \Gamma_{X_G^{U}, c} \underline{A}_G^{U}  \to R \Gamma_{X_G^U, c} \underline{B}_G^U, \]
which exists because $\pi_{U, U_H}$ is proper. The corollary now follows from the proposition and the commutativity of the following diagram for any $t \in \End_{\mathbf{D}(\bbZ)}(R \Gamma_{X_G^{U}} \underline{B}_G^{U})$:
\[ \xymatrix{ R \Gamma_{X_G^U, c} \underline{B}_G^U \ar[r]^F  & R \Gamma_{X_G^U} \underline{B}_G^U \ar[r]^{t} & R \Gamma_{X_G^U} \underline{B}_G \ar[d]^s \\
R \Gamma_{X_H^{U_H}, c} \underline{A}_H^{U_H} \ar[r]^F \ar[u]^{i_c} & R \Gamma_{X_H^{U_H}} \underline{A}_H^{U_H} \ar[r]^{s t i^\sim} \ar[u]^{i^\sim} & R \Gamma_{X_H^{U_H}} \underline{A}_H^{U_H}. } \]

\end{proof}
\subsubsection{Borel--Serre compactifications and restriction to parabolic subgroups}

We continue to suppose that $G$ is a connected linear algebraic group over the number field $F$. Let $U \in \cJ_G$. According to Borel--Serre \cite{Bor73}, we can add boundary strata to $X_G^U$ in order to obtain compact manifolds with corners. We now discuss some elements of this theory.

Let $\mathfrak{P}$ denote the set of $F$-rational parabolic subgroups of $G$ (which includes $G$ itself). For each $P \in \mathfrak{P}$, we have defined the group $A_P$ and observed that the quotient $e(P) = A_P \backslash X_G$ admits a canonical structure of space of type $S-\bbQ$ for $\Res^F_\bbQ P$, with respect to which the map $X_G \to e(P)$ is $P(F \otimes_\bbQ \bbR)$-equivariant. Accordingly, we define
\[ \overline{X}_G = \coprod_{P \in \mathfrak{P}} e(P), \]
endowed with the structure of smooth manifold with corners described in \cite[\S 7.1]{Bor73}. For each $P \in \mathfrak{P}$, the subset $X(P) = \coprod_{Q \supset P} e(Q)$ is an open subset of $\overline{X}_G$, the structure of which can be described explicitly, see \cite[\S 5]{Bor73}. In particular, $e(G) = X_G \subset \overline{X}_G$ is an open submanifold. If $g \in G(F)$, then there is a natural isomorphism $X(P) \to X(P^g)$; the action of $G(F)$ on $X_G$ extends naturally to $\overline{X}_G$ in a way compatible with these isomorphisms, see \cite[Proposition 7.6]{Bor73}.

We define $\overline{\mathfrak{X}}_G = G(F) \backslash \left[ G(\bbA_F^\infty) \times \overline{X}_G \right]$, where as in the previous section, $G(\bbA_F^\infty)$ gets the discrete topology in the formation of the quotient. For any $U \in \cJ_G$, we define 
\[ \overline{X}_G^U = \overline{\mathfrak{X}}_G/U = G(F) \backslash \left[ G(\bbA_F^\infty) / U \times \overline{X}_G \right]. \]
As in the previous section, we can choose representatives $g_1, \dots, g_s \in G(\bbA_F^\infty)$ for the finite double quotient $G(F) \backslash G(\bbA_F^\infty) / U$ and calculate
\[ \overline{X}_G^U = \coprod_{i=1}^s \Gamma_{g_i, U} \backslash \overline{X}_G. \]
For each $g \in G(\bbA_F^\infty)$, the neat arithmetic subgroup $\Gamma_{g_i, U} \subset G(\bbR)$ acts freely on $\overline{X}_G^U$, and the quotient $\Gamma_{g_i, U} \backslash \overline{X}_G$ is compact (\cite[Theorem 9.3]{Bor73}).

We define $\partial \overline{\mathfrak{X}}_G = \overline{\mathfrak{X}}_G - \mathfrak{X}_G$ and $\partial \overline{X}_G^U = \overline{X}_G^U - X_G^U$. Then we have similarly
\[ \partial \overline{X}_G^U = \partial \overline{\mathfrak{X}}_G / U = \coprod_{i=1}^s \Gamma_{g_i, U} \backslash \partial\overline{X}_G. \]
Suppose given a finite set $S$ of finite places of $F$ and a fixed open compact subgroup $U_S \subset G_S$. For any $\bbZ[U_S]$-module $A$, we will write $\underline{A}_G \in \Sh_{G^S \times U_S}(\overline{\mathfrak{X}}_G)$ for its pullback from $\Sh_{G^S \times U_S}(\text{pt})$. Since the pullback of $\underline{A}_G$ to the $G^\infty$-invariant open submanifold $\mathfrak{X}_G \subset \overline{\mathfrak{X}}_G$ agrees with the equivariant sheaf previously denoted as $\underline{A}_G$, we hope that this will not cause confusion.

We can use the Borel--Serre compactification to define a Hecke action on the compactly supported cohomology of the spaces $X_G^U$. More precisely, we can define for any $\bbZ[U_S]$-module $A$ a homomorphism 
\[ \cH(G^S, U^S) \to \End_{\mathbf{D}(\bbZ)}(R \Gamma_{X_G^U, c} \underline{A}_G^U) \]
which is compatible with the natural morphism
\[ R \Gamma_{X_G^U, c} \underline{A}_G^U \to R \Gamma_{X_G^U} \underline{A}_G^U. \]
To do this, let us write $j_G : \mathfrak{X}_G \to \overline{\mathfrak{X}}_G$ for the natural open immersion and $j_G^U : X_G^U \to \overline{X}_G^U$ for the corresponding open immersion at finite level. We can take $R \Gamma_{X_G^U, c} \underline{A}_G^U = R \Gamma_{\overline{X}_G^U} j^U_{G, !} \underline{A}_G^U$. It now suffices to observe that $j_{G}$ induces a functor $j_{G, !} : \Sh_{G^\infty}(\mathfrak{X}_G) \to \Sh_{G^\infty}(\overline{\mathfrak{X}}_G)$ and that there is a canonical isomorphism
\[ R\Gamma_U^\sim R \Gamma_{\overline{\mathfrak{X}}_G} j_{G, !} \underline{A}_G \cong R \Gamma_{\overline{X}_G^U} j_{G, !}^U \underline{A}_G^U; \]
this follows easily from the observation that $U$ acts freely on $\overline{\mathfrak{X}}_G$, as in Corollary \ref{cor_correct_cohomology_groups}. The Hecke actions we have defined are related by Verdier duality as follows.
\begin{proposition}\label{prop_hecke_action_and_verdier_duality}
Let $R$ be a Noetherian ring and let $A \in \Mod(U_S, R)$ be finite free as $R$-module; let $B = \Hom_R(A, R)$.
\begin{enumerate}
\item There is a natural Verdier duality isomorphism in $\mathbf{D}(R)$:
\begin{equation}\label{eqn_verdier_duality_isomorphism} R Hom_R ( R \Gamma_{X_G^U, c} \underline{A}_G^U, R ) \cong R \Gamma_{X_G^U} \underline{B}_G^U. 
\end{equation}
\item Let $S$ be a Noetherian $R$-algebra and let $A_S = A \otimes_R S$, $B_S = B \otimes_R S$. Then there are natural isomorphisms
\[ R \Gamma_{X_G^U, c} \underline{A_S}_G^U \cong (R \Gamma_{X_G^U, c} \underline{A}_G^U ) \otimes^\bbL_R S\]
and
\[ R \Gamma_{X_G^U} \underline{B_S}_G^U \cong (R \Gamma_{X_G^U} \underline{B}_G^U ) \otimes^\bbL_R S. \]
\item For $g \in G^S$, the Verdier duality isomorphism (\ref{eqn_verdier_duality_isomorphism}) identifies the transpose of the operator $[U g^{-1} U]$ on the left hand side with the operator $[U g U]$ on the right hand side.
\end{enumerate}
\end{proposition}
\begin{proof}
The first part follows from the usual Verdier duality isomorphism \cite{Ver95}. We have used the fact that the derived sheaf Hom $RHom_R$ in $\Sh(X_G^U, R)$ satisfies $RHom_R(\underline{A}_G^U, R) \cong \underline{B}_G^U$. The second part follows as in \cite[Chapitre 2, 4.12]{SGA4demi}; note that the functors $\Gamma_{X_G^U}$ and $\Gamma_{X_G^U, c}$ have finite cohomological dimension (\cite[Proposition 3.2.3]{Kas94}).
The third part follows from the explicit formula of Lemma \ref{lem_hecke_action_as_expected} and its analogue for cohomology with compact support and the functoriality of Verdier duality.
\end{proof}
Now suppose that $G$ is reductive and that $P$ is a maximal proper parabolic subgroup of $G$. Then $e(P) \subset \partial \overline{X}_G$ is an open submanifold, and we write $j_P : \mathfrak{X}_P \to \partial \overline{\mathfrak{X}}_G$ for the induced $P^\infty$-equivariant open immersion. This leads to an exact functor $j_{P, !} : \Sh_{P^S \times U_{P, S}}(\mathfrak{X}_P) \to \Sh_{P^S \times U_{P, S}}(\partial \overline{\mathfrak{X}}_G)$. By passage to quotient, we obtain an open immersion $j_P^U : X_P^{U_P} \to \partial \overline{X}_G^U$.
\begin{proposition}\label{prop_restriction_of_cohomology_to_boundary} Let $G$ be a reductive group and let $P \subset G$ be a maximal proper parabolic subgroup. Let $U \in \cJ_{G, U_S}$ and let $B$ be a $\bbZ[U_S]$-module. We consider $B$ also as a $U_{P, S}$-module by restriction.
\begin{enumerate}
\item We have a canonical isomorphism 
\[ R \Gamma_{U_P}^\sim R \Gamma_{\partial \overline{\mathfrak{X}}_G} j_{P, !} \underline{B}_P \cong R \Gamma_{X_P^{U_P}, c} \underline{B}_P. \]
\item There are natural morphisms
\[ p :  R \Gamma_{U_P}^\sim R \Gamma_{\partial \overline{\mathfrak{X}}_G} j_{P, !} \underline{B}_P \to R \Gamma_U^\sim R \Gamma_{\partial \overline{\mathfrak{X}}_G} \underline{B}_G \]
in $\mathbf{D}(\bbZ)$ and
\[ q : R \Gamma_U R \Gamma_{\partial \overline{\mathfrak{X}}_G} \underline{B}_G \to  r_P^\ast R \Gamma_{U_P} R \Gamma_{\mathfrak{X}_P}  \underline{B}_P \]
in $\mathbf{D}(\cH(G^S, U^S))$. The morphism $q^\sim p$ in $\Hom_{\mathbf{D}(\bbZ)}(R \Gamma_{{X}_P^{U_P}, c} \underline{B}_P^{U_P}, R \Gamma_{{X}_P^{U_P}} \underline{B}_P^{U_P})$ is the canonical one \(arising from the `forget supports' map $R \Gamma_c \to R \Gamma$\), and we obtain a commutative diagram of $\bbZ$-modules:
\[ \xymatrix{ \cH(G^S, U^S) \ar[d]_{r_P} \ar[r]^-{T_G} & \End_{\mathbf{D}(\bbZ)}(R \Gamma_{\partial \overline{X}_G^U} \underline{B}_G^U ) \ar[d]^{t \mapsto q^\sim t p} \\
\cH(P^S, U_P^S) \ar[r]_-{T_P} & \Hom_{\mathbf{D}(\bbZ)}(R \Gamma_{{X}_P^{U_P}, c} \underline{B}_P^{U_P}, R \Gamma_{{X}_P^{U_P}} \underline{B}_P^{U_P}). } \]
\end{enumerate} 
\end{proposition}
\begin{proof}
For the first part, it is enough to note that there is a canonical isomorphism 
\[ R \Gamma_{U_P}^\sim R \Gamma_{\partial \overline{\mathfrak{X}}_G} j_{P, !} \underline{B}_P \cong R \Gamma_{\partial \overline{X}_G^U} j^U_{P, !} \underline{B}_P, \]
because $U$ acts freely on $\partial \overline{\mathfrak{X}}_G$. The isomorphism with $R \Gamma_{X_P^{U_P}, c} \underline{B}_P$ then follows from the fact that $j^U_{P, !}$ takes injectives to soft sheaves, which are $\Gamma_c$-acyclic, see \cite[Proposition III.7.2]{Ive86}.

We now construct the morphisms $p$ and $q$. First, $p$ is the morphism
\[ R \Gamma_{\partial \overline{X}_G^U} j_{P, !}^U j_{P}^{U, \ast} \underline{B}_G \to R \Gamma_{\partial \overline{X}_G^U} \underline{B}_G \]
which arises from the natural map $j_{P, !}^U j_{P}^{U, \ast} \underline{B}_G^U \to \underline{B}_G^U$ (note that $j_P^\ast \underline{B}_G = \underline{B}_P$). Next, $q$ is obtained by applying Corollary \ref{cor_derived_hecke_morphism_for_parabolic_restriction} to the morphism $\Res^{G^S \times U_S}_{P^S \times U_{P, S}} R \Gamma_{\partial \overline{\mathfrak{X}}_G} \underline{B}_G \to R \Gamma_{\mathfrak{X}_P} \underline{B}_P$ induced by pullback. To complete the proof of the proposition, it remains to check that $q^\sim p$ is the morphism induced by the natural `forget supports' transformation $R \Gamma_c \to R \Gamma$. However, it follows from the definitions that $q^\sim p$ is equal to the composite
\[ R \Gamma_{{X}_P^{U_P}, c} \underline{B}_P \cong R \Gamma_{\partial \overline{X}_G^U} j_{P, !}^U \underline{B}_P \to R \Gamma_{\partial \overline{X}_G^U} \underline{B}_G \to R \Gamma_{X_P^{U_P}} \underline{B}_P, \]
where the map in the middle is induced by $j^U_{P, !} j_{P}^{U, \ast} \underline{B}_G \to \underline{B}_G$ and the last by pullback. This is the correct map.
\end{proof}

\begin{corollary}\label{cor_restriction_of_hecke_to_boundary}
Let notation and assumptions be as in Proposition \ref{prop_restriction_of_cohomology_to_boundary}. Fix a Levi decomposition $P = M \ltimes N$ and suppose that $U_P = U_M \cdot U_N$ is decomposed. Let $B$ be a $\bbZ[U_S]$-module equipped with a decomposition $\Res^{U_S}_{U_{M, S}} B = A \oplus K$, where $A \subset B^{U_{N, S}}$. Then there exists a commutative diagram:
\[ \xymatrix{ \cH(G^S, U^S) \ar[r]^-{T_G} \ar[d]_{r_P} & \End_{\mathbf{D}(\bbZ)}(R \Gamma_{\partial \overline{X}_G^U} \underline{B}_G^U) \ar[d] \\
\cH(P^S, U_P^S) \ar[r]^-{T_P} \ar[d]_{r_M} & \Hom_{\mathbf{D}(\bbZ)}(R \Gamma_{X_P^{U_P}, c} \underline{B}_P^{U_P}, R \Gamma_{X_P^{U_P}} \underline{B}_P^{U_P}) \ar[d] \\
\cH(M^S, U_M^S) \ar[r]^-{T_M} &  \Hom_{\mathbf{D}(\bbZ)}(R \Gamma_{X_M^{U_M}, c} \underline{A}_M^{U_M}, R \Gamma_{X_M^{U_M}} \underline{A}_M^{U_M}) } \]
\end{corollary}
\begin{proof}
This follows immediately by combining Corollary \ref{cor_parabolic_hecke_action_on_compactly_supported_cohomology} and Proposition \ref{prop_restriction_of_cohomology_to_boundary}.
\end{proof}
An important fact is that the boundary $\partial \overline{\mathfrak{X}}_G$ admits a $G^\infty$-invariant stratification, with strata indexed by conjugacy classes of rational parabolic subgroups of $G$. More precisely, let $P$ be a rational parabolic subgroup of $G$. Then there is a $G^\infty$-equivariant isomorphism
\[ \Ind_{P^\infty}^{G^\infty} \mathfrak{X}_P \cong P(F) \backslash \left[ G(\bbA_F^\infty) \times e(P) \right], \]
and the induced map $\Ind_{P^\infty}^{G^\infty} \mathfrak{X}_P \to \partial \overline{\mathfrak{X}}_G$ is a $G^\infty$-equivariant locally closed immersion. (We define $\Ind_{P^\infty}^{G^\infty} \mathfrak{X}_P = G^\infty \times_{P^\infty} \mathfrak{X}_P$, as in Proposition \ref{prop_shapiros_lemma_for_spaces}.) We then have the following lemma.
\begin{lemma}\label{lem_boundary_stratification_by_induced_strata}
Let $P_1, \dots P_s$ be representatives of the distinct $G(F)$-conjugacy classes of proper rational parabolic subgroups of $G$. Then:
\begin{enumerate}
\item The natural maps $j_{P_i} : \Ind_{P_i^\infty}^{G^\infty} \mathfrak{X}_{P_i}  \to \partial \overline{\mathfrak{X}}_G$ are locally closed immersions, and the induced map $\coprod_i \Ind_{P_i^\infty}^{G^\infty} \mathfrak{X}_{P_i} \to \partial \overline{\mathfrak{X}}_G$ is a continuous bijection.
\item For each $U \in \cJ_G$, the quotient maps $j_{P_i}^U : (\Ind_{P_i^\infty}^{G^\infty} \mathfrak{X}_{P_i})/U \to \partial \overline{X}^U_G$ are locally closed immersions, and the induced map $\coprod_i [\Ind_{P_i^\infty}^{G^\infty} \mathfrak{X}_{P_i}]/U \to \partial \overline{X}^U_G$ is a continuous bijection.
\end{enumerate}
\end{lemma}
\begin{proof}
The second part follows from the first. For the first, we need to show this map is bijective. We simply calculate
\[ \partial\overline{\mathfrak{X}}_G  = G(F) \backslash  [G^\infty \times \partial\overline{X}_G ]= \coprod_i \coprod_{P' \sim P_i} G(F) \backslash [G^\infty \times e(P')], \]
where the second disjoint union is over rational parabolics $P'$ which are $G(F)$-conjugate to $P_i$. Since a parabolic subgroup is its own normalizer, this becomes
\[ \coprod_i P_i(F) \backslash [G^\infty \times e(P)] = \coprod_i \Ind_{P_i^\infty}^{G^\infty} \mathfrak{X}_{P_i}, \]
as desired.
\end{proof}

\subsection{Derived Hecke algebras and the idempotents associated to maximal ideals}\label{sec_idempotents_in_derived_category}

We now introduce some more `automorphic' notation. Let $F$ be a number field, and $G$ a connected reductive group over $F$. Fix a prime $p$ and a choice of finite extension $E/\bbQ_p$ with ring of integers $\cO$, uniformizer $\pi$, and residue field $k = \cO/(\pi)$. Let $S$ be a finite set of finite places of $F$, containing the $p$-adic places. Let $\underline{G}$ be an integral model of the group $G$ such that $\underline{G}_{\cO_{F, S}}$ is reductive, and let $U^S = \prod_{v \not\in S} \underline{G}(\cO_{F})$. 

We write $\bbT^S = \cH(G^S, U^S) \otimes_\bbZ \cO$; then $\bbT^S$ is a commutative $\cO$-algebra, because $U^S$ is a product of hyperspecial maximal compact subgroups. When we wish to emphasize the ambient group $G$, we will write $\bbT^S = \bbT^S_G$. If $C^\bullet$ is a perfect complex of $\cO$-modules (which in this context just means that $H^\ast(C^\bullet)$ is a finite $\cO$-module) equipped with a homomorphism $\bbT^S \to \End_{\mathbf{D}(\cO)}(C^\bullet)$ of $\cO$-algebras, then we will write $\bbT^S(C^\bullet)$ for the quotient of $\bbT^S$ which injects into $\End_{\mathbf{D}(\cO)}(C^\bullet)$. Thus $\bbT^S(C^\bullet)$ is a commutative finite $\cO$-algebra, equipped with a surjective map
\[ \bbT^S(C^\bullet) \to \bbT^S(H^\ast(C^\bullet)) \]
which has nilpotent kernel (because $C^\bullet$ is perfect; see \cite[Lemma 2.5]{Tho14}).

Being a finite $\cO$-algebra, we can decompose $\bbT^S(C^\bullet)$ as a product $\bbT^S(C^\bullet) = \prod_\ffrm \bbT^S(C^\bullet)_\ffrm$ over the finitely many maximal ideals $\ffrm \subset \bbT^S(C^\bullet)$. For each such maximal ideal there is a corresponding idempotent $e_\ffrm \in \bbT^S(C^\bullet) \subset \End_{\mathbf{D}(\cO)}(C^\bullet)$, which is the projector onto the factor $\bbT^S(C^\bullet)_\ffrm$. The derived category $\mathbf{D}(\cO)$ is idempotent complete, so we deduce the existence of a direct sum decomposition $C^\bullet = C^\bullet_\ffrm \oplus D^\bullet$ in $\mathbf{D}(\cO)$. The summand $C^\bullet_\ffrm$ is defined uniquely up to unique isomorphism in $\mathbf{D}(\cO)$, and the composite map 
\[ \bbT^S(C^\bullet)_\ffrm \to \bbT^S(C^\bullet_\ffrm) \]
is an isomorphism. Similarly, there is a canonical identification $H^\ast(C^\bullet)_\ffrm \cong H^\ast(C^\bullet_\ffrm)$. (For a similar but more detailed discussion, see \cite[\S 2.4]{Tho14}.) 

Fix $U \in \cJ_{G, U^S}$, and let $A$ be an $\cO[U_S]$-module, finite over $\cO$. Then $R \Gamma_{X_G^U} \underline{A}_G^U$ is a perfect complex of $\cO$-modules, equipped with a canonical homomorphism 
\[ \bbT^S \to  \End_{\mathbf{D}(\cO)}(R \Gamma_{X_G^U} \underline{A}_G^U ). \]
The algebras $\bbT^S(R \Gamma_{X_G^U} \underline{A}_G^U)$ are the derived Hecke algebras referred to in the introduction of this paper. In the forthcoming sections, we will use the decomposition of the complex $R \Gamma_{X_G^U} \underline{A}_G^U $ according to maximal ideals of this algebra in order to study their associated Galois representations.

The following results will be useful later.
\begin{lemma}\label{Mittag-Leffler}
Let $M$ and $N$ be perfect complexes of $\cO$-modules. Then the natural map
\[\Hom_{\mathbf{D}(\cO)}(M,N) \rightarrow \varprojlim_r \Hom_{\mathbf{D}(\cO/\pi^r\cO)}(M\otimes_{\cO}^\mathbb{L}\cO/\pi^r\cO,N\otimes_{\cO}^\mathbb{L}\cO/\pi^r\cO)\] is an isomorphism.
\end{lemma}
\begin{proof}
For $r = \infty$ (resp. $r \in \mathbb{N}$) let $\mathbf{K}_r$ denote the category with objects bounded complexes of finite free $\cO$-modules (resp. $\cO/\pi^r\cO$-modules) and morphisms given by morphisms of complexes modulo homotopy. The obvious functors $\mathbf{K}_\infty \rightarrow \mathbf{D}(\cO)$, $\mathbf{K}_r \rightarrow \mathbf{D}(\cO/\pi^r\cO)$ are fully faithful, so it suffices to prove that the natural map 
\[\Hom_{\mathbf{K}_\infty}(M,N) \rightarrow \varprojlim_r \Hom_{\mathbf{K}_r}(M\otimes^\mathbb{L}_{\cO}\cO/\pi^r\cO,N\otimes^\mathbb{L}_{\cO}\cO/\pi^r\cO)\] is an isomorphism for all $M, N \in \mathbf{K}_\infty$. This is the content of \cite[Lemma 2.13, (iii)]{Tho14}.
\end{proof}
\begin{lemma}\label{lem_inverse_limit_of_mod_p_hecke_is_hecke}
Let $M$ be a perfect complex of $\cO$-modules endowed with a homomorphism $\bbT_G^S \to \End_{\mathbf{D}(\cO)}(M)$. Let $\bbT_\infty = \bbT_G^S(M)$, and for each $N \geq 1$, let $\bbT_N = \bbT_G^S(M\otimes^\mathbb{L}_\cO \cO/(\pi^N))$. Then the natural map $\bbT_\infty \to \plim_N \bbT_N$ is an isomorphism.
\end{lemma}
\begin{proof}
The map is injective, by Lemma \ref{Mittag-Leffler}. It is surjective because each map $\bbT_\infty \to \bbT_N$ and $\bbT_{N+1} \to \bbT_N$ is surjective, and $\bbT_\infty$ is compact. It is therefore an isomorphism.
\end{proof}
\begin{lemma}\label{lem_square_zero_endomorphisms}
If $\mathbf{C}$ is a triangulated category, $A \to B \to C \to A[1]$ is an exact triangle in $\mathbf{C}$, and $s, t : B \to B$ are morphisms making the diagram
\[ \xymatrix{ A \ar[r] \ar[d]^0 & B\ar[d]^s \ar[r] & C\ar[d]^0  \ar[r] & A[1] \ar[d]^0 \\
A \ar[r] & B \ar[r] & C \ar[r] & A[1] } \]
(and its analogue with $s$ replaced by $t$) commute, then $st = ts = 0 $ in $\End_{\mathbf{C}}(B)$.
\end{lemma}
\begin{proof}
The proof is an easy diagram chase (apply the functor $\Hom(B, -)$).
\end{proof}

\section{The boundary cohomology of the $\GL_n$ locally symmetric space}\label{sec_i_dont_see_any_boundary}

We fix throughout this section a base number field $F$, a prime $p$, and an integer $n \geq 1$. Let $G = \GL_{n, F}$. We fix as well a finite set $S$ of finite places of $F$, containing the $p$-adic places, and set $U^S = \prod_{v \notin S}\GL_n(\cO_{F_v}) \subset G(\bbA_F^{\infty, S})$. Finally, we fix a finite extension $E/\bbQ_p$, and let $\cO$ denote the ring of integers of $E$, $\pi \in E$ a choice of uniformizer, and $k = \cO/(\pi)$ the residue field.

If $m \geq 1$ is an integer, then the Hecke algebra $\bbT^S_{\GL_m} = \cH(\GL_m(\bbA_F^S), \prod_{v \not\in S} \GL_m(\cO_{F_v})) \otimes_\bbZ \cO$ is a commutative $\cO$-algebra, generated by the elements $T_v^i$, $i = 1, \dots, m$, $v \not\in S$, where
\[ T_v^i = \left[ \GL_m(\cO_{F_v}) \diag(\underbrace{\varpi_v, \dots, \varpi_v}_i, \underbrace{1, \dots, 1}_{m-i}) \GL_m(\cO_{F_v}) \right], \]
together with $(T_v^m)^{-1}$. We will say that a perfect complex $C^\bullet$ of $\cO$-modules equipped with a map $\bbT^S_{\GL_m} \to \End_{\mathbf{D}(\cO)}(C^\bullet)$ is of $S$-Galois type if for each maximal ideal $\ffrm \subset \bbT^S_{\GL_m}(C^\bullet)$, there exists a continuous semi-simple representation $\overline{\rho}_\ffrm : G_{F, S} \to \GL_m(\bbT^S(C^\bullet)/\ffrm)$ such that for every finite place $v\not\in S$ of $F$, $\overline{\rho}_\ffrm(\Frob_v)$ has characteristic polynomial
\begin{equation}\label{eqn_characteristic_polynomial_of_Frobenius}X^n - T_v^1 X^{n-1} + \dots + (-1)^j q_v^{j(j-1)/2} T_v^j X^{n-j} + \dots + (-1)^n q_v^{n(n-1)/2} T_v^n \in (\bbT^S(C^\bullet)/\ffrm)[X].
\end{equation}
If $C^\bullet$ is of $S$-Galois type and a representation $\overline{\rho}_\ffrm$ is absolutely reducible, then we say that the maximal ideal $\ffrm$ is Eisenstein. If $C^\bullet$ is of $S$-Galois type and every maximal ideal $\ffrm \subset \bbT^S_{\GL_m}(C^\bullet)$ is Eisenstein, we say that the complex $C^\bullet$ itself is Eisenstein. Note that these conditions hold for a given complex $C^\bullet$ if and only if they hold for the cohomology $H^\ast(C^\bullet)$.

We now suppose for the rest of \S \ref{sec_i_dont_see_any_boundary} that the following hypothesis holds:
\begin{itemize}
\item[$\spadesuit$] For every integer $1 \leq m \leq n$ and for every $U \in\cJ_{\GL_{m, F}, \prod_{v \not\in S} \GL_m(\cO_{F_v})}$, the complex $R \Gamma_{X_{\GL_{m, F}}^U}k$ is of $S$-Galois type.

\end{itemize}
If $F$ is an imaginary CM or totally real field, then $\spadesuit$ holds, by Corollary \ref{cor_intro_scholze_result}.
\begin{lemma}\label{lem_compactly_supported_coh_of_galois_type}
This hypothesis is equivalent to: for every integer $1 \leq m \leq n$ and for every open compact subgroup $U \in\cJ_{\GL_{m, F}, \prod_{v \not\in S} \GL_m(\cO_{F_v})}$, the complex $R \Gamma_{X_{\GL_{m, F}}^U, c}k$ is of $S$-Galois type.
\end{lemma}
\begin{proof}
We show that our hypothesis implies the given condition on the cohomology with compact support; the other direction is similar. We can assume that $m = n$. All maximal ideals occur in the support of cohomology groups. By Proposition \ref{prop_hecke_action_and_verdier_duality}, there is a perfect Poincar\'e duality pairing of finite-dimensional $k$-vector spaces
\[ \langle \cdot, \cdot \rangle_U : H^\ast_c(X_G^U, k) \times H^\ast(X_G^U, k) \to k \]
satisfying the equation $\langle x, [U g U] y \rangle_U = \langle [U g^{-1} U]x, y \rangle_U$ for any $g \in G^S$. For unramified Hecke operators $T_v^i$, $v \not\in S$, this implies that the action of $T_v^i$ on $H^\ast(X_G^U, k)$ is dual to the action of $T_v^{m-i} (T_v^m)^{-1}$ on $H^\ast_c(X_G^U, k)$. We must therefore show that for any maximal ideal $\ffrm$ of $\bbT_G^S$ in the support of $H^\ast(X_G^U, k)$, there exists a continuous semi-simple representation $\overline{\sigma}_\ffrm : G_{F, S} \to \GL_m(\bbT^S(H^\ast(X_G^U, k))/\ffrm)$ such that for each finite place $ v\not\in S$ of $F$, $\overline{\sigma}_\ffrm(\Frob_v)$ has characteristic polynomial
\[ X^m - T_v^{m-1} (T_v^m)^{-1} X^{m-1} + \dots + (-1)^j q_v^{j(j-1)/2} T_v^{m-j} (T_v^m)^{-1} X^{m-j} + \dots + (-1)^m q_v^{m(m-1)/2} (T_v^m)^{-1}. \]
A calculation shows that we can take $\overline{\sigma}_\ffrm \cong \overline{\rho}_\ffrm^\vee \otimes \epsilon^{1-m}$.
\end{proof}
Subject to $\spadesuit$, we will prove the following theorem:
\begin{theorem}\label{thm_boundary_cohomology_is_eisenstein}
For every $U \in \cJ_{G, U^S}$ and for every smooth $\cO[U_S]$-module $A$, finite as $\cO$-module, the complex $R \Gamma_{\partial \overline{X}_G^U}(\underline{A}_G^U)$ is Eisenstein. In particular, for every non-Eisenstein maximal ideal $\ffrm \subset \bbT^S(R \Gamma_{X_G^U} \underline{A}_G^U)$, the natural morphism 
\[ (R \Gamma_{X_G^U, c} \underline{A}_G^U)_\ffrm \to (R \Gamma_{X_G^U} \underline{A}_G^U)_\ffrm \]
in $\mathbf{D}(\cO)$ is a quasi-isomorphism.
\end{theorem}
The proof of Theorem \ref{thm_boundary_cohomology_is_eisenstein} will be an exercise in understanding the structure of the Borel--Serre boundary. We begin with some preliminary reductions.
\begin{lemma}
To prove Theorem \ref{thm_boundary_cohomology_is_eisenstein}, it is enough to treat the case where $A = k$ is the trivial $k[U_S]$-module.
\end{lemma}
\begin{proof}
Since $A$ is a finite $\cO$-module, we can find a normal subgroup $U' \subset U$ such that $U' \in \cJ_{G, U^S}$ and $U'_S$ acts trivially on $A$. We will show that $\Supp_{\bbT^S} R \Gamma_{X_G^U} \underline{A}_G^U \subset \Supp_{\bbT^S} R \Gamma_{X_G^{U'}} k$. There is a $\bbT^S_G$-equivariant spectral sequence
\[ E_2^{i, j} = H^i(U/U', H^j(X_G^{U'}, \underline{A}_G^{U'})) \Rightarrow H^{i+j}(X_G^U, \underline{A}_G^U), \]
which shows that $\Supp_{\bbT^S} R \Gamma_{X_G^U} \underline{A}_G^U \subset \Supp_{\bbT^S} R \Gamma_{X_G^{U'}} \underline{A}_G^{U'}$. But $\underline{A}_G^{U'}$ is the constant sheaf associated to a finite $\cO$-module, and the result follows by the theorem on universal coefficients.
\end{proof}
\begin{lemma}
To prove Theorem \ref{thm_boundary_cohomology_is_eisenstein}, it is enough to show that for each proper standard parabolic subgroup $P \subset G$, the complex $R \Gamma_{[ \Ind_{P^\infty}^{G^\infty} \mathfrak{X}_P ]/U} k$ is Eisenstein.
\end{lemma}
The same argument shows that we are free to replace $U_S$ by any normal open compact subgroup. 
\begin{proof}
By Poincar\'e duality (i.e.\ by the same argument as in the proof of Lemma \ref{lem_compactly_supported_coh_of_galois_type}), the vanishing of $R \Gamma_{[ \Ind_{P^\infty}^{G^\infty} \mathfrak{X}_P ]/U} k$ at non-Eisenstein maximal ideals implies that of the compactly supported cohomology $R \Gamma_{[ \Ind_{P^\infty}^{G^\infty} \mathfrak{X}_P ]/U, c} k$. By Lemma \ref{lem_boundary_stratification_by_induced_strata} and the long exact sequence in cohomology with compact supports associated to the inclusion of an open subspace, we deduce the corresponding result for the full boundary $\partial \overline{X}_G^U$.
\end{proof}
Let us therefore fix a proper partition $n = n_1 + \dots + n_s$, and let $P \subset G$ denote the corresponding standard parabolic subgroup, $M \cong \GL_{n_1} \times \dots \times \GL_{n_s}$ its standard Levi subgroup. We will now show that $R \Gamma_{[ \Ind_{P^\infty}^{G^\infty} \mathfrak{X}_P ]/U} k$ is Eisenstein. Since every rational proper parabolic subgroup of $G$ is conjugate to one of this form, this will complete the proof of Theorem \ref{thm_boundary_cohomology_is_eisenstein}. 

We isomorphisms in $\mathbf{D}(\cH(G^S, U^S))$:
\begin{equation}\label{eqn_boundary_coh_1} \begin{split} R \Gamma_U \Res^{G^\infty}_{G^S \times U_S} R \Gamma_{\Ind_{P^\infty}^{G^\infty} \mathfrak{X}_P} k & \cong  R \Gamma_U \Res^{G^\infty}_{G^S \times U_S} H^0(P(F) \backslash G^\infty, k) \\
 & \cong R \Gamma_{U^S} H^0(P(F) \backslash G^\infty / U_S, k) \\ & \cong  R \Gamma_{U^S} \Ind_{P^S}^{G^S} H^0(P(F) \backslash P^S \times G_S / U_S, k). \end{split}
\end{equation}
Our assumptions imply that $P^S U_P^S = G^S$. We can therefore apply Corollary \ref{cor_derived_hecke_morphism_for_parabolic_induction} to deduce that the complex in (\ref{eqn_boundary_coh_1}) is quasi-isomorphic with 
\begin{equation}\label{eqn_boundary_coh_2} \begin{split} \bigoplus_{g \in P(F) \backslash G_S / U_S} r_P^\ast R \Gamma_{U_P^S} & \Ind^{P^S}_{P(F) \cap g U_S g^{-1}} k \\ &\cong \bigoplus_{g \in P(F) \backslash G_S / U_S} r_P^\ast R \Gamma_{U_P^S \times (P_S \cap g U_S g^{-1})} \Res^{P^\infty}_{P^S \times (P_S \cap g U_S g^{-1})} R \Gamma_{\mathfrak{X}_P} k. \end{split}
\end{equation}
The index set in the direct sum is finite because the quotient $P_S \backslash G_S / U_S$ is finite and $P(F)$ is dense in $P_S$ ($P$ is an $F$-rational variety). By varying $U_S$, we can therefore reduce to showing that the complex $r_P^\ast R \Gamma_{U_P} \Res^{P^\infty}_{P_S \times U_{P, S}} R \Gamma_{\mathfrak{X}_P} k$ is Eisenstein. To show this, we write $\pi : \mathfrak{X}_P \to \mathfrak{X}_M$ for the canonical projection and observe that there is a quasi-isomorphism in $\Mod(\cH(P^S, U_P^S))$:
\begin{equation} R{\Gamma_{U_P}} R \Gamma_{\Res^{P^\infty}_{P^S \times U_{P, S}}\mathfrak{X}_P} k \cong R \Gamma_{U_P^S \times U_{M, S}} R\Gamma_{\Inf_{M^S \times U_{M,S}}^{P^S \times U_{M, S}} \mathfrak{X}_M} R \pi_\ast^{U_{N, S}} k. 
\end{equation}
The sheaf $\pi_\ast k \in \Sh_{P^S \times U_{P, S}} \mathfrak{X}_M$ is constant on connected components, with stalk at a point $[(\overline{p}, x)] \in \mathfrak{X}_M$ given by the formula
\begin{equation} (\pi_\ast k)_{[\overline{p}, x]} = H^0(P(F) \backslash P(F) p N(\bbA^\infty_F), k) \cong \Ind_{ N(F)}^{N^\infty} k, \gamma pn \mapsto p n p^{-1}.
\end{equation}
By strong approximation, there is an isomorphism of $\bbZ[U_N]$-modules, where $\Gamma_{1, U_N} = N(F) \cap U_N$ as usual:
\[ \Ind_{N(F)}^{N^\infty} k \cong \Ind_{\Gamma_{1, U_N}}^{U_N} k. \]
Since $U_{N, S}$ acts freely on the set $\Gamma_{1, U_N} \backslash U_N$, $\Ind_{\Gamma_{1, U_N}}^{U_N} k$ is an injective $\bbZ[U_{N, S}]$-module, and the natural map $\pi_\ast^{U_{N, S}} k \to R \pi_\ast^{U_{N, S}}$ in $\mathbf{D}(\Sh_{P^S \times U_{M, S}}\mathfrak{X}_M)$ is a quasi-isomorphism. We deduce the existence of a quasi-isomorphism 
\begin{equation}\label{eqn_higher_derived_invariants_vanish} R{\Gamma_{U_P}} R \Gamma_{\Res^{P^\infty}_{P^S \times U_{P, S}}\mathfrak{X}_P} k \cong R{\Gamma_{U_P^S \times U_{M, S}}} R \Gamma_{\Inf_{M^S \times U_{M, S}}^{P^S \times U_{M, S}} \mathfrak{X}_M} \pi_\ast^{U_{N, S}} k. 
\end{equation}
We now construct a morphism $k \to \pi_\ast k$ in $\Sh_{P^S \times U_{P, S}} \mathfrak{X}_M$. Since the constant sheaf $k$ is pulled back from a point, it is equivalent to give a $P^S \times U_{P, S}$-equivariant map
\[ k \to H^0(\mathfrak{X}_M, \pi_\ast k) = H^0(\mathfrak{X}_P, k) \cong H^0(P(F) \backslash P^\infty, k), \]
which we take to be the inclusion of the constant functions. Taking derived $U_{N, S}$-invariants, we obtain a natural map $R 1_\ast^{U_{N, S}} k \to R \pi_\ast^{U_{N, S}} k \cong \pi_\ast^{U_{N, S}} k$, hence a map
\begin{equation}  R \Gamma_{U_P^S \times U_{M, S}} R \Gamma_{\Inf_{M^S \times U_{M, S}}^{P^S \times U_{M, S}} \mathfrak{X}_M} R 1_\ast^{U_{N, S}} k \to  R \Gamma_{U_P^S \times U_{M, S}} R \Gamma_{\Inf_{M^S \times U_{M, S}}^{P^S \times U_{M, S}} \mathfrak{X}_M} \pi_\ast^{U_{N, S}} k \cong  R{\Gamma_{U_P}} R \Gamma_{\Res^{P^\infty}_{P^S \times U_{P, S}}\mathfrak{X}_P} k,
\end{equation}
the last isomorphism by (\ref{eqn_higher_derived_invariants_vanish}). We claim that this is a quasi-isomorphism of complexes of $\cH(P^S, U_P^S)$-modules. It suffices to check this after applying forgetful functors, which reduces us to showing that the natural map 
\[ R 1_\ast^{U_N^S} R 1_\ast^{U_{N, S}} k \cong R 1_\ast^{U_N} k \to R 1_\ast^{U_{N, S}} \pi_\ast^{U_{N, S}} k = R \pi_\ast^{U_N} k \]
is a quasi-isomorphism. After taking cohomology and looking at stalks, we must show that the maps
\[ H^i(U_N, k) \to H^i(U_N, \Ind_{\Gamma_{1, U_N}}^{U_N} k) \cong H^i(\Gamma_{1, U_N}, k) \]
are isomorphisms. This is part of the following lemma.
\begin{lemma}\label{lem_dubious_homotopy_lemma_for_nilpotent_groups}
Let $\Gamma_{1, U_{N, S}} = N(F) \cap U_{N, S}$ (intersection inside $N_S$). Then the natural maps in (discrete) group cohomology
\[ H^\ast(U_{N, S}, k) \to H^\ast(\Gamma_{1, U_{N, S}}, k) \to H^\ast(\Gamma_{1, U_N}, k) \]
are all isomorphisms, while $H^i(U_N^S, k) = 0$ for $i > 0$.
\end{lemma}
\begin{proof}
Let $S_p \subset S$ denote the set of places dividing $p$. Nilpotent groups have the congruence subgroup property, so the natural map $\Gamma_{1, U_N} \to U_{N, S_p}$ identifies $U_{N, S_p}$ with the $p$-profinite completion of $\Gamma_{1, U_N}$, and hence (\cite[Ch. VI, 5.6]{Bou72}) the natural map $H^\ast(U_{N, S_p}, \bbF_p) \to H^\ast(\Gamma_{1, U_N}, \bbF_p)$ is an isomorphism. Let $T$ be the set of places of $S$ which are prime to $p$. Then the group $U_{N, T}$ is uniquely $p$-divisible, hence is $\bbZ[1/p]$-complete in the terminology of \emph{op. cit.}, hence satisfies $H^i(U_{N, T}, \bbZ) = H^i(U_{N, T}, \bbZ[1/p])$ and $H^i(U_{N, T}, \bbF_p) = 0$ for $i > 0$, by \cite[Ch. V, 3.3]{Bou72}. The K\"unneth theorem in group cohomology then implies that the natural map
\[ H^\ast(U_{N, S}, \bbF_p) \to H^\ast(\Gamma_{1, U_N}, \bbF_p) \]
is an isomorphism. Essentially the same argument shows that the same holds when $\Gamma_{1, U_N}$ is replaced by $\Gamma_{1, U_{N, S}}$. The group $U_{N}^S$ is uniquely $p$-divisible and nilpotent, so another application of \cite[Ch. V, 3.3]{Bou72} shows that it has trivial $\bbF_p$-cohomology.
\end{proof}
Let us now take stock. We have shown that there is a quasi-isomorphism
\[  R{\Gamma_{U_P}} R\Gamma_{\Res^{P^\infty}_{P^S \times U_{P, S}}\mathfrak{X}_P} k \cong R \Gamma_{U_P^S \times U_{M, S}} R \Gamma_{\Inf_{M^S \times U_{M, S}}^{P^S \times U_{M, S}} \mathfrak{X}_M} R 1_\ast^{U_{N, S}} k \]
of complexes of $\cH(P^S, U_{P}^S)$-modules, and we wish to show that these complexes become Eisenstein after applying the exact functor $r_P^\ast$. It is enough to show that for each $i \geq 0$, the complex
\[ r_P^\ast R \Gamma_{U_P^S \times U_{M, S}} R \Gamma_{\Inf_{M^S \times U_{M, S}}^{P^S \times U_{M, S}} \mathfrak{X}_M} R^i 1_\ast^{U_{N, S}} k \]
is Eisenstein. The sheaf $R^i 1_\ast^{U_{N, S}} k \in \Sh_{P^S \times U_{M, S}} \mathfrak{X}_M$ can be calculated explicitly as follows: it is pulled back from the sheaf on a point associated to the $k[U_{M, S}]$-module $A = H^i(U_{N, S}, k)$. In the remainder of this section, we will show that there is a quasi-isomorphism
\[ R \Gamma_{U_P^S \times U_{M, S}} R \Gamma_{\Inf_{M^S \times U_{M, S}}^{P^S \times U_{M, S}} \mathfrak{X}_M} R^i 1_\ast^{U_{N, S}} k \cong r_M^\ast R \Gamma_{U_M} R \Gamma_{\mathfrak{X}_M} \underline{A}_M \]
of complexes of $\cH(P^S, U_P^S)$-modules, and that these last complexes become Eisenstein after applying the exact functor $r_P^\ast$. This will complete the proof of Theorem \ref{thm_boundary_cohomology_is_eisenstein}.

There is a natural morphism 
\[ r_M^\ast R \Gamma_{U_M} R \Gamma_{\mathfrak{X}_M} \underline{A}_M  \to R \Gamma_{U_P \times U_{M, S}} R \Gamma_{\Inf_{M^S \times U_{M, S}}^{P^S \times U_{M, S}} \mathfrak{X}_M} \underline{A}_M, \]
by Corollary \ref{cor_universality_of_derived_functor}. It is a quasi-isomorphism; indeed, we can check this after applying the exact forgetful functor $(\cdot)^\sim$, which reduces us to showing that the natural map
\[ \underline{A}_M \to R 1_\ast^{U_N^S} \Inf_{U_M^S \times U_{M, S}}^{U_P^S \times U_{M, S}} \underline{A}_M \]
of complexes of sheaves in $\mathbf{D}(\Sh_{M^S \times U_{M, S}} \mathfrak{X}_M)$ is a quasi-isomorphism. This can be checked on stalk cohomology, where it reduces to the assertion that $H^i(U_N^S, k) = 0$ if $i > 0$, which is part of Lemma \ref{lem_dubious_homotopy_lemma_for_nilpotent_groups}. 

We now show that the complex
\[ r_P^\ast r_M^\ast   R \Gamma_{U_M} R \Gamma_{\mathfrak{X}_M} \underline{A}_M = \cS^\ast R \Gamma_{U_M} R \Gamma_{\mathfrak{X}_M} \underline{A}_M \]
 is Eisenstein. After possibly shrinking $U$, we can assume that the action of $U_{M, S}$ on $H^i(U_{N, S}, k)$ induced by conjugation is trivial, implying that the sheaf $\underline{A}_M$ is in fact constant on $\mathfrak{X}_M$. This reduces us to showing that the complex 
\[ r_P^\ast r_M^\ast R \Gamma_{U_M} R \Gamma_{\mathfrak{X}_M} k \]
 is Eisenstein. After further shrinking $U$, we can assume that $U_M = U_1 \times \dots \times U_s$ for neat open compact subgroups $U_i \subset \GL_{n_i}(\bbA_F^\infty)$. In this case, we have a commutative diagram
\[ \xymatrix{\cH(G^S, \prod_{v\not\in S} \GL_n(\cO_{F_v})) \ar[r] & \cH(M^S, \prod_{v\not\in S} \prod_{i=1}^s \GL_{n_i}(\cO_{F_v})) \ar[r] \ar[d] & \End_k(H^\ast(X_M^{U_M}, k)) \ar[d] \\
& \otimes_{i=1}^s \cH(\GL_{n_i}^S, \prod_{v \not\in S} \GL_{n_i}(\cO_{F_v})) \ar[r] &  \End_k(\otimes_{i=1}^s H^\ast(X_{\GL_{n_i}}^{U_i}, k)). } \]
We now use the following lemma:
\begin{lemma}\label{lem_calculation_of_satake_transform_for_GL_n}
For each place $v \not\in S$ of $F$, there is a commutative diagram
\[ \xymatrix{ \cH( \GL_n(F_v), \GL_n(\cO_{F_v}) )  \ar[r] \ar[d] & \bbR[ Y_1^{\pm 1}, \dots, Y_n^{\pm 1} ]^{S_n} \ar[d] \\
\otimes_{i=1}^s \cH( \GL_{n_i}(F_v), \GL_{n_i}(\cO_{F_v}) )  \ar[r] & \otimes_{i=1}^s \bbR[Z_{i, 1}^{\pm 1}, \dots, Z_{i, n_i}^{\pm 1}]^{S_{n_i}}. } \]
The horizontal arrows are induced by the usual (normalized) Satake isomorphisms. The left vertical arrow is the unnormalized Satake transform $\cS_v = r_M \circ r_P$. The right hand arrow is defined by the bijective correspondence for each $i = 1, \dots, s$:
 \[ \{ Y_{n_1 + \dots + n_{i-1} + 1}, \dots, Y_{n_1 + \dots + n_i} \} \leftrightarrow q_v^{(n_{i+1} + \dots + n_{s}) - (n_1 + \dots + n_{i-1})} \{ Z_{i, 1}, \dots, Z_{i, n_i} \}. \]
\end{lemma}
\begin{proof}
 The proof is very similar to the proof of Proposition-Definition \ref{prop_calculation_of_unnormalized_satake_transform_in_unitary_case} below, and is therefore omitted.
\end{proof}
We can now complete the proof of Theorem \ref{thm_boundary_cohomology_is_eisenstein}. If $\ffrm$ is a maximal ideal of $\bbT^S_M$ in the support of $H^\ast(X_M^{U_M}, k)$ with residue field $k$ (which we can always assume after possibly enlarging the field of scalars), then we can find for each $i = 1, \dots, s$ a maximal ideal $\ffrm_i$ of $\bbT^S_{\GL_{n_i}}$ with residue field $k$ and appearing in the support of $H^\ast(X_{\GL_{n_i}}^{U_i}, k)$ such that $\ffrm$ is the product of $\ffrm_1, \dots, \ffrm_s$, in the obvious sense. 

Applying finally our hypothesis $\spadesuit$, we can find for each $i = 1, \dots, s$ a semi-simple Galois representation
\[ \overline{\rho}_{\ffrm_i} : G_{F, S} \to \GL_{n_i}( \bbT^S_{\GL_{n_i}} / \ffrm_i) \]
such that for each finite place $v\not\in S$ of $F$, $\overline{\rho}_{\ffrm_i}(\Frob_v)$ has characteristic polynomial
\[ X^{n_i} - T_v^1 X^{n_i-1} + \dots  + (-1)^{n_i} q_v^{n_i(n_i-1)/2} T_v^{n_i} \in (\bbT^S_{\GL_{n_i}} / \ffrm_i)[X]. \]
If $\cS^\ast(\ffrm)$ denotes the pullback of $\ffrm$ to $\bbT^S_G$, we see that $\cS^\ast(\ffrm)$ is in the support of $H^\ast(X_M^{U_M}, k)$ as $\bbT^S_G$-module, and (using Lemma \ref{lem_calculation_of_satake_transform_for_GL_n} and the fact that the normalized Satake isomorphism for $\GL_n$ is characterized by the formula $T_v^i \mapsto q^{i(n-i)/2} e_i(Y_1, \dots, Y_n)$, where $e_i$ denotes the $i^\text{th}$ symmetric polynomial) that the Galois representation
\[ \overline{\rho}_{\cS^\ast(\ffrm)} = \bigoplus_{i=1}^s \overline{\rho}_{\ffrm_i} \otimes \epsilon^{-(n_{i+1} + \dots + n_s)} \]
satisfies the desired relation (\ref{eqn_characteristic_polynomial_of_Frobenius}). We observe that this Galois representation is reducible, by construction. Since every maximal ideal of $\bbT^S_G$ which is in the support of $H^\ast(X_M^{U_M}, k)$ is of the form $\cS^\ast(\ffrm)$ for some $\ffrm \subset \bbT^S_M$, this shows that $\cS^\ast H^\ast(X_M^{U_M}, k)$ is Eisenstein, as desired.

\section{The boundary cohomology of the $U(n, n)$ locally symmetric space}\label{sec_main_arguments}

In this section, we will prove the main theorems of this paper. 
\subsection{Groups and local systems}
Let $F$ be an imaginary CM number field with totally real subfield $F^+$, and let $c \in \Gal(F/F^+)$ denote complex conjugation. Let $d = [F^+ : \bbQ]$, and let $n \geq 2$ be an integer. Let $\Psi_n$ denote the $n \times n$ matrix with 1's on the anti-diagonal and 0's elsewhere, and 
\[ J_n = \left( \begin{array}{cc} 0 & \Psi_n \\ - \Psi_n & 0 \end{array} \right). \]
Then $J_n$ defines a perfect Hermitian pairing $\langle \cdot, \cdot \rangle : \cO_F^{2n} \times \cO_F^{2n} \to \cO_F$, given by the formula $\langle x, y \rangle = {}^t x J_n y^c$. We write $\underline{G}$ for the group over $\cO_{F^+}$ given by the formula for any $\cO_{F^+}$-algebra $R$:
\[ \underline{G}(R) = \{ g \in \GL_{2n}(\cO_F \otimes_{\cO_{F^+}} R) \mid {}^t g J_n g^c = J_n \}. \]
We write $\underline{P} \subset \underline{G}$ for the closed subgroup which leaves invariant the subspace $\cO_F^{n} \oplus 0^n \subset \cO_F^{2n}$, and $\underline{M} \subset \underline{P}$ for the closed subgroup which leaves invariant the direct factors $\cO_F^n \oplus 0^n$ and $0^n \oplus \cO_F^n$. We write $\underline{T} \subset \underline{G}$ for the standard diagonal torus and $\underline{B} \subset \underline{G}$ for the standard upper-triangular subgroup. We write $\underline{S} \subset \underline{T}$ for the subtorus consisting of matrices with elements in $\cO_{F^+}$. We will denote base change to the $F$-fibre of these groups by omitting the underline.

Thus $P$ is a parabolic subgroup of $G$, and $M$ is the unique Levi subgroup of $P$ containing $T$. The torus $S$ is a maximal $F^+$-split torus of $G$, $T = Z_G(S)$, and $G = U(n, n)$ is quasi-split. We have the equalities
\[ \dim_\bbR X_G = 2 d n^2,\text{ } \dim_\bbR X_P = 2 d n^2 - 1,\text{ } \dim_\bbR X_M = d n^2 - 1. \]
We set $D = d n^2$. 
\begin{lemma} Let notation be as above.
\begin{enumerate}
\item If $v$ is a place of $F^+$ which is unramified in $F$, then $\underline{G}_{\cO_{F^+_v}}$ is reductive, and hence $G_{F^+_v}$ is unramified.
\item Let $\underline{N} \subset \underline{P}$ denote the closed subgroup which acts trivially on the factors $\cO_F^n \oplus 0^n$ and $0^n \oplus \cO_F^n$. Then $\underline{P} \cong \underline{M} \ltimes \underline{N}$. There is an isomorphism $\underline{M} \cong \Res^{\cO_F}_{\cO_{F^+}} \GL_n$. 
\end{enumerate}
\end{lemma}
\begin{proof}
The first part follows easily from the definitions; indeed, one can check that if $w$ is a place of $F$ lying above $F^+$, then $\underline{G}_{\cO_{F_w}}$ is isomorphic to $\GL_{2n}$. For the second part, we make things explicit. Let $(\cdot)^\ast$ denote the anti-involution of $\Res^{\cO_F}_{\cO_{F^+}} \GL_n$ given by $A^\ast = \Psi_n {}^t A^c \Psi_n^{-1}$. (Explicitly, $A^\ast$ is given by reflection in the anti-diagonal of $A$ and conjugation of coefficients.) Then $\underline{P}$ can be identified with the subgroup of $\Res^{\cO_F}_{\cO_{F^+}} \GL_{2n}$ consisting of matrices of the form
\[ g = \left( \begin{array}{cc} D^{-\ast} & B \\ 0 & D \end{array} \right), \]
with $B = B^\ast$ and no condition on $D$. The subgroup $\underline{N} \subset \underline{P}$ is given by the condition $D = 1$, and the subgroup $\underline{M}$ by the condition $B = 0$. It is easy to see that the natural map $\underline{M} \times \underline{N} \to \underline{P}$ given by multiplication of components is an isomorphism. We identify $\underline{M}$ with $\Res^{\cO_F}_{\cO_{F^+}} \GL_n$ via the map $g \mapsto D$.
\end{proof}
For each place $v$ of $F^+$ which is unramified in $F$, the group $U_v = \underline{G}(\cO_{F^+_v}) \subset G(F^+_v)$ is a hyperspecial maximal compact subgroup. Moreover, the subgroups $M$ and $P = M \ltimes N$ satisfy the conditions of Lemma \ref{lem_hecke_restriction_to_parabolic} and Lemma \ref{lem_hecke_integration_along_fibres} with respect to $U_v$, and $U_{M, v} = \underline{M}(\cO_{F^+_v}) = \GL_n(\cO_F \otimes_{\cO_{F^+}} \cO_{F_v^+})$. We thus have the map 
\begin{equation}\label{eqn_recall_unnormalized_satake_map} \cS_v = r_{M_v} \circ r_{P_v} : \cH(G(F_v^+), U_v) \to \cH(M(F_v^+), U_{M, v}). 
\end{equation}
In our situation, it can be given explicitly as in the following two propositions.
\begin{propdef}\label{prop_calculation_of_G_Hecke_algebras}
Let $v$ be a place of $F^+$ which is unramified in $F$, and let $w$ be a place of $F$ dividing $v$.
\begin{enumerate}
\item Suppose that $v$ splits in $F$. Then $G(F_v^+) \cong \GL_{2n}(F_w)$, and the Satake isomorphism gives a canonical isomorphism
\[ \cH(G(F_v^+), U_v) \otimes_\bbZ \bbR \cong \bbR[Y_1^{\pm 1}, \dots, Y_{2n}^{\pm 1}]^{S_{2n}}. \]
For each $i = 1, \dots, 2n$, we write $T_{G, w, i}$ for the element of $\cH(G(F_v^+), U_v) \otimes_\bbZ \bbZ[q_v^{-1}]$ which corresponds under the Satake transform to the element $q_w^{i(2n-i)/2} e_i(Y_1, \dots, Y_{2n})$.
\item Suppose that $v$ is inert in $F$. The Satake isomorphism gives a canonical isomorphism
\[ \cH(G(F_v^+), U_v) \otimes_\bbZ \bbR \cong \bbR[X_1^{\pm 1}, \dots, X_n^{\pm 1}]^{S_n \ltimes (\bbZ / 2 \bbZ)^n}. \]
The unramified endoscopic transfer from $G(F_v^+)$ to $\GL_{2n}(F_w)$ is dual to the map
\[ \bbR[Y_1^{\pm 1}, \dots, Y_{2n}^{\pm 1}]^{S_{2n}} \to \bbR[X_1^{\pm 1}, \dots, X_n^{\pm 1}]^{S_n \ltimes (\bbZ/2\bbZ)^n} \]
which puts the set $\{ Y_1, \dots, Y_{2n} \}$ in bijection with $\{ X_1^{\pm 1}, \dots, X_n^{\pm 1} \}$. For each $i = 1, \dots, 2n$, we write $T_{G, w, i}$ for the element of $\cH(G(F_v^+), U_v) \otimes_\bbZ \bbZ[q_v^{-1}]$ which corresponds under the Satake transform to the image of $q_w^{i(2n-i)/2} e_i(Y_1, \dots, Y_{2n})$.
\end{enumerate}
\end{propdef}
\begin{proof}
For concreteness, we recall the definition of the normalized Satake isomorphism. We temporarily let notation be as at the beginning of \S \ref{sec_application_when_G_unramified}. Thus $F$ is a finite extension of $\bbQ_p$, $\underline{G}$ is a reductive group over $\cO_F$, $\underline{S}$ is a maximal $\cO_F$-split torus of $\underline{G}$, $\underline{T}$ is the maximal torus which centralizes $\underline{S}$, and $\underline{B}$ is a Borel subgroup containing $\underline{T}$. Let $\underline{N} \subset \underline{B}$ denote the unipotent radical. The Satake isomorphism is then the isomorphism (\cite[p. 147]{Car79}):
\[ \cN : \cH( G(F), \underline{G}(\cO_F) ) \otimes_\bbZ \bbR \to \cH( T(F), \underline{T}(\cO_F) )^{W(G, T)} \otimes_\bbZ \bbR \]
given by the formula $f \mapsto (t \mapsto \delta_B(t)^{1/2} \int_{n \in N(F)} f(t n) dn)$. (We use the notation $\cN$ to emphasize that this is the `normalized' Satake isomorphism, in contrast to the `unnormalized' Satake transform $\cS$ that we hvae used elsewhere in this paper, in which the factor $\delta_B^{1/2}$ does not appear.) Here $\delta_B(t)$ is the character  $\delta_B : T(F) \to \bbR_{>0}$ given by the formula $\delta_B(t) = | \det_F \mathrm{Ad}_{\Lie N}(t) |_F$, where $|\cdot|_F$ is the usual normalized absolute value on $F$ (satisfying the formula $|\varpi|_F = (\# k_F)^{-1} = q^{-1}$ for $\varpi \in \cO_F$ a uniformizer).

It follows from the proof given in \emph{loc.~cit.~}that $\cN$ in fact defines an isomorphism 
\[ \cN : \cH( G(F), \underline{G}(\cO_F) ) \otimes_\bbZ \bbZ[q^{\pm \frac{1}{2}}] \cong \cH( T(F), \underline{T}(\cO_F) )^{W(G, T)} \otimes_\bbZ  \bbZ[q^{\pm \frac{1}{2}}]. \]
If the character $\delta_B$ takes values in $q^{2 \bbZ}$, then it even defines an isomorphism over $\bbZ[q^{-1}]$. More generally, if $\rho_G \in X^\ast(G)$ is a character such that $t \mapsto \delta_B(t)^{1/2} |\rho_G(t)|_F^{1/2}$ takes values in $q^{\bbZ}$, then we get an isomorphism
\[ \cN' :  \cH( G(F), \underline{G}(\cO_F) ) \otimes_\bbZ \bbZ[q^{-1}] \to \cH( T(F), \underline{T}(\cO_F) )^{W(G, T)} \otimes_\bbZ \bbZ[q^{-1}]\]
given by the formula $\cN'(f) = (t \mapsto |\rho_G(t)|_F^{1/2} \cN(f)(t) )$.

We now return to the notation of the proposition. To complete the first part, we must check that the element $T_{G, w, i}$, which a priori lies in $\cH( G(F_v^+), U_v ) \otimes_\bbZ \bbZ[q_v^{\pm \frac{1}{2}}]$, in fact lies in $\cH( G(F_v^+), U_v ) \otimes_\bbZ \bbZ[q_v^{-1}]$. We could use the stronger fact, used already in \S \ref{sec_i_dont_see_any_boundary}, that $T_{G, w, i}$ is actually equal to one of the canonical basis elements of $\cH( G(F_v^+), U_v )$. Alternatively, we can apply the above formalism with $\rho_G = \det^{2n-1}$. Then we find that
\[ \begin{split} \cN'(T_{G, w, i}) & = | \varpi_w^{i(2n-1)}|_{F_w} q_w^{i(2n-i)/2} e_i(Y_1, \dots, Y_{2n})\\ & = q_w^{-i(i-1)/2} e_i(Y_1, \dots, Y_{2n}) \in  \cH( T(F), \underline{T}(\cO_F) )^{W(G, T)} \otimes_\bbZ \bbZ[q_w^{-1}], \end{split} \]
hence $T_{G, w, i} \in \cH( G(F_v^+), U_v ) \otimes_\bbZ \bbZ[q_v^{-1}]$.

We now prove the second part of the proposition. By definition, the unramified endoscopic transfer is the map on unramified Hecke algebras dual to the standard unramified base change map defined, for example, in \cite[\S 4.1]{Min11}. This is easily seen to correspond under the respective Satake isomorphisms to the map appearing in the statement of the proposition above. To finish the proof, we must again show that $T_{G, w, i} \in \cH( G(F_v^+), U_v ) \otimes_\bbZ \bbZ[q_v^{-1}]$. We observe that since $q_w = q_v^2$, the image of $q_w^{i(2n-i)/2} e_i(Y_1, \dots, Y_{2n})$ lies in $\cH( T(F), \underline{T}(\cO_F) )^{W(G, T)} \otimes_\bbZ \bbZ[q_v^{-1}]$. It is easy to check that for the unramified unitary group, the character $\delta_B$ takes values in $q_v^{2 \bbZ}$, and hence the normalized Satake isomorphism is itself defined over $\bbZ[q_v^{-1}]$. These facts together imply the result.
\end{proof}
\begin{propdef}\label{prop_calculation_of_unnormalized_satake_transform_in_unitary_case} Let $v$ be a place of $F^+$ which is unramified in $F$, and let $w$ be a place of $F$ dividing $v$.
\begin{enumerate}
\item Suppose that $v$ splits in $F$. The unnormalized Satake transform (\ref{eqn_recall_unnormalized_satake_map}) corresponds under the Satake isomorphism to the map
\[ \bbR[Y_1^{\pm 1}, \dots, Y_{2n}^{\pm 1}]^{S_{2n}} \to \bbR[W_1^{\pm 1}, \dots, W_n^{\pm 1}, Z_1^{\pm 1}, \dots, Z_n^{\pm 1}]^{S_n \times S_n} \]
which puts the set $\{ Y_1, \dots, Y_{2n} \}$ in bijection with $\{ q_v^{n/2} Z_n^{-1}, \dots, q_v^{n/2} Z_1^{-1}, q_v^{-n/2} W_1, \dots, q_v^{-n/2} W_n \}$. For each $i = 1, \dots, n$, let $T_{M, w, i} \in \cH(M(F_v^+), U_{M, v}) \otimes_\bbZ \bbZ[q_v^{-1}]$ correspond to $q_w^{i(n-i)/2} e_i(W_1, \dots, W_n)$, and let $T_{M, w^c, i} \in \cH(M(F_v^+), U_{M, v})$ correspond to $q_w^{i(n-i)/2} e_i(Z_1, \dots, Z_n)$. 
\item Suppose instead that $v$ is inert in $F$. Then the unnormalized Satake transform (\ref{eqn_recall_unnormalized_satake_map}) corresponds under the Satake isomorphism to the map
\[ \bbR[X_1^{\pm 1}, \dots, X_n^{\pm 1}]^{S_n \ltimes (\bbZ/2\bbZ)^n}  \to \bbR[W_1^{\pm 1}, \dots, W_n^{\pm 1}]^{S_n} \]
which puts the set $\{ X_1, \dots, X_n \}$ in bijection with the set $\{ q_w^{-n/2} W_1, \dots, q_w^{-n/2} W_n \}$. For each $i = 1, \dots, n$, we let $T_{M, w, i} \in \cH(M(F_v^+), U_{M, v}) \otimes_\bbZ \bbZ[q_v^{-1}]$ correspond under the Satake isomorphism to the element $q_w^{i(n-i)/2} e_i(W_1, \dots, W_n)$. 
\end{enumerate}
\end{propdef}
\begin{proof}
In either case, we have a diagram
\[ \xymatrix{ \cH(G(F_v^+), U_v) \otimes_\bbZ \bbR \ar[r]^{\cS_v} \ar[rd]_{\cN_G} & \cH(M(F_v^+), U_{M, v}) \otimes_\bbZ \bbR \ar[d]^{\cN_M} \\
& \cH(T(F_v^+), U_{T, v}) \otimes_\bbZ \bbR, } \]
where the maps $\cN_G$, $\cN_M$ are the Satake isomorphisms defined in proof of Proposition-Definition \ref{prop_calculation_of_G_Hecke_algebras}. This diagram commutes up to multiplication by the ratio of the modulus characters for $G$ and $M$, by the transitivity of the formation of constant terms. More precisely, let $B_M = B \cap M$. Then for any $f \in \cH(G(F_v^+), U_v)$, $t \in T(F_v^+)$, we have the formula
\[ (\cN_M \cS_v f)(t) = \delta_{B_M}(t)^{1/2} \delta_B(t)^{-1/2} (\cN_G f)(t). \]
A calculation now gives the claimed formulae for the Satake transform. To finish the proof of the proposition, we must show the rationality of the element $T_{M, w, i}$ in each case. This step is essentially the same as in Proposition-Definition \ref{prop_calculation_of_G_Hecke_algebras}, so we omit the details.
\end{proof}
With the above definitions, if $w$ is a finite place of $F$ unramified over the place $v$ of $F^+$ then we define a polynomial in $\cH(G(F_v^+), U_v)[q_v^{-1}][X]$
\begin{equation}\label{eqn_G_unramified_hecke_polynomial} P_{G, w}(X) = X^{2n} - T_{G, w, 1} X^{2n - 1} + \dots + (-1)^j q_w^{j(j-1)/2} T_{G, w, j} X^{2n - j} + \dots + q_w^{n(2n-1)} T_{G, w, 2n} 
\end{equation}
and polynomials in $\cH(M(F_v^+), U_v)[q_v^{-1}][X]$
\begin{equation}\label{eqn_M_unramified_Hecke_polynomial} P_{M, w}(X) = X^n - T_{M, w, 1} X^{n-1} + \dots + (-1)^j q_w^{j(j-1)/2} T_{M, w, j} X^{n-j} + \dots + q_w^{n(n-1)/2} T_{M, w, n} 
\end{equation}
and
\begin{equation}\label{eqn_M_conjugate_unramified_Hecke_polynomial}  P_{M, w}^\vee(X) = (-1)^n (q_w^{n(n-1)/2} T_{M, w, n})^{-1} X^n P_w(X^{-1}). 
\end{equation}
Then the relation $\cS_v P_{G, w}(X) = P_{M, w}(X) q_w^{n(2n-1)} P_{M, w^c}^\vee(q_w^{1-2n} X)$ holds inside $\cH(M(F_v^+), U_v)[q_v^{-1}][X]$.

Now let $p$ be a prime, and let $S$ be a finite set of finite places of $F^+$, containing the $p$-adic places and the places which are ramified in $F$. We assume that the primes of $F^+$ above $p$ are unramfied in $F$; this implies that the group $\underline{G}_{\cO_{F_v^+}}$ is reductive for each place $v | p$ of $F^+$, so we can use the construction in \S \ref{sec_representations_of_u} to describe families of local systems on the spaces $X_G^U$. Let $E$ be a finite extension of $\bbQ_p$ which contains the image of every embedding $F \hookrightarrow \overline{\bbQ}_p$, let $\cO$ denote its ring of integers, $\pi$ a choice of uniformizer, and $k$ its residue field. We now describe a parameterization of certain $p$-adic local systems on the spaces $X_G^U$ and $X_M^{U_M}$. For convenience, we let $I_p$ denote the set of embeddings $\tau : F^+ \hookrightarrow E$, and choose for each $\tau \in I_p$ an embedding $\widetilde{\tau} : F \hookrightarrow E$ extending it. We let $\widetilde{I}_p$ denote the set of such embeddings. Let $\underline{T}_n \subset \Res^{\cO_F}_{\cO_{F^+}} \GL_n = \underline{M}$ denote the standard diagonal maximal torus. The fixed Levi embedding $\underline{M} \hookrightarrow \underline{G}$ induces an isomorphism $\underline{T}_n \cong \underline{T}$, and we will use this isomorphism to relate the parameterization of local systems on $X_G^U$ and $X_M^{U_M}$.

Fix a place $v \in S_p$, and let $\tau \in I_p$ be an embedding inducing $v$. Then the choice of $\widetilde{\tau}$ determines canonical isomorphisms $M \otimes_{F^+, \tau} E \cong \GL_n \times \GL_n$ and $T_n \otimes_{F^+, \tau} E \cong \GL_1^n \times \GL_1^n$, hence $X^\ast(T_{n, E, \tau}) \cong \bbZ^n \times \bbZ^n$. An element $(\lambda_{\widetilde{\tau}}, \lambda_{\widetilde{\tau} c}) \in \bbZ^n \times \bbZ^n$ lies in the dominant subset $X^\ast(T_{n, E, \tau})^+$ if and only if it satisfies the conditions
\[ \lambda_{\widetilde{\tau}, 1} \geq \lambda_{\widetilde{\tau}, 2} \geq \dots \geq \lambda_{\widetilde{\tau}, n} \]
\[ \lambda_{\widetilde{\tau} c, 1}  \geq \lambda_{\widetilde{\tau}c, 2} \geq \dots \geq \lambda_{\widetilde{\tau}c, n}, \]
i.e.\ if and only if it lies in the subset $\bbZ^n_+ \times \bbZ^n_+$, where we define
\[ \bbZ^n_+ = \{ (x_1, \dots, x_n) \in \bbZ^n \mid x_1 \geq x_2 \geq \dots \geq x_n \}. \]
(We will also use the notation $\bbZ^n_{++} \subset \bbZ^n_+$ to refer to the set of tuples for which these inequalities are strict.) In \S \ref{sec_representations_of_u}, we associated to a dominant pair $(\lambda_{\widetilde{\tau}}, \lambda_{\widetilde{\tau} c})$ an $\cO[\underline{M}(\cO_{F^+_v})]$-module $A(\underline{M}; \lambda_{\widetilde{\tau}}, \lambda_{\widetilde{\tau} c})$, finite free as $\cO$-module. Given a tuple
\[ \boldsymbol{\lambda} = (\lambda_{\widetilde{\tau}}) \in (\bbZ^n_+)^{\Hom(F, E)} = (\bbZ_+^n\times \bbZ_+^n)^{\Hom(F^+, E)}, \]
we define $A(\underline{M}, \boldsymbol{\lambda}) = \otimes_{\tau \in I_p} A(\underline{M}; \lambda_{\widetilde{\tau}}, \lambda_{\widetilde{\tau c}})$, the tensor product being taken over $\cO$. Then $A(\underline{M}; \boldsymbol{\lambda})$ is an $\cO[\underline{M}(\cO_{F^+} \otimes_\bbZ \bbZ_p)]$-module, finite free as $\cO$-module.

Now choose again a place $v \in S_p$, and let $\tau \in I_p$ be an embedding inducing $v$. Then the choice of $\widetilde{\tau}$ induces canonical isomorphisms $G \otimes_{F^+, \tau} E \cong \GL_{2n}$ and $T \otimes_{F^+, \tau} E \cong \GL_1^{2n}$, hence $X^\ast(T_{E, \tau}) \cong \bbZ^{2n}$. An element $a_\tau \in X^\ast(T_{E, \tau})$ lies in the dominant subset $X^\ast(T_{E, \tau})^+$ if and only if it satisfies the conditions
\[ a_{\tau, 1} \geq a_{\tau, 2} \geq \dots \geq a_{\tau, {2n}}. \]
Under the isomorphism $X^\ast(T_{E, \tau}) \cong X^\ast(T_{n, E, \tau})$, we have
\begin{equation}\label{eqn_dominant_weight_dictionary} (\lambda_{\widetilde{\tau}, 1}, \dots, \lambda_{\widetilde{\tau}, n}, \lambda_{\widetilde{\tau} c, 1}, \dots, \lambda_{\widetilde{\tau} c, n}) \leftrightarrow a_\tau = (-\lambda_{\widetilde{\tau}c, n}, \dots, -\lambda_{\widetilde{\tau}c, 1}, \lambda_{\widetilde{\tau}, 1}, \dots, \lambda_{\widetilde{\tau}, n}). 
\end{equation}
In particular, the subset $X^\ast(T_{E, \tau})^+ \subset X^\ast(T_{n, E, \tau})^+$ is described by the single extra condition $- \lambda_{\widetilde{\tau} c, 1} \geq \lambda_{\widetilde{\tau}, 1}$. We have associated to each $a_\tau \in X^\ast(T_{E})^+$ an $\cO[\underline{G}(\cO_{F_v^+})]$-module $A(\underline{G}; a_\tau)$, finite free as $\cO$-module. Given a tuple
\[ \mathbf{a} = (a_\tau) \in (\bbZ^{2n}_+)^{\Hom(F^+, E)}, \]
we define $A(\underline{G}; \mathbf{a}) = \otimes_{\tau \in I_p} A(\underline{G}; a_\tau)$, the tensor product being taken over $\cO$. Then $A(\underline{G}; a)$ is an $\cO[\underline{G}(\cO_{F^+} \otimes_\bbZ \bbZ_p)]$-module, finite free as $\cO$-module, and we have the following lemma, which follows from Corollary \ref{cor_local_tensor_product_of_weyl_modules}:
\begin{lemma}\label{lem_global_tensor_product_of_weyl_modules}
Fix an element $\mathbf{a} \in (\bbZ^{2n}_+)^{\Hom(F^+, E)}$ corresponding to $\boldsymbol{\lambda} \in (\bbZ^n_+)^{\Hom(F, E)}$ under \(\ref{eqn_dominant_weight_dictionary}\). Let $U \in \cJ_G$ be such that $U_p = U_{M, p} \ltimes U_{N, p}$ is decomposed. Then there is a direct sum decomposition
\[ \Res^{U_p}_{U_{M, p}} A(G; \mathbf{a}) = A(M; \boldsymbol{\lambda}) \oplus K \]
of $\cO[U_{M, p}]$-modules, with $A(M; \boldsymbol{\lambda}) \subset A(G; \mathbf{a})^{U_{N, p}}$.
\end{lemma}
Let $\boldsymbol{\lambda} \in (\bbZ_+^n)^{\Hom(F, E)}$. Although the weight $\boldsymbol{\lambda}$ may not be dominant for $G$, it becomes so after twisting. More precisely, let $\mathbf{1} \in (\bbZ^n_+)^{\Hom(F, E)}$ be the element with all entries equal to 1; it is the highest weight of the determinant of the standard representation of $\Res^F_\bbQ \GL_{n, F}$. For any $w \in \bbZ$, the weight $\boldsymbol{\lambda} + w \cdot \mathbf{1}$ is associated to the tensor product of $A(M; \boldsymbol{\lambda})$ with this determinant character, raised to the power $w$. For $w$ sufficiently negative (more precisely, for $w \leq -\sup_{\tau \in I_p}(\lambda_{\widetilde{\tau}c, 1} + \lambda_{\widetilde{\tau}, 1})/2$), we have $\boldsymbol{\lambda} + w \cdot \mathbf{1} \in (\bbZ_+^{2n})^{\Hom(F^+, E)}$. In this connection, we have the following useful lemma.
\begin{lemma}\label{lem_twisting_preserves_hecke_algebra}
For any $k \in \bbZ$ and $U \in \cJ_{M}$, there is a canonical isomorphism in $\mathbf{D}(\cO)$:
\[ R \Gamma_{X_M^U} \underline{A(M, \boldsymbol{\lambda})}_M^U \cong R \Gamma_{X_M^{U}} \underline{A(M, \boldsymbol{\lambda} + k\cdot \mathbf{1})}_M^U.  \]
This isomorphism intertwines the action of $T_{M, w, i}$ on the left hand side and $q_w^{ki} T_{M, w, i}$ on the right hand side. Consequently, there is an isomorphism
\[ \bbT^S_M(R \Gamma_{X_M^U} \underline{A(M; \boldsymbol{\lambda})}_M^U ) \cong \bbT^S_M(R \Gamma_{X_M^{U}} \underline{A(M; \boldsymbol{\lambda} + k \cdot \mathbf{1})}_M^U )  \]
which sends $T_{M, w, i}$ to $\epsilon(\Frob_w)^{-k} T_{M, w, i} = q_w^{ki} T_{M, w, i}$.
\end{lemma}
\begin{proof}
The $k^\text{th}$ power of the cyclotomic character defines a class $[\epsilon^k]$ in $H^0(X_M^U, \underline{A(M; k \cdot \mathbf{1})}_M^U)$. More precisely, we can interpret $H^0(X_M^U, \underline{A(M; k \cdot \mathbf{1})}_M^U)$ as the set of $U$-equivariant sections of the map
\[ M(F) \backslash M^\infty \times A(M, k \cdot \mathbf{1}) \to M(F) \backslash M^\infty. \]
Any character of the form $\chi = \epsilon^k \psi$, where $\psi$ satisfies $\psi( \Art_F (\det U)) = 1$, determines such a section by the formula
\[ m^\infty \mapsto (m^\infty, \chi( \Art_F (\det m^\infty) ) ). \]
 There is an isomorphism $A(M; \boldsymbol{\lambda}) \otimes_\cO A(M; k \cdot \mathbf{1}) \cong A(M; \boldsymbol{\lambda} + k \cdot \mathbf{1})$, and the first isomorphism of the lemma is defined by cup product with the class $[\epsilon^k]$. The remainder of the lemma now follows from the relation $T_{M, w, i}(c \cup [\chi]) = \chi(\Frob_w)^i (T_{M, w, i}(c) \cup [\chi])$.
\end{proof}
\subsection{Application to Galois representations, I}\label{sec_application_to_galois_I}

In this section, we prove our first main theorem about the existence of Galois representations (Theorem \ref{thm_intro_theorem_no_assumption_on_p}). It will be convenient to introduce some notation. Let $S$ be a finite set of finite places of $F^+$, containing the $p$-adic places, and let $U_S \subset \prod_{v \in S} \underline{G}(\cO_{F_v^+})$ be an open compact subgroup. If $N \geq 1$, $\mathbf{a}\in (\bbZ_+^{2n})^{\Hom(F^+, E)}$ and $U \in \cJ_{G, U_S}$, then we define ideals
\[ J_{U, \mathbf{a}, N} = \ker\left( \bbT_G^S \to \End_{\mathbf{D}(\cO)}(R \Gamma_{X_G^U} \underline{A(G; \mathbf{a})}_G^U\otimes_\cO \cO/(\pi^N)) \right) \]
and
\[ J_{U,  \mathbf{a}, c, N} = \ker\left( \bbT_G^S \to \End_{\mathbf{D}(\cO)}(R \Gamma_{X_G^U, c} \underline{A(G; \mathbf{a})}_G^U\otimes_\cO \cO/(\pi^N)) \right).\]
Thus we have, for example, 
\[ \bbT_G^S / J_{U, \mathbf{a}, N} = \bbT_G^S(R \Gamma_{X_G^U} \underline{A(G; \mathbf{a})}_G^U \otimes_\cO \cO/(\pi^N)). \]
\begin{definition}\label{defsuffsmall}
A compact open subgroup $U \subset G(\bbA^\infty_{\Fp})$ is \emph{small} if for some rational prime $q \ne p$, $U_q = \prod_{v|q}U_v$ is contained in the principal congruence subgroup $U(q)$ given by the kernel of the map \[\prod_{v|q}G(\cO_{F_v^+}) \rightarrow \prod_{v|q}G(\cO_{F_v^+}/\mathbf{q}\cO_{F_v^+}),\] where $\mathbf{q} = q$ if $q$ is odd and $\mathbf{q} = 4$ otherwise.

Similarly, we say that a compact open subgroup $U \subset \GL_n(\bbA^\infty_F)$ is \emph{small} if for some rational prime $q \ne p$, $U_q = \prod_{v|q}U_v$ is contained in the principal congruence subgroup $U(q)$ given by the kernel of the map \[\prod_{v|q}\GL_n(\cO_{F_v}) \rightarrow \prod_{v|q}G(\cO_{F_v}/\mathbf{q}\cO_{F_v}).\] 
\end{definition}
\begin{remark}
\begin{enumerate}\label{smallremark}
\item We introduce this condition in order to be able to apply \cite[Theorem 4.1.1]{Sch14}.
\item For a rational prime $q$, every root of unity $\zeta \in \overline{\bbZ}_q$ which is congruent to $1$ mod $\mathbf{q}$ is equal to $1$ (see \cite[0.6]{Pin90}).
\item If $U \subset G(\bbA^\infty_{\Fp})$ or $U \subset \GL_n(\bbA^\infty_F)$ is small then it is neat. This follows from the fact mentioned above, since for $v|q$ the eigenvalues of $g_v \in U_v$ under a faithful representation are congruent to $1$ mod $\mathbf{q}$.  
\item If $U \subset G(\bbA^\infty_{\Fp})$ is small then the image of $U$ under an algebraic group homomorphism $i:\Res^{\Fp}_\bbQ G \rightarrow H$ is contained in a neat compact open subgroup of $H(\bbA_\bbQ^\infty)$. We use the fact that for $g \in U$ the eigenvalues of $i(g)_q$ acting under a faithful representation of $H$ are congruent to $1$ mod $\mathbf{q}$. We take a faithful representation of $H$ and let $K$ be the compact open subgroup of $H(\bbA_\bbQ^\infty)$ obtained by taking the inverse image under this representation of the level $\mathbf{q}$ principal congruence subgroup. Then $U$ and $\cap_{g \in U}g K g^{-1}$ generate a neat compact open subgroup of $H(\bbA_\bbQ^\infty)$.
\end{enumerate}
\end{remark}

Our starting point for the proof of Theorem \ref{thm_intro_theorem_no_assumption_on_p} will be the following result:
\begin{theorem}\label{thm_existence_of_G_Galois_representations}
Fix $N \geq 1$. For each small $U \in \cJ_{G, U_S}$, there exists an ideal $J_{U, N} \subset \bbT_G^S$ satisfying the following conditions:
\begin{enumerate}
\item There exists a continuous group determinant of dimension 2n 
\[ D_{G, U} : G_{F, S} \to \bbT^S_G/J_{U,N} \]
 such that for every finite place $w\not\in S$ of $F$, $D_{G, U}(\Frob_w)$ has characteristic polynomial $P_{G, w}(X)$ (see equation (\ref{eqn_G_unramified_hecke_polynomial})).
\item For any $\mathbf{a} \in (\bbZ_+^{2n})^{\Hom(F^+, E)}$, we have $J_{U, N} \subset J_{U, \mathbf{a}, N}$ and $J_{U, N} \subset J_{U, \mathbf{a}, c, N}$.
\end{enumerate}
\end{theorem}
\begin{proof}
We can find an open normal subgroup $V_p \subset U_p$ such that the action of $V_p$ on $A(G; \mathbf{a}) / (\pi^N)$ is trivial for all $\mathbf{a} \in (\bbZ_+^{2n})^{\Hom(F^+, E)}$. Let $V = U^p V_p$. We define 
\[ J_{U, N} = \ker( \bbT_G^S \to \End_{\mathbf{D}(U/V, \cO/(\pi^N))}(R \Gamma_{X_G^V, c} \cO/(\pi^N) ). \]
The existence of $D_{G, U}$ thus follows from Theorem \ref{thmclass}, to be proved in \S \ref{sec_james} below, and \cite[Corollary 5.1.11]{Sch14}. There is a canonical isomorphism in $\mathbf{D}(\cO/(\pi^N))$:
\[R\Gamma_{X_G^U,c}\left(\underline{A(G;\mathbf{a})_G^U}\otimes_{\cO}\cO/(\pi^N)\right) = R\Gamma_{U/V} \left( R\Gamma_{X_G^V,c}\cO/(\pi^N) \otimes^\bbL_{\cO/(\pi^N)[U/V]} A(G;\mathbf{a})/(\pi^N) \right).\]
This implies the inclusion $J_{U, N} \subset J_{U, \mathbf{a}, c, N}$. The inclusion $J_{U, N} \subset J_{U, \mathbf{a}, N}$ follows by Verdier duality, in the guise of Proposition \ref{prop_hecke_action_and_verdier_duality}.
\end{proof}
We use this to prove the following result. 
\begin{proposition}\label{prop_existence_of_galois_in_general_case}
Let $\boldsymbol{\lambda} \in (\bbZ_{+}^{n})^{\Hom(F, E)}$, and let $\ffrm \subset \bbT^S_M(R \Gamma_{X_M^{U_M}} \underline{A(M; \boldsymbol{\lambda})}_M^{U_M})$ be a non-Eisenstein maximal ideal. Let $N \geq 1$ be an integer.
For every continuous character $\chi : G_{F, S} \to \cO^\times$ of finite odd order, prime to $p$, there exists an ideal 
\[ I_{U, \chi} \subset \bbT^S_M(R \Gamma_{X_M^{U_M}} \underline{A(M; \boldsymbol{\lambda})}_M^{U_M} \otimes_\cO \cO/(\pi^N))_\ffrm\] of square 0 and a continuous group determinant 
\[ D_{M, U, \chi} : G_{F, S} \to \bbT^S_M(R \Gamma_{X_M^{U_M}} \underline{A(M; \boldsymbol{\lambda})}_M^{U_M}\otimes_\cO \cO/(\pi^N))_\ffrm/I_{U, \chi} \]
of dimension $2n$ such that for every finite place $w\not\in S$ of $F$, $D_{M, U, \chi}(\Frob_w)$ has characteristic polynomial given by $\chi(\Frob_w)^{n} P_{M, w}(\chi(\Frob_w)^{-1}X)\chi(\Frob_{w^c})^{-n} q_w^{n(2n-1)} P_{M, w^c}^\vee(q_w^{1-2n} \chi(\Frob_{w^c}) X)$.
\end{proposition}
\begin{proof}
By Lemma \ref{lem_twisting_preserves_hecke_algebra}, we can and do replace $\boldsymbol{\lambda}$ by $\boldsymbol{\lambda} + w \cdot \mathbf{1}$, where
\[ w = -\sup_{\tau \in I_p} \lfloor (\lambda_{\tau, 1} + \lambda_{\tau c, 1})/2 \rfloor. \]
Then the weight $\mathbf{a} \in (\bbZ_{+}^{2n})^{\Hom(F^+, E)}$ corresponding to $\boldsymbol{\lambda}$ is dominant, so the coefficient module $A(G; \mathbf{a})$ is defined. We first treat the case $\chi = 1$. There is an exact triangle in $\mathbf{D}(\cO)$:
\[ \xymatrix@1{ R \Gamma_{X_G^U, c} \underline{A(G; \mathbf{a})}_G^{U} \otimes_\cO \cO/(\pi^N) \ar[r] & R \Gamma_{X_G^U} \underline{A(G; \mathbf{a})}_G^{U} \otimes_\cO \cO/(\pi^N)\ar[r] & R\Gamma_{\partial \overline{X}_G^U} \underline{A(G; \mathbf{a})}_G^U \otimes_\cO \cO/(\pi^N). } \]
Define 
\[ J_{U, \mathbf{a}, \partial, N} = \ker\left( \bbT_G^S \to \End_{\mathbf{D}(\cO)}( R\Gamma_{\partial \overline{X}_G^U} \underline{A(G; \mathbf{a})}_G^U\otimes_\cO \cO/(\pi^N)) \right). \]
It follows from Lemma \ref{lem_square_zero_endomorphisms} that 
\[ J_{U, N}^2 \subset J_{U, \mathbf{a},  N} \cdot J_{U, \mathbf{a}, c, N} \subset J_{U, \mathbf{a}, \partial, N}. \]
By Theorem \ref{thm_existence_of_G_Galois_representations}, we find that there exists a continuous group determinant $D_{G, U, \partial} : G_{F, S} \to \bbT_G^S/( J_{U, \mathbf{a}, \partial, N}, J_{U, N})$ such for each finite place $w \not\in S$ of $F$, $D_{G, U, \partial}(\Frob_w)$ has characteristic polynomial $P_{G, w}(X)$.

By Corollary \ref{cor_restriction_of_hecke_to_boundary} and Theorem \ref{thm_boundary_cohomology_is_eisenstein}, there is a commutative diagram
\[ \xymatrix{ \bbT_G^S \ar[r] \ar[d] & \End_{\mathbf{D}(\cO/(\pi^N))}((R \Gamma_{\partial \overline{X}_G^U} \underline{A(G; \mathbf{a})}_G^U\otimes_\cO \cO/(\pi^N))_{\cS^\ast(\ffrm)}) \ar[d] \\
\bbT_M^S \ar[r] & \End_{\mathbf{D}(\cO/(\pi^N))}( (R \Gamma_{X_M^{U_M}} \underline{A(M; \boldsymbol{\lambda})}_M^{U_M}\otimes_\cO \cO/(\pi^N))_\ffrm),} \]
and hence a canonical map
\begin{equation}\label{eqn_canonical_surjection_of_hecke} \bbT_G^S(R \Gamma_{\partial \overline{X}_G^U} \underline{A(G; \mathbf{a})}_G^U\otimes_\cO \cO/(\pi^N))_{\cS^\ast(\ffrm)} \to \bbT_M^S(R \Gamma_{X_M^{U_M}} \underline{A(M; \boldsymbol{\lambda})}_M^{U_M}\otimes_\cO \cO/(\pi^N))_\ffrm, 
\end{equation}
which is induced by the unnormalized Satake transform. The proof of the proposition is completed in this case on taking $D_{M, U, 1}$ to be the image of $D_{G, \partial, U}$ under the map (\ref{eqn_canonical_surjection_of_hecke}) and $I_{\chi, U}$ to be the image of the ideal $J_{U, N}$.

We now treat the case of an arbitrary character $\chi$. We can find a normal subgroup $V_M \in \cJ_{M, U_{M, S}}$ of $U_M$ such that the index $[U_M : V_M]$ is prime to $p$ and the character $\chi \circ \Art_F : \bbA_F^\infty \to \cO^\times$ is trivial on $\det V_M$. We can find $V \in \cJ_{G, U_S}$ such that $V \cap M^\infty = V_M$ (so the notation is consistent). In this case we can describe a class $[\chi] \in H^0(X_{M}^{V_M}, \cO)$ as in Lemma \ref{lem_twisting_preserves_hecke_algebra} on which $\GL_n(\bbA_F^\infty)$ acts by the character $\chi \circ \Art_F \circ \det$. Pullback and cup product with the class $[\chi]$ then defines a map 
\[ F_\chi : R \Gamma_{X_M^{U_M}} \underline{A(M; \boldsymbol{\lambda})}_M^{U_M}\otimes_\cO \cO/(\pi^N) \to  R \Gamma_{X_M^{V_M}} \underline{A(M; \boldsymbol{\lambda})}_M^{V_M}\otimes_\cO \cO/(\pi^N) \]
which is an isomorphism onto a direct summand $A^\bullet$ in $\mathbf{D}(\cO/(\pi^N))$ which is invariant under the action of $\bbT_M^S$. We obtain a commutative diagram
\[ \xymatrix{ \bbT_M^S \ar[d]_-{T_{M, w, i} \mapsto \chi(\Frob_w)^{-i} T_{M, w, i}}^\cong \ar[r] & \End_{\mathbf{D}(\cO/(\pi^N))}(R \Gamma_{X_M^{U_M}}  \underline{A(M; \boldsymbol{\lambda})}_M^{U_M}\otimes_\cO \cO/(\pi^N)) \ar[d]_\cong^-{F_\chi (\cdot) F_\chi^{-1}} \\
 \bbT_M^S \ar[r] & \End_{\mathbf{D}(\cO/(\pi^N))}(A^\bullet) \\
 \bbT_M^S \ar[u]_-= \ar[r] & \ar[u] \End_{\mathbf{D}(\cO/(\pi^N))}(R \Gamma_{X_M^{V_M}}\underline{A(M; \boldsymbol{\lambda})}_M^{V_M}\otimes_\cO \cO/(\pi^N)).} \]
This diagram gives a morphism of Hecke algebras
\[ f_\chi : \bbT_M^S(R  \Gamma_{X_M^{V_M}}\underline{A(M; \boldsymbol{\lambda})}_M^{V_M}\otimes_\cO \cO/(\pi^N)) \to \bbT_M^S(A^\bullet) \to \bbT_M^S(R \Gamma_{X_M^{U_M}}  \underline{A(M; \boldsymbol{\lambda})}_M^{U_M}\otimes_\cO \cO/(\pi^N)) \]
which sends the operator $T_{M, w, i}$ to $\chi(\Frob_w)^{-i} T_{M, w, i}$. The proof is completed in this case on taking the ideal $I_{U, \chi}$ to be the image under $f_\chi$ of the ideal $\cS(J_{V, N})$ in $\bbT_M^S(R \Gamma_{X_M^{U_M}}  \underline{A(M; \boldsymbol{\lambda})}_M^{U_M}\otimes_\cO \cO/(\pi^N))$, and $D_{M, U, \chi}$ to be the image under $f_\chi$ of $D_{M, V, 1}$.
\end{proof}
Using Proposition \ref{prop_existence_of_galois_in_general_case}, we can prove our first main theorem.
\begin{theorem}\label{thm_existence_of_galois_in_general_case}
Let $\boldsymbol{\lambda} \in (\bbZ_{+}^n)^{\Hom(F, E)}$, and let $U \in \cJ_{G, U_S}$ be small. Let $\ffrm \subset \bbT^S_M(R \Gamma_{X_M^{U_M}} \underline{A(M; \boldsymbol{\lambda})}_M^{U_M})$ be a non-Eisenstein maximal ideal. Then there exists an ideal $I \subset \bbT^S_M(R \Gamma_{X_M^{U_M}} \underline{A(M; \boldsymbol{\lambda})}_M^{U_M})_\ffrm$ satisfying $I^4 = 0$ and a continuous representation
\[ \rho_\ffrm : G_{F, S} \to \GL_n(\bbT^S_M(R \Gamma_{X_M^{U_M}} \underline{A(M; \boldsymbol{\lambda})}_M^{U_M})_\ffrm/I) \]
such that for each place $w\not\in S$ of $F$, we have the equality
\[ \det( X \cdot 1_n - \rho_\ffrm(\Frob_w)) = P_{M, w}(X) \]
inside $(\bbT^S_M(R \Gamma_{X_M^{U_M}} \underline{A(M; \boldsymbol{\lambda})}_M^{U_M})_\ffrm/I)[X]$.
\end{theorem}
\begin{proof}
Given Proposition \ref{prop_existence_of_galois_in_general_case}, exactly the same `separation of parts' argument as in \cite[5.3]{Sch14} implies that for each $N \geq 1$, there is an ideal 
\[ I_N \subset \bbT^S_M(R \Gamma_{X_M^{U_M}} \underline{A(M; \boldsymbol{\lambda})}_M^{U_M} \otimes_\cO \cO/(\pi^N))_\ffrm = \bbT_N, \] 
say, satisfying $I_N^4 = 0$, and a continuous group determinant $D_{M, U, N} : G_{F, S} \to \bbT_N/I_N$
of dimension $n$ and with the expected characteristic polynomial.

Let $\bbT_\infty = \bbT^S_M(R \Gamma_{X_M^{U_M}} \underline{A(M; \boldsymbol{\lambda})}_M^{U_M})_\ffrm$. By Lemma \ref{lem_inverse_limit_of_mod_p_hecke_is_hecke}, the natural map $\bbT_\infty \to \plim_N \bbT_N$ is an isomorphism. Let $I_\infty$ denote the kernel of the map $\bbT_\infty \to \prod_N \bbT_N / I_N$. Then we have $I_\infty^4 = 0$, and by \cite[Example 2.32]{Che14}, the group determinants $D_{M, U, N}$ glue into a group determinant $D_{M, \infty} : G_{F, S} \to \bbT_\infty / I_\infty$.

In order to obtain a true representation $\rho_\ffrm$ at the end, instead of just the group determinant $D_{M, \infty}$, we recall that the deforming the determinant $\overline{D}_\ffrm = \det(X \cdot 1_n - \overline{\rho}_\ffrm)$ is equivalent to deforming $\overline{\rho}_\ffrm$, because of the assumption that the residual representation $\overline{\rho}_\ffrm$ is absolutely irreducible (see \cite[Theorem 2.22]{Che14}).
\end{proof}
To deduce Theorem \ref{thm_intro_theorem_no_assumption_on_p} of the introduction from Theorem \ref{thm_existence_of_galois_in_general_case}, we need only observe that for any small $V \in \cJ_{\GL_n, \prod_{v \not\in S} \GL_n(\cO_{F_v})}$, we can find a small $U \in \cJ_{G, \prod_{v \not\in S} \underline{G}(\cO_{F_v^+})}$ such that $U_M = V$.
\subsection{Application to Galois representations, II}\label{sec_application_to_galois_II}
In this section, we prove our second main theorem (Theorem \ref{thm_intro_theorem_p_regular}). We will make use of the following result of Lan--Suh (\cite[Theorem 10.1]{Lan13}). Let $U_p = \underline{G}(\cO_{F^+} \otimes_\bbZ \bbZ_p)$.
\begin{theorem}\label{thm_lan_suh_vanishing}
Suppose that $p$ is unramified in $F$, and choose $\mathbf{a} \in (\bbZ_{+}^{2n})^{\Hom(F^+, E)}$, $U \in \cJ_{G, U_p}$ small. Suppose that $\mathbf{a}$ satisfies the following further conditions:
\begin{equation}
\text{For each }\tau \in \Hom(F^+, E), a_{\tau, 1} > a_{\tau, 2} > \dots > a_{\tau, 2n};
\end{equation}
\begin{equation}\label{eqn_condition_for_lan_suh_vanishing}
\text{ and }dn(n+1) + \sum_{\tau} \sum_i (a_{\tau, i} - 2 \lfloor a_{\tau, 2n} / 2 \rfloor ) < p.
\end{equation}
Then 
\[ H^i(X_G^U, \underline{A(G; \boldsymbol{\lambda})}_G^U) = H^i(X_G^U, \underline{A(G; \boldsymbol{\lambda})}_G^U \otimes_\cO k) = 0 \]
for each $0 \leq i \leq D - 1 $.
\end{theorem}
\begin{proof}
We show how to deduce this from \cite{Lan13}. The necessary conditions on $\mathbf{a} \in (\bbZ_{++}^{2n})^{\Hom(F^+, E)}$ appearing in \cite[Theorem 10.1]{Lan13} are $| \mathbf{a} |_{\text{re}, +} < p$, $|\mathbf{a}|'_\text{comp} \leq p - 2$, and for all $\tau \in I_p$, $a_{\tau, 1} - a_{\tau, 2n} < p$. According to \cite[Definition 9.7]{Lan13}, we have
\[ |\mathbf{a}|'_\text{comp} = 1 + D + |\mathbf{a}|_L = 1 + dn^2 + |\mathbf{a}|_L, \]
where (\cite[Definition 3.2]{Lan12}) 
\[ |\mathbf{a}|_L = \sum_{\tau \in I_p} \sum_{i=1}^{2n} \left( a_{\tau, i} - 2 \lfloor a_{\tau, 2n} / 2 \rfloor \right). \]
 According to \cite[(7.22)]{Lan12} and \cite[Definition 3.9]{Lan12}, we have
\[ |\mathbf{a}|_{\text{re}, +} = |\mathbf{a}|_{\text{re}} + dn = D + |\mathbf{a}|_L + dn = dn(n+1) + |\mathbf{a}|_L. \]
After rearranging, the condition $|\mathbf{a}|_{\text{re},+} < p$ becomes (\ref{eqn_condition_for_lan_suh_vanishing}) above, and it is easy to see that this implies the other two conditions.
\end{proof}
\begin{corollary}\label{cor_quasi_isomorphism_of_boundary_and_compactly_supported_cohomology}
Let $N \geq 1$. With assumptions as in Theorem \ref{thm_lan_suh_vanishing}, the map $R \Gamma_{\partial \overline{X}_G^U} \underline{A(G;\mathbf{a})}_G^U \otimes_\cO \cO/(\pi^N) \to R \Gamma_{X_G^U, c} \underline{A(G;\mathbf{a})}_G^U\otimes_\cO \cO/(\pi^N)[-1]$ induces an isomorphism
\[ \tau_{\leq D-2} R \Gamma_{\partial \overline{X}_G^U} \underline{A(G;\mathbf{a})}_G^U \otimes_\cO \cO/(\pi^N)\cong \tau_{\leq D-1} (R \Gamma_{X_G^U, c} \underline{A(G; \mathbf{a})}_G^U\otimes_\cO \cO/(\pi^N))[-1] \]
in $\mathbf{D}(\cO/(\pi^N))$. 
\end{corollary}
In the situation of Corollary \ref{cor_quasi_isomorphism_of_boundary_and_compactly_supported_cohomology}, we can prove the following refinement of Proposition \ref{prop_existence_of_galois_in_general_case}.
\begin{proposition}\label{prop_existence_of_galois_in_regular_case}
Let $\boldsymbol{\lambda} \in (\bbZ_{+}^{n})^{\Hom(F, E)}$, let $U \in \cJ_{G, U_p}$ be small, and let $\ffrm \subset \bbT^S_M(R \Gamma_{X_M^{U_M}} \underline{A(M; \boldsymbol{\lambda})}_M^{U_M})$ be a non-Eisenstein maximal ideal. Suppose that $p$ is unramified in $F$ and that $\boldsymbol{\lambda}$ satisfies the following conditions:
\begin{equation}
\text{For all }\tau \in \Hom(F, E), \lambda_{\tau, 1} > \lambda_{\tau, 2} > \dots > \lambda_{\tau, n}; 
\end{equation}
\begin{equation}\label{eqn_lan_suh_condition_for_levi_subgroup} \text{and } dn(n + 6 + \sup_\tau (\lambda_{\wt, 1} + \lambda_{\wt c, 1})) + \sum_{\tau \in I_p} \sum_{i=1}^n \left( \lambda_{\wt , i} - \lambda_{\wt c, i} - 2 \lambda_{\wt, n} \right) < p. 
\end{equation}
Let $N \geq 1$. For every continuous character $\chi : G_{F, S} \to \cO^\times$ of finite odd order, prime to $p$, there exists a continuous group determinant 
\[ D_{M, U, \chi} : G_{F, S} \to \bbT^S_M(R \Gamma_{X_M^{U_M}} \underline{A(M; \boldsymbol{\lambda})}_M^{U_M}\otimes_\cO \cO/(\pi^N))_\ffrm\]
 of dimension $2n$ such that for every finite place $w\not\in S$ of $F$, $D_{M, U, \chi}(\Frob_w)$ has characteristic polynomial given by $\chi(\Frob_w)^{n} P_{M, w}(\chi(\Frob_w)^{-1}X)\chi(\Frob_{w^c})^{-n} q_w^{n(2n-1)} P_{M, w^c}^\vee(q_w^{1-2n} \chi(\Frob_{w^c}) X)$.
\end{proposition}
\begin{proof}
By Lemma \ref{lem_twisting_preserves_hecke_algebra}, we can and do replace $\boldsymbol{\lambda}$ by $\boldsymbol{\lambda} + w \cdot \mathbf{1}$, where
\[ w = -\sup_{\tau \in I_p} \lfloor (\lambda_{\wt, 1} + \lambda_{\wt c, 1})/2 \rfloor - 1. \]
Then the inequality (\ref{eqn_lan_suh_condition_for_levi_subgroup}) implies that the weight $\mathbf{a} \in (\bbZ_{++}^{2n})^{\Hom(F^+, E)}$ associated to $\boldsymbol{\lambda}$ satisfies the conditions of Corollary \ref{cor_quasi_isomorphism_of_boundary_and_compactly_supported_cohomology}. In the remainder of the proof we just treat the case $\chi = 1$, since the modifications in the case $\chi \neq 1$ are exactly the same as in the proof of Proposition \ref{prop_existence_of_galois_in_general_case}.

By Corollary \ref{cor_restriction_of_hecke_to_boundary}, Theorem \ref{thm_boundary_cohomology_is_eisenstein}, and Corollary \ref{cor_quasi_isomorphism_of_boundary_and_compactly_supported_cohomology}, there is a commutative diagram
\[ \xymatrix{  & \End_{\mathbf{D}(\cO/(\pi^N))}( (R \Gamma_{X_G^U, c} \underline{A(G; \boldsymbol{\lambda})}_G^U \otimes_\cO \cO/(\pi^N) )_{\cS^\ast(\ffrm)} )\ar[d]\\
\bbT_G^S \ar[r] \ar[ru] \ar[d] & \End_{\mathbf{D}(\cO/(\pi^N))}( (\tau_{\leq D-1} R \Gamma_{X_G^U, c} \underline{A(G; \boldsymbol{\lambda})}_G^U \otimes_\cO \cO/(\pi^N))_{\cS^\ast(\ffrm)}) \ar[d] \\
\bbT_M^S \ar[r] & \End_{\mathbf{D}(\cO/(\pi^N))}( \tau_{\leq D-2}(R \Gamma_{X_M^{U_M}} \underline{A(M; \boldsymbol{\lambda})}_M^{U_M} \otimes_\cO \cO/(\pi^N))_{\ffrm}). } \]
We have $\dim_\bbR X_M = D - 1$, and the top degree cohomology of $X_M^{U_M}$ is 0 (as $X_M^{U_M}$ is non-compact). It follows that the natural map
\[ \tau_{\leq D-2} R \Gamma_{X_M^{U_M}} \underline{A(M; \boldsymbol{\lambda})}_M^{U_M}  \otimes_\cO \cO/(\pi^N)\to  R \Gamma_{X_M^{U_M}} \underline{A(M; \boldsymbol{\lambda})}_M^{U_M} \otimes_\cO \cO/(\pi^N) \]
is a quasi-isomorphism. We find that the unnormalized Satake transform induces a map
\[ \bbT_G^S(R \Gamma_{X_G^U, c} \underline{A(G; \boldsymbol{\lambda})}_G^U \otimes_\cO \cO/(\pi^N))_{\cS^\ast(\ffrm)} \to \bbT_M^S( R \Gamma_{X_M^{U_M}} \underline{A(M; \boldsymbol{\lambda})}_M^{U_M} \otimes_\cO \cO/(\pi^N))_\ffrm. \]
The proposition now follows from the existence of this map and Theorem \ref{thm_existence_of_G_Galois_representations}.
\end{proof}
We finally obtain our second main theorem.
\begin{theorem}\label{thm_existence_of_galois_in_regular_case}
Suppose that $p$ is unramified in $F$ and let $U_p = \underline{G}(\cO_{F^+} \otimes_\bbZ \bbZ_p)$. Let $\boldsymbol{\lambda} \in (\bbZ^n_{+})^{\Hom(F, E)}$, and let $U \in \cJ_{G, U_p}$ be small. Let $\ffrm \subset \bbT^S_M(R \Gamma_{X_M^{U_M}} \underline{A(M; \boldsymbol{\lambda})}_M^{U_M})$ be a non-Eisenstein maximal ideal. Suppose that there exists a set $\widetilde{I}_p = \{ \widetilde{\tau} \mid \tau \in \Hom(F^+, E) \}$ such that $\boldsymbol{\lambda}$ satisfies the following conditions:
\begin{equation}
\text{For all }\tau : F \hookrightarrow E, \lambda_{\tau, 1} > \lambda_{\tau, 2} > \dots > \lambda_{\tau, n};
\end{equation}
\begin{equation}
\text{ and }[F^+ : \bbQ] n(n + 6 + \sup_{\tau \in \Hom(F^+, E)} (\lambda_{\widetilde{\tau}, 1} + \lambda_{\widetilde{\tau} c, 1})) + \sum_{\tau \in \Hom(F^+, E)} \sum_{i=1}^n \left( \lambda_{\widetilde{\tau}, i} - \lambda_{\widetilde{\tau} c, i} - 2 \lambda_{\widetilde{\tau}, n} \right) < p. 
\end{equation}
Then there exists a continuous representation
\[ \rho_\ffrm : G_{F, S} \to \GL_n(\bbT^S_M(R \Gamma_{X_M^{U_M}} \underline{A(M; \boldsymbol{\lambda})}_M^{U_M})_\ffrm) \]
such that for each place $w\not\in S$ of $F$, we have the equality
\[ \det( X \cdot 1_n - \rho_\ffrm(\Frob_w)) =  P_{M, w}(X) \]
inside $\bbT^S_M(R \Gamma_{X_M^{U_M}} \underline{A(M; \boldsymbol{\lambda})}_M^{U_M})_\ffrm[X]$, where $P_{M, w}(X)$ is as defined by (\ref{eqn_M_unramified_Hecke_polynomial}).
\end{theorem}
\begin{proof}
The deduction of Theorem \ref{thm_existence_of_galois_in_regular_case} from Proposition \ref{prop_existence_of_galois_in_regular_case} is essentially the same as the deduction of Theorem \ref{thm_existence_of_galois_in_general_case} from Proposition \ref{prop_existence_of_galois_in_general_case}, although slightly easier (since there is no longer any nilpotent ideal to worry about). We therefore omit the details. 
\end{proof}
\subsection{The proof of Theorem \ref{thm_existence_of_G_Galois_representations}, by $p$-adic interpolation}\label{sec_james}

The rest of this paper is devoted to the proof of Theorem \ref{thmclass} below, which was used in the proof of Theorem \ref{thm_existence_of_G_Galois_representations}. Let $U_p = \prod_{v\in S_p}U_v$ where the $U_v$ are compact open subgroups of $G(F^+_v)$. Let $V_p$ be a normal compact open subgroup of $U_p$. Fix $U = U_p U^p \in \mathcal{J}_{G,U_p}$ small and set $V = V_p U^p$. Note that $V \in \mathcal{J}_{G,V_p}$ is also small. Fix $N \ge 1$, and set $\Lambda = \cO/(\pi^N)$. Denote by $\underline{\Lambda}_{U/V}$ the $U/V$-equivariant sheaf on $X_G^V$ given by pulling back the constant sheaf $\underline{\Lambda}$ on $X_G^U$. There is a canonical homomorphism (arising from the diagram (\ref{eqn_completed_coh_hecke_functor_diagram}))\[\bbT^S=\bbT^S_G \rightarrow \End_{\mathbf{D}(\Lambda[U/V])}(R\Gamma_{X_G^V,c}\underline{\Lambda}_{U/V}).\] 
We write $\bbT^S(R\Gamma_{X_G^V,c}\underline{\Lambda}_{U/V})$ for the image of this map. Note that this image is a finite (commutative) $\Lambda$-algebra. We are going to show that $\bbT^S(R\Gamma_{X_G^V,c}\underline{\Lambda}_{U/V})$ is a quotient of a Hecke algebra acting faithfully on spaces of classical cuspidal automorphic forms for $G$ of (varying) regular weight. First we define this `classical' Hecke algebra, as in the statement of \cite[Theorem 4.3.1]{Sch14}. 

We denote by $X_G^{U,\mathrm{alg}}/\overline{\bbQ}$ the algebraic model over $\overline{\bbQ}$ for $X_G^U$ provided by \cite[Theorem 1]{Fal84}, and denote by $X_G^{U,*}/\overline{\bbQ}$ its minimal compactification. We write $X_G^{U,*,\mathrm{ad}}/C$ (and the same thing without $*$) for the adic space obtained by base changing the appropriate scheme to $C$ (a fixed complete and algebraically closed extension of $\overline{\bbQ}_p$, with a fixed embedding $\overline{\bbQ}\subset C$) and then taking the associated adic space over $\mathrm{Spa}(C,\cO_C)$. We can define this `adification' by first taking the associated rigid analytic variety over $C$ \cite[9.3.4]{BGR} and then apply the functor $r$ of \cite[Proposition 4.3]{Hub94}. We write $\mathcal{I}$ for the subsheaf of $\cO_{X_G^{U,*,\mathrm{ad}}}$ corresponding to the boundary $X_G^{U,*,\mathrm{ad}} \backslash X_G^{U,\mathrm{ad}}$. 

We fix an embedding $(\Res_\bbQ^{F^+} G, X_G) \rightarrow (\Sp_{2g} , D_{\Sp_{2g}})$, and following Scholze \cite[before Theorem 4.1.1]{Sch14} write $X_G^{\overline{U},*}$ for the scheme theoretic image of $X_G^{U,*}$ in $X_{\Sp_{2g}}^{U', \ast}$, where $U'$ is any sufficiently small subgroup of $\Sp_{2g}(\bbA_\bbQ^\infty)$ such that $U' \cap G(\bbA_{F^+}^\infty) = U$. Since $U$ is small, $U'$ can be assumed to be neat (see \ref{smallremark}). We write $\omega_U$ for the ample line bundle on $X_G^{\overline{U},*}$ obtained by pullback from the natural ample line bundle $X_{\Sp_{2g}}^{U', \ast}$, and we also write $\omega_U$ for the (ample) pullback to $X_G^{U,*}$ and the line bundle on the associated adic space. Since $X_G^{U,*}$ is normal and the boundary has codimension $\ge 2$, $\omega_U$ is the unique line bundle on $X_G^{U,*}$ extending $\omega_U|_{X_G^{U,\mathrm{alg}}}$.
\begin{theorem}\label{thmclass}
Fix some integer $m \ge 1$. Let $\bbT^S_{cl}$ be the image of the map \[\bbT^S \rightarrow \prod_{\substack{U \in \cJ_{G, U^p} \\ k \geq 1}} \End_C(H^0(X_G^{U,*,\mathrm{ad}},\omega_U^{\otimes mk}\otimes\mathcal{I})).\] 
Then the surjective map $\bbT^S \rightarrow \bbT^S(R\Gamma_{X_G^V,c}\underline{\Lambda}_{U/V})$ factors over $\bbT^S_{cl}$.
\end{theorem}

\subsubsection{Comparison theorems}
We write $\bbV = \cO_C/(\pi^N) = \Lambda \otimes_\cO \cO_C$ and $A$ for the $\cO_C^a$-algebra $\bbV^a$, as in \S \ref{sec_almost_smooth}. We are going to compare various complexes in the derived categories $\mathbf{D}^+_\text{sm}(U_p,\Lambda)$ and $\mathbf{D}^+_\text{sm}(U_p,A)$ (as defined in \S \ref{sec_smooth_eqvt} and \S \ref{sec_almost_smooth}).

We first put ourselves in the situation of \S \ref{cc}. We set $X_0 = \overline{X}^U_G$ and let the tower $X_n$ be given by $X_n = \overline{X}^{U_{p,n}U^p}_G$, where $U_{p,n}$ runs over a cofinal system of compact open normal subgroups of $U_p$. Set $U_n = U_{p,n}U^p$. We set $X = \varprojlim_n X_n$. We have a functor (Definition \ref{def_derived_free_action})
\[R\Gamma_X: \mathbf{D}^+(\Sh(X_0,\Lambda)) \rightarrow \mathbf{D}^+_\text{sm}(U_p,\Lambda). \] 
Denote by $j_n$ the open embedding $j_n: X^{U_n}_G \hookrightarrow \overline{X}^{U_n}_G$, and let $K^{\mathrm{top}} = R\Gamma_X(j_{0,!}\Lambda) \in \mathbf{D}^+_\text{sm}(U_p,\Lambda)$.

Next, we work in a number of different settings where we can apply the formalism of \S \ref{ccnoprojlim}. We set $Y_n = X^{U_n}_G$ and $Y_n^\mathrm{alg} = X^{U_n,\mathrm{alg}}_G$. As above, we write $X_n^\mathrm{alg} = X^{\overline{U}_n,\mathrm{alg},*}_G/\overline{\bbQ}$ for the scheme theoretic image of the minimal compactification of $X^{U_n}_G$ in the minimal compactification of a Siegel modular variety of suitable level. Denote by $j'_n$ the open immersion $j'_n: Y^{\mathrm{alg}}_n\hookrightarrow X^{\mathrm{alg}}_n$. We also have associated adic spaces $Y_n^\mathrm{ad} = (Y_{n,C}^\mathrm{alg})^\mathrm{ad}$ and $X_n^\mathrm{ad} = (X_{n,C}^\mathrm{alg})^\mathrm{ad}$.

Note that for $n \ge m$ and $U_n' \subset U_m'$ subgroups of $\Sp_{2g}(\bbA_\bbQ^\infty)$ such that $U_n \subset U_n' \cap G(\bbA_{F^+}^\infty)$ and $U_m \subset U_m' \cap G(\bbA_{F^+}^\infty)$ we have a commutative diagram:

\[\xymatrix{X_G^{U_n,*} \ar[d]\ar[r] & X_{\Sp_{2g}}^{U_n', \ast}\ar[d]\\ X_G^{U_m,*} \ar[r]& X_{\Sp_{2g}}^{U_m', \ast}  }.\] If we choose $U_m'$ sufficiently small, the bottom horizontal map in the above diagram has scheme theoretic image  $X_m^\mathrm{alg}$, and we can then find $U_n'$ such that the top horizontal map has scheme theoretic image $X_n^\mathrm{alg}$. We can also view $X_m^\mathrm{alg}$ as the scheme theoretic image of the composite $X_G^{U_n,*} \rightarrow  X_G^{U_m,*} \rightarrow X_{\Sp_{2g}}^{U_m', \ast}$. Therefore we have an induced map between scheme theoretic images $X_n^\mathrm{alg} \rightarrow X_m^\mathrm{alg}$ \cite[Tag 01R9]{stacks-project}.

We let $X_n^*$ be the topological space given by the complex points of $X_{n,\bbC}^\mathrm{alg}$. We write $j'_n$ for the maps $Y_n \hookrightarrow X_n^*$ and $Y_n^\mathrm{ad}\hookrightarrow X_n^\mathrm{ad}$ induced by the algebraic $j'_n$. We also write $\pi_n$ for all of the projection maps $X_n^* \rightarrow X_0^*$, $Y_n \rightarrow Y_0$, etc.~ 

The formalism of \S \ref{ccnoprojlim} applies to the tower of spaces $(X^?_n)$, where $?$ is $*$, $\mathrm{alg}$, $\mathrm{ad}$ or nothing, and we obtain categories $S^?$. We denote the associated functors from $S^?$ to $\Mod_\text{sm}(U_p,\Lambda)$ (denoted $\widetilde{\Gamma}$ in Definition \ref{def_fibered_cc}) by $\widetilde{\Gamma}^?$, with right derived functors $R\widetilde{\Gamma}^?$.
\begin{lemma}\label{topcomp}
There is a natural isomorphism $K^{\mathrm{top}}\cong R\widetilde{\Gamma}^*(j'_!\Lambda)$.
\end{lemma}
\begin{proof}
This follows from the discussion in \S \ref{sec_compare_fibered_cc}: both complexes are naturally quasi-isomorphic to a complex $R\widetilde{\Gamma}_c(\Lambda)$ defined using the tower $(Y_n)_{n \geq 0}$.
\end{proof}
We denote by $K^{\mathrm{alg}}$ the complex $R\widetilde{\Gamma}^\mathrm{alg}(j'_!\Lambda)$ and denote by $K^{\mathrm{ad}}$ the complex $R\widetilde{\Gamma}^\mathrm{ad}(j'_!\Lambda)$.
\begin{lemma}\label{triplecomp}
There are natural isomorphisms: $K^{\mathrm{alg}} \cong K^{\mathrm{top}} \cong K^{\mathrm{ad}}$ in $\mathbf{D}^+_\text{sm}(U_p, \Lambda)$.
\end{lemma}
\begin{proof}
This follows from the discussion in \S \ref{sec_compare_fibered_cc}.
\end{proof}
We also consider the functor $\widetilde{\Gamma}_\bbV^\mathrm{ad}: S^\mathrm{ad}_\bbV \rightarrow \Modsm(U_p,\bbV)$ together with its right derived functor $R\widetilde{\Gamma}_\bbV^\mathrm{ad}$. Here $S^\mathrm{ad}_\bbV$ is defined in the same way as $S^\mathrm{ad}$, but with coefficients in $\bbV$. For $\mathcal{F} \in \Sh_{\text{\'et}}(X_n^\mathrm{ad},\bbV)$ there are natural isomorphisms \[\Gamma(X_n^\mathrm{ad}, \mathcal{F})\otimes_\Lambda \bbV \cong \Gamma(X_n^\mathrm{ad}, \mathcal{F}\otimes_\Lambda \bbV).\] 
Since direct limits commute with tensor products we also obtain a natural isomorphism \[\widetilde{\Gamma}^\mathrm{ad}(\mathcal{F})\otimes_\Lambda \bbV \cong \widetilde{\Gamma}_\bbV^\mathrm{ad}(\mathcal{F}\otimes_\Lambda \bbV).\]
\begin{lemma}\label{flatbasechange}
For $\mathcal{F} \in S^\mathrm{ad}$ the natural isomorphism $\widetilde{\Gamma}^\mathrm{ad}(\mathcal{F})\otimes_\Lambda \bbV \cong \widetilde{\Gamma}_\bbV^\mathrm{ad}(\mathcal{F}\otimes_\Lambda \bbV)$ extends to an isomorphism $R\widetilde{\Gamma}^\mathrm{ad}(\mathcal{F})\otimes_\Lambda \bbV \cong R\widetilde{\Gamma}_\bbV^\mathrm{ad}(\mathcal{F}\otimes_\Lambda \bbV)$.
\end{lemma}
\begin{proof}
It suffices to show that for an injective object $\mathcal{I} \in  S^\mathrm{ad}$ the higher derived functors $R^i\widetilde{\Gamma}^\mathrm{ad}(\mathcal{I}\otimes_\Lambda \bbV)$ vanish for $i > 0$. We have 
\[R^i\widetilde{\Gamma}^\mathrm{ad}(\mathcal{I}\otimes_\Lambda \bbV) = \varinjlim_n H^i(X_n^\mathrm{ad}, \mathcal{I}_n\otimes_\Lambda \bbV) = \varinjlim_n H^i(X_n^\mathrm{ad},\mathcal{I}_n)\otimes_\Lambda \bbV = 0,\] 
since $\mathcal{I}_n$ is injective. We are using the fact that we can compute cohomology of a sheaf of $\bbV$-modules after applying the forgetful functor to sheaves of $\Lambda$-modules and \cite[Rapport, 4.9.1]{SGA4demi}.
\end{proof}
\begin{lemma}\label{lem_ps_compare} The natural maps $\bbV \rightarrow \cO^+_{Y_n^{ad}}/(\pi^N)$ in $\Sh_{\text{\'et}}(Y_n^{ad},\Lambda)$ induce a map \[K^{ad}\otimes_\Lambda \bbV = R\widetilde{\Gamma}_\bbV^\mathrm{ad}(j'_!\bbV) \rightarrow R\widetilde{\Gamma}_\bbV^\mathrm{ad}((j'_{n,!}\cO^+_{Y_n^{ad}}/(\pi^N))_{n\ge 0})\] in $\mathbf{D}_\text{sm}(U_p,\bbV)$ which becomes an isomorphism in $\mathbf{D}_\text{sm}(U_p,A)$. 
\end{lemma}
\begin{proof}
This follows from Scholze's comparison theorem: the induced maps on cohomology are the natural maps
\[\varinjlim_n H^i_{\text{\'et},c}(Y_n^\mathrm{ad},\Lambda)\otimes_\Lambda A \rightarrow \varinjlim_n H^i_{\text{\'et},c}(Y_n^\mathrm{ad},\cO^{+}_{Y_n^{ad}}/(\pi^N))^a \]
which are isomorphisms by \cite[Theorem 3.13]{CDM} (and induction on $N$). 
\end{proof}

\subsubsection{Hecke operators}
We now define a Hecke action on the complex $R\Gamma_X(j_!\Lambda)$. We set $X_n^S = \overline{X}_G^{U_{p,n}U^p_S}$ and $X^S = \varprojlim_n X_n^S$, so that $X^S$ is a $U_p\times G^S$-space, where $U_p$ has the profinite topology and $G^S$ has the discrete topology. We have $X = (X^S)^{U^S}$. Therefore, by Proposition \ref{prop_top_heck_action_on_cohomology}, we have a natural map 
\[ T_G : \mathcal{H}(G^S,U^S) \rightarrow \mathrm{End}_{\mathbf{D}_{\text{sm}}(U_p,\Lambda)}(R\Gamma_X(j_!\Lambda))=\mathrm{End}_{\mathbf{D}_{\text{sm}}(U_p,\Lambda)}(K^\mathrm{top}).\]
For $g \in G^S$, we can describe explicitly the element $\theta(g) = T_G([U^SgU^S])$ in $\mathrm{End}_{\mathbf{D}_\text{sm}(U_p, \Lambda)}(R\Gamma_X(j_!\Lambda))$. Set $X' = (X^S)^{U^S \cap gU^Sg^{-1}}$, and consider the two maps $p_1,p_2: X' \rightarrow X$, where $p_1$ is the natural projection and $p_2$ is given by the (right) action of $g$ followed by the natural projection. The maps $p_i$ are finite \'{e}tale and we have a natural isomorphism $p_1^*j_!\Lambda \cong p_2^*j_!\Lambda$. We also set $X'_n = (X^S_n)^{U^S \cap gU^Sg^{-1}}$ and denote the two projection maps $X'_n \rightarrow X_n$ by $p_1,p_2$. Lemma \ref{lem_hecke_action_as_expected_top} implies:
\begin{lemma}
The endomorphism  $\theta(g)$ is given by the composition of natural maps: \[R\Gamma_X(j_!\Lambda) \overset{p_2^*}{\rightarrow} R\Gamma_{X'}(p_2^*j_!\Lambda) \cong R\Gamma_{X'}(p_1^*j_!\Lambda) \cong R\Gamma_X(p_{1,*}p_1^*j_!\Lambda) \rightarrow R\Gamma_X(j_!\Lambda)\] where the final map is the trace, defined by the adjunction $(p_{1,*}=p_{1,!},p_1^*=p_1^!)$.
\end{lemma}
The description of the above lemma can be translated into a description of the Hecke operators as a limit of Hecke operators at finite level. We apply the formalism of \S \ref{ccnoprojlim} to the towers of spaces $X_n$ and $X_n'$, with associated categories $S$ and $S'$. Let $j_!\Lambda \rightarrow \mathcal{I}^\bullet$ be an injective resolution in $S$.  Since we have an isomorphism  $p_2^*j_!\Lambda \cong p_1^*j_!\Lambda$, we have a map of complexes (unique up to homotopy) $p_2^*\mathcal{I}^\bullet \rightarrow p_1^*\mathcal{I}^\bullet$, and an induced map of complexes $p_2^*\mathcal{I}_n^\bullet \rightarrow p_1^*\mathcal{I}_n^\bullet$. Now for each $n$ we have maps (compatible as $n$ varies) \[\Gamma(X_n,\mathcal{I}_n^\bullet) \overset{p_2^*}{\rightarrow} \Gamma(X'_n,p_2^*\mathcal{I}_n^\bullet) \rightarrow \Gamma(X'_n,p_1^*\mathcal{I}_n^\bullet) \cong \Gamma(X_n,p_{1,*}p_1^*\mathcal{I}_n^\bullet) \rightarrow \Gamma(X_n,\mathcal{I}_n^\bullet).\] Taking the limit over $n$ gives the map $\theta(g)$, under the equivalence of Lemma \ref{lem_fibered_free_equiv}.
\begin{remark}
The following observation will be useful: to define the trace map $p_{1,*}p_1^*\mathcal{I}_n^\bullet \rightarrow \mathcal{I}_n^\bullet$ we only need to use the fact that $p_1$ is \'{e}tale over $Y_n$. Indeed, we have a map of sheaves on $Y_n$: $p_{1,*}p_1^*\Lambda \rightarrow \Lambda$, and applying $j_!$ gives a map of sheaves on $X_n$: $p_{1,*}p_1^*j_!\Lambda \rightarrow j_!\Lambda$. These maps lift (uniquely up to homotopy) to a map of complexes \[p_{1,*}p_1^*\mathcal{I}^\bullet \rightarrow \mathcal{I}^\bullet\] and this induces the desired compatible system of trace maps \[p_{1,*}p_1^*\mathcal{I}_n^\bullet \rightarrow \mathcal{I}_n^\bullet.\]
\end{remark}

Given the above remark, the description of $\theta(g)$ we have given applies immediately to define endomorphisms $\theta(g)$ of $K^{\mathrm{alg}}$ and $K^{\mathrm{ad}}$. Under the comparison isomorphisms of Lemma \ref{topcomp} and Lemma \ref{triplecomp}, we obtain the action of $\mathcal{H}(G^S,U^S)$ on $K^{\mathrm{ad}} \cong R\Gamma_X(j_!\Lambda)$. 

Similarly, this description gives endomorphisms $\theta(g)$ of $R\widetilde{\Gamma}^\mathrm{ad}((j_!\cO^+_{Y_n^{ad}}/(\pi^N))_{n\ge 0})$, and hence of $R\widetilde{\Gamma}^\mathrm{ad}((j_!\cO^+_{Y_n^{ad}}/(\pi^N))_{n\ge 0})^a$ such that the isomorphism \[\left(K^{ad}\otimes_\Lambda \bbV\right)^a \cong \left(R\widetilde{\Gamma}^\mathrm{ad}((j_!\cO^+_{Y_n^{ad}}/(\pi^N))_{n\ge 0})\right)^a\] in $\mathbf{D}_\text{sm}(U_p,A)$ is $\theta(g)$-equivariant.

\subsubsection{Comparison with Cech cohomology}
Finally, we are going to compute $R\widetilde{\Gamma}^\mathrm{ad}((j_!\cO^+_{Y_n^{ad}}/(\pi^N))_{n\ge 0})$ (and its Hecke action) using a Cech complex. Recall that by \cite[Theorem 4.1.1]{Sch14} (and the assumption that $U$ is small) there is a perfectoid space $X_\infty^{\mathrm{ad}}$ over $\mathrm{Spa}(C,\cO_C)$ with 
\[X_\infty^{\mathrm{ad}}\sim \varprojlim_n X_n^{\mathrm{ad}}.\] 

\begin{definition}\label{adcover}
A $U_p$-admissible cover of the perfectoid space $X_\infty^\mathrm{ad}$ is an open cover $\mathcal{V} = (V_i)_{i \in I}$ of  $X_\infty^{ad}$ by finitely many affinoid perfectoids (with affinoid perfectoid multiple intersections), such that:
\begin{enumerate}
\item There exists $n_0$ such that each $V_i$ is the inverse image of an affinoid open $V_{i,n_0}$ in $X_{n_0}^\mathrm{ad}$.
\item For $\gamma \in U_p$ and $i \in I$, $(V_i)\gamma\in \mathcal{V}$. 
\end{enumerate}
For $n \ge n_0$ we denote by $\mathcal{V}_n$ the affinoid cover of $X_n^\mathrm{ad}$ given by the inverse images of the $V_{i,n_0}$, and denote by $\widetilde{\mathcal{V}}_n$ the affinoid cover of $X_G^{U_n,*,\mathrm{ad}}$ obtained by pullback.
\end{definition}

We now recall some more of the results contained in \cite[Theorem 4.1.1]{Sch14}. Denote by $\mathscr{F}\ell$ the flag variety over $C$ which is associated with $(\Sp_{2g} , D_{\Sp_{2g}})$, i.e.\ the flag variety of totally isotropic subspaces of $C^{2g}$. There is a $(\Res^{F^+}_\bbQ G)(\bbQ_p)$-equivariant Hodge--Tate period map 
\[\pi_{HT} :X_\infty^\mathrm{ad}\rightarrow \mathscr{F}\ell.\]
For a subset $J \subset \{1,\ldots,2g\}$ of cardinality $g$ we denote by $s_J$ the corresponding Pl\"{u}cker coordinate on $\mathscr{F}\ell$ and denote by $\mathscr{F}\ell_J$ the open affinoid subspace of $\mathscr{F}\ell$ defined by $|s_{J'}|\le |s_J|$ for all $J'$. Now for $J \subset \{1,\ldots,2g\}$ of cardinality $g$ we denote by $X_{\infty,J}^\mathrm{ad}$ the preimage $\pi_{HT}^{-1}(\mathscr{F}\ell_J)$ in $X_\infty^\mathrm{ad}$. By \cite[Theorem 4.1.1]{Sch14} the subsets $X_{\infty,J}^\mathrm{ad}$ provide an affinoid perfectoid cover of $X_\infty^\mathrm{ad}$. Moreover, they satisfy the first condition in Definition \ref{adcover}. This means that there are only finitely many $U_p$ translates of each $X_{\infty,J}^\mathrm{ad}$ (they are each stabilised by a compact open subgroup of $U_p$), and so by adjoining these translates we obtain a $U_p$-admissible cover $(X_{\infty,J}\cdot\gamma)_{\gamma,J}$ of $X_\infty^\mathrm{ad}$. This is the only $U_p$-admissible cover we will use in practice.

\begin{lemma}\label{coverlemma} For $\gamma \in U_p$ and $J \subset \{1,\ldots,2g\}$ of cardinality $g$, the open affinoid subspace $\mathscr{F}\ell_J\cdot\gamma$ of $\mathscr{F}\ell$ is defined by the inequalities $|(\gamma')^{-1} s_{J'}|\le |\gamma^{-1} s_J|$ for all $\gamma' \in U_p$ and $J' \subset \{1,\ldots,2g\}$ of cardinality $g$.
\end{lemma}
\begin{proof}
We have $x \in \mathscr{F}\ell_J\cdot \gamma$ if and only if $|\gamma^{-1}s_J(x)| = |s_J(x\gamma^{-1})| \ge |s_{J'}(x\gamma^{-1})|$ for all $J'$. It suffices to show that if $x \in \mathscr{F}\ell_J$ then $|s_J(x)| \ge |s_{J'}(x\gamma^{-1})|$ for all $J'$ and $\gamma$ (i.e.~we set $\gamma = 1$ and $\gamma' = \gamma$).

The action of $U_p$ on the coordinates $s_J$ is given by the action of elements of $\GL_{2g}(\bbZ_p)$ (in fact they are elements of $\Sp_{2g}(\bbZ_p)$) on basis elements of $\bigwedge^g\mathrm{Std}$, where $\mathrm{Std}$ is the standard $2g$-dimensional representation of $\GL_{2g}/\bbZ_p$. In particular, we have $\gamma^{-1} s_{J'} = \sum_{J''} a_{J''}s_{J''}$ with $a_{J''} \in \bbZ_p$. So if $x \in \mathscr{F}\ell_J$ we have \[|s_{J'}(x\gamma^{-1})| = |\sum_{J''} a_{J''}s_{J''}(x)| \le \max_{J''}(|s_{J''}(x)|) = |s_J(x)|.\]
\end{proof}

We write $I$ for the (finite) index set of the cover $\mathcal{V}$ of $X_\infty^\mathrm{ad}$, and for $i \in I$ write $s_i$ for the section $\gamma_i^{-1} s_{J_i}$ of $\omega$ on $X_\infty^\mathrm{ad}$. For $n \ge n_0$ we write $V_{i,n} \in \widetilde{\mathcal{V}}_n$ for the open affinoid in $X_G^{U_n,*,\mathrm{ad}}$ obtained from $V_i$.

\begin{lemma}
For $n \ge n_0$ sufficiently large, there exist sections $s_j^{(i)}$ (for all pairs $i,j \in I$) of $\omega$ on $V_{i,n} \in \widetilde{\mathcal{V}}_n$ such that \begin{itemize}
\item For each $i, j \in I$, \[\left|\frac{s_j - s_j^{(i)}}{s_i^{(i)}}\right| \le |\pi^N|\] on $V_{i,n}$.
\item For each $i, j \in I$, the section $s^{(i)}_i$ is invertible and $s^{(i)}_j / s^{(i)}_i \in H^0(V_{i,n}, \cO^+)$.
\item For each $i, j \in I$, the subset $V_{ij,n} = V_{i,n} \cap V_{j,n} \subset V_{i,n}$ is defined by the equation $| s^{(i)}_j / s^{(i)}_i | = 1$.
\item For each $i, j, k \in I$, the inequality 
\[ | s^{(i)}_k / s^{(i)}_j - s^{(j)}_k / s^{(i)}_j | \leq |\pi^N| \]
holds in $V_{ij,n}$.
\end{itemize}
In particular, the cover $\widetilde{\mathcal{V}}_n$ and the sections $s_j^{(i)}$ satisfy the assumptions of \cite[Lemma 2.1.1]{Sch14}.
\end{lemma}
\begin{proof}
We first choose $n$ large enough so that there are sections $s_j^{(i)}$ approximating the $s_j$: more precisely we suppose that \[\left|\frac{s_j - s_j^{(i)}}{s_i}\right| \le |\pi^N|\] on $V_i$. This is possible because the map \[\varinjlim_n H^0(V_{i,n},\omega) \rightarrow H^0(V_i,\omega) \] has dense image \cite[Theorem 4.1.1, (i)]{Sch14}. The itemized conditions now follow from Lemma \ref{coverlemma} and the fact that the sections $s_j^{(i)}$ approximate $s_j$.
\end{proof}
\begin{remark}
To illustrate the process of adjoining $U_p$-translates to an affinoid cover we discuss the case of $\bbP^1_{\bbQ_p}$ with its right action of $\GL_2(\bbZ_p)$. We begin with the affinoid cover given by $\{|z|\le 1\}$ and $\{|z|\ge 1\}$. These are the complements of the (open) residue discs $\mathrm{red}^{-1}(0)$ and $\mathrm{red}^{-1}(\infty)$ where $\mathrm{red}$ is the reduction map to $\bbP^1_{\bbF_p}$. These affinoids are stable under the action of $1+pM_2(\bbZ_p)$, and the translates by $\GL_2(\bbF_p)$ are the $p+1$ affinoids given by the complements of the residue discs $\mathrm{red}^{-1}(x)$ for $x \in \bbP^1_{\bbF_p}$. Since the action of $\GL_2(\bbZ_p)$ on $\bbP^1_{\bbQ_p}$ extends to an action on $\bbP^1_{\bbZ_p}$, the formal model of $\bbP^1_{\bbQ_p}$ obtained from this cover by $p+1$ affinoids is again the formal completion of $\bbP^1_{\bbZ_p}$ along the special fibre.

If we first apply \cite[Lemma 2.1.1]{Sch14} to the cover $X_{n,J}^\mathrm{ad} = \Spa(R_{n,J},R_{n,J}^+)$ of $X_n^\mathrm{ad}$ (for $n\ge n_0$) we obtain a formal model $\mathfrak{X}_n$ for $X_n^\mathrm{ad}$ over $\cO_C$ with an affine cover by $\mathrm{Spf}(R_{n,J}^+)$. For $\gamma \in U_p$ and $J' \subset \{1,\ldots,2g\}$ of cardinality $g$, the intersection $X_{n,J}^\mathrm{ad}\cap X_{n,J'}^\mathrm{ad}\cdot\gamma \subset X_{n,J}^\mathrm{ad}$ is defined by \[\left|\frac{\gamma^{-1}s_{J'}}{s_J}\right|=1.\] This in turn defines an affine open formal subscheme \[\mathrm{Spf}\left(R_{n,J}^+\left\langle \frac{s_J}{\gamma^{-1}s_{J'}} \right\rangle \right)\subset\mathrm{Spf}(R_{n,J}^+),\] and so we see that the cover $(X_{n,J}^\mathrm{ad}\cdot\gamma)_{J,\gamma}$ of $X_n^\mathrm{ad}$ is the generic fibre of a cover of $\mathfrak{X}_n$ by formal affine opens. Therefore the formal model for $X_n^\mathrm{ad}$ given by applying \cite[Lemma 2.1.1]{Sch14} to the cover $(X_{n,J}^\mathrm{ad}\cdot\gamma)_{J,\gamma}$ of $X_n^\mathrm{ad}$ is the same as the formal model coming from the cover $(X_{n,J}^\mathrm{ad})_{J}$. The same remarks apply to the formal model we obtain for $X_G^{U_n,*,\mathrm{ad}}$.
\end{remark}

Let $\mathcal{F} \in S^\mathrm{ad}$ be a system of sheaves of $\bbV$-modules, and let $\mathcal{V}$ be a $U_p$-admissible cover of $X_\infty^\mathrm{ad}$. We define a Cech complex $\widetilde{C}^\bullet(\mathcal{V},\mathcal{F})$ with entries in $\Modsm(U_p,\bbV)$ by \[\widetilde{C}^p(\mathcal{V},\mathcal{F}) = \varinjlim_{n \ge n_0} {C}^p(\mathcal{V}_n,\mathcal{F}_n)\] where  ${C}^\bullet(\mathcal{V}_n,\mathcal{F}_n)$ is the usual Cech complex for the sheaf $\cF_n$ on $X_n^\mathrm{ad}$ with respect to the cover $\cV_n$, endowed with its natural $U_p/U_{p,n}$-action (for example, $g \in U_p$ maps a section in $\mathcal{F}_n(V_{n,i})$ to a section in $\mathcal{F}_n(V_{n,i}g^{-1})$).
\begin{lemma}\label{Cech}
Let $\mathcal{F} \in S^\mathrm{ad}$ and let $\mathcal{V}$ be a $U_p$-admissible cover of $X_\infty^\mathrm{ad}$. Then there is a natural map in $\mathbf{D}_\text{sm}(U_p,\bbV)$: \[\widetilde{C}^\bullet(\mathcal{V},\mathcal{F}) \rightarrow R\widetilde{\Gamma}_\bbV^\mathrm{ad}(\mathcal{F}).\]
If $\mathcal{F} = (j_!\cO_{Y_n^\mathrm{ad}}^+/(\pi^N))_{n \ge 0}$ then this map becomes an isomorphism in $\mathbf{D}_\text{sm}(U_p, A)$.
\end{lemma}
\begin{proof}
The first part is \cite[Tag 03AX]{stacks-project}: let $\mathcal{F}\rightarrow \mathcal{I}^\bullet$ be an injective resolution. Define a double complex \[A^{p,q} = \widetilde{C}^p(\mathcal{V},\mathcal{I}^q),\] and denote by $sA^\bullet$ the total complex. The natural maps $\widetilde{\Gamma}^\mathrm{ad}(\mathcal{I}^q) \rightarrow A^{0,q}$ induce a quasi-isomorphism $\alpha:\widetilde{\Gamma}^\mathrm{ad}(\mathcal{I}^\bullet) \rightarrow sA^\bullet$. The maps $\widetilde{C}^p(\mathcal{V},\mathcal{F}) \rightarrow A^{p,0}$ induce a map of complexes \[\widetilde{C}^\bullet(\mathcal{V},\mathcal{F}) \rightarrow sA^\bullet.\] Composing this map with the inverse of $\alpha$ gives the desired map in the derived category. 

The second statement is proved as in \cite[Theorem 4.2.1]{Sch14}: we just need to check that the natural map described above induces an almost isomorphism on cohomology groups, and this can be checked after forgetting the action of $U_p$.
\end{proof}

\subsubsection{The end of the proof}
\begin{lemma}\label{lem_completed_almost_factors}
The map $\bbT^S \to \End_{\mathbf{D}_{\text{sm}}(U, A)}(R\widetilde{\Gamma}^\mathrm{ad}((j_!\cO^+_{Y_n^{ad}}/(\pi^N))_{n \ge 0})^a)$ factors through $\bbT^S_{cl}$.
\end{lemma}
\begin{proof}
This is implied by Lemma \ref{Cech}, following the proof of \cite[Theorem 4.3.1]{Sch14}. Indeed, $\bbT^S$ acts on each term of the Cech complex $\widetilde{C}^\bullet(\mathcal{V},(j_!\cO^+_{Y_n^{ad}}/(\pi^N))_{n \ge 0})$, so it suffices to show that it acts via $\bbT^S_{cl}$ on each term $H^0(\cap_{i \in J}V_i,j_!\cO^+/(\pi^N))^a$. We set $\mathcal{V}$ to be the $U_p$-admissible cover $(X_{\infty,J}\cdot \gamma)_{J,\gamma}$, and proceed exactly as in \emph{loc.~cit.}, using the sections $\gamma^{-1}s_{J}$ of the line bundle $\omega$.
\end{proof}
\begin{corollary}\label{almostfactors}
The map $\bbT^S \to \End_{\mathbf{D}(U/V, A)}(R\Gamma_{X_G^V,c}\underline{\Lambda}_{U/V}\otimes_{\Lambda}A)$ factors through $\bbT^S_{cl}$.
\end{corollary}
\begin{proof}
By Lemma \ref{lem_ps_compare} we have an isomorphism $R\Gamma(V_p,R\widetilde{\Gamma}^\mathrm{ad}((j_!\cO^+_{Y_n^{ad}}/(\pi^N))_{n \ge 0})^a) \cong R\Gamma(V_p,(K^\mathrm{ad}\otimes_\Lambda \bbV)^a)$. By Lemma \ref{triplecomp} this is isomorphic to $R\Gamma(V_p,(K^\mathrm{top}\otimes_\Lambda \bbV)^a)$, and by Lemma \ref{lem_commute_a_and_invts} this is isomorphic to $R\Gamma(V_p,K^\mathrm{top}\otimes_\Lambda \bbV)^a$. Finally, by Lemma \ref{lem_flat_coeff_change} this is isomorphic to $R\Gamma(V_p,K^\mathrm{top})\otimes_{\Lambda} A$, and the Corollary follows from the Lemma \ref{lem_completed_almost_factors} and Lemma \ref{lem_invts_quo}.
\end{proof}

\begin{lemma}\label{lem_fg_complex}
$R\Gamma_{X_G^V,c}\underline{\Lambda}_{U/V}$ is isomorphic to a bounded complex of finitely generated $\Lambda[U/V]$-modules.
\end{lemma}
\begin{proof}
We know that the $R^i\Gamma_{X_G^V,c}\underline{\Lambda}_{U/V}$ are all finitely generated and non-zero for only finitely many $i$. So our statement follows from \cite[\S5, Lemma 1]{Mumford}.
\end{proof}
We can now finish the proof of Theorem \ref{thmclass}. We restate the result:
\begin{theorem*}
The map $\bbT^S \to \End_{\mathbf{D}(U/V, \Lambda)}(R\Gamma_{X_G^V,c}\underline{\Lambda}_{U/V})$ factors through $\bbT^S_{cl}$.
\end{theorem*}
\begin{proof}
We apply Corollary \ref{almostfactors}. With this in mind, it suffices to show that the natural map \[\End_{\mathbf{D}(U/V, \Lambda)}(R\Gamma_{X_G^V,c}\underline{\Lambda}_{U/V}) \rightarrow \End_{\mathbf{D}(U/V, A)}(R\Gamma_{X_G^V,c}\underline{\Lambda}_{U/V}\otimes_{\Lambda}A)\] is injective. Combining Lemma \ref{lem_fg_complex}, Lemma \ref{lem_derived_homtensor} and Lemma \ref{lem_almost_hom} we see that this map is a map of $\Lambda$-modules of the form \[M \rightarrow N ,\] where $M$ is a $\Lambda$-module, $N$ is a $\bbV$-module, and the induced map $M\otimes_\Lambda\bbV\rightarrow N$ is an almost isomorphism. In particular, if we write $K$ for the kernel of $M \rightarrow N$, then $K\otimes_\Lambda\bbV$ is almost zero. This implies that $K$ is zero, as desired. Indeed, if $K\otimes_\Lambda\bbV$ is almost zero then $K/\pi K \otimes_\Lambda\bbV$ is simultaneously almost zero and a free $\bbV/\pi\bbV$-module, and is therefore zero. So $K/\pi K = 0$ which implies $K = 0$.

\end{proof}
\bibliographystyle{alpha}

\end{document}